\documentclass{amsart}

\usepackage[utf8]{inputenc}
\usepackage[T1]{fontenc}
\usepackage[english]{babel}

\usepackage{amsfonts}
\usepackage{amsmath}
\usepackage{amssymb}
\usepackage{amstext}
\usepackage{amsthm}

\usepackage{geometry}
\usepackage{fullpage}

\usepackage[shortlabels]{enumitem}
\usepackage{graphicx}
\usepackage{hyperref}
\usepackage{mathrsfs}
\usepackage{mathtools}

\usepackage{tikz}
\usepackage{tikz-cd}
\usetikzlibrary{arrows,backgrounds,chains,decorations.pathreplacing,patterns,shapes, calc, positioning}

\allowdisplaybreaks

\newcommand{\wt}{\operatorname{wt}}
\newcommand{\DD}{\mathsf{D}}
\newcommand{\RR}{\mathsf{R}}

\newcommand{\area}{\operatorname{area}}
\newcommand{\dinv}{\operatorname{dinv}}
\newcommand{\PF}{\operatorname{PF}}
\newcommand{\touch}{\operatorname{touch}}
\newtheorem{theorem}{Theorem}[section]
\newtheorem{lemma}[theorem]{Lemma}
\newtheorem{proposition}[theorem]{Proposition}
\newtheorem{corollary}[theorem]{Corollary}
\newtheorem{conjecture}[theorem]{Conjecture}

\theoremstyle{definition}
\newtheorem{definition}[theorem]{Definition}
\newtheorem{example}[theorem]{Example}

\theoremstyle{remark}

\theoremstyle{plain}

\title{A proof of the Theta conjecture}

\author{Haijun Bai}
\address{Center for Combinatorics, LPMC, Nankai University, Tianjin 300071, P. R. China}
\email{haojunbai@mail.nankai.edu.cn}
\thanks{}

\author{Feng Gu}
\address{Center for Combinatorics, LPMC, Nankai University, Tianjin 300071, P. R. China}
\email{fgu@mail.nankai.edu.cn}
\thanks{}

\author{Zhuowei Lin}
\address{Center for Combinatorics, LPMC, Nankai University, Tianjin 300071, P. R. China}
\email{zwlin0825@163.com}
\thanks{}

\author{Dun Qiu}
\address{Center for Combinatorics, LPMC, Nankai University, Tianjin 300071, P. R. China}
\email{qiudun@nankai.edu.cn}
\thanks{Dun Qiu is supported in part by the National Natural Science Foundation of China (No. 12671390),
the Natural Science Foundation of Tianjin (No. 24JCZDJC01390), and the Fundamental Research Funds for the Central Universities (No. 105-63263093).}

\author{Minhao Zhang}
\address{Center for Combinatorics, LPMC, Nankai University, Tianjin 300071, P. R. China}
\email{carpentercheungmath@gmail.com}
\thanks{}

\begin{document}

\subjclass[2020]{Primary: 05E05; Secondary: 05E10, 33D52, 16T30}

\begin{abstract}
We prove the touching refinement of the Theta conjecture of
Interdonato and Iraci, giving a formula for
$\Theta_{e_\ell}\Theta_{e_k}\nabla E_{d,j}$
as a weighted sum over parking functions with decorated rises
and valleys and touch number $j$.
Summing over $j$ gives the Theta conjecture.
As corollaries, we prove the valley Delta conjecture of
Haglund, Remmel, and Wilson and its generalized version
proposed by Qiu and Wilson.
Our proof compares algebraic and combinatorial recurrences
using an extension of the plethystic substitution
$X\mapsto(q-1)X$ to quasisymmetric functions.
The algebraic recurrence follows from an endomorphism
constructed by D'Adderio, Interdonato, Iraci and Pagaria
on the Dyck path algebra of Carlsson and Mellit.
The combinatorial recurrence is obtained by introducing
signed labels and deleting the main diagonal.
The two recurrences have the same initial conditions,
which identifies the corresponding generating functions.
\end{abstract}

\maketitle

\section{Introduction}
The nabla operator $\nabla$, introduced by Bergeron and
Garsia~\cite{BG99}, acts diagonally on the modified Macdonald
polynomials.  Garsia and Haiman~\cite{GH96} conjectured that
$\nabla e_d$ is the bigraded Frobenius characteristic of the
diagonal coinvariant ring.  Haiman~\cite{Haiman02} proved this
conjecture using the geometry of the Hilbert scheme.

Haglund, Haiman, Loehr, Remmel and
Ulyanov~\cite{HHLRU05} proposed a combinatorial expression for
$\nabla e_d$ in terms of parking functions, namely the shuffle
conjecture.  Carlsson and Mellit~\cite{CM18} proved it, obtaining
\begin{equation}\label{eqRefinedS13Shuffle}
  \nabla e_d
  =
  \sum_{\pi\in\PF_d}
  q^{\dinv(\pi)}t^{\area(\pi)}x^\pi,
\end{equation}
where $\PF_d$ denotes the set of parking functions of size $d$.

Haglund, Remmel and Wilson~\cite{HRW18}
proposed the Delta conjecture, giving 
two conjectural combinatorial formulas for $\Delta'_{e_{d-1}}e_{d+k}$:
one using decorated rises and the other using decorated
contractible valleys.

We write $\PF_{d+k+\ell}^{\ast k,\bullet\ell}$ for the set of
parking functions with $k$ decorated rises and $\ell$ decorated
contractible valleys.  Their index sets are denoted by
$\RR(\pi)$ and $\DD(\pi)$, respectively.
These decorations and the statistics $\dinv$ and $\area$
are recalled in Subsection~\ref{subsecCombinatorialObjects}.

\begin{conjecture}[Haglund et al.~\cite{HRW18}]
\label{conjIntroDelta}
For integers $d\geq1$ and $k\geq0$, we have
\[
  \Delta'_{e_{d-1}}e_{d+k}
  =
  \sum_{\pi\in\PF_{d+k}^{\ast k,\bullet0}}
  q^{\dinv(\pi)}t^{\area(\pi)}x^\pi
  =
  \sum_{\pi\in\PF_{d+k}^{\ast0,\bullet k}}
  q^{\dinv(\pi)}t^{\area(\pi)}x^\pi.
\]
\end{conjecture}

The two expressions are the rise version and the valley version.
At $k=0$ both reduce to~\eqref{eqRefinedS13Shuffle},
since $\Delta'_{e_{d-1}}e_d=\nabla e_d$.

The rise version has been extended in two directions.
D'Adderio, Iraci and Vanden Wyngaerd~\cite{DIVW21} formulated a
compositional refinement of the rise version.
D'Adderio and Mellit~\cite{proofDelta} proved this refinement, thereby
proving the rise version of Conjecture~\ref{conjIntroDelta}.
Haglund et al.~\cite{HRW18} proposed the extended
Delta conjecture for $\Delta_{h_m}\Delta'_{e_{d-1}}e_{d+k}$,
using partial parking functions with $m$ zero labels,
$d+k$ positive labels, and $k$ decorated rises.
Blasiak, Haiman, Morse, Pun and Seelinger~\cite{BHMPS23}
proved this conjecture.

For the valley version, Qiu and
Wilson~\cite{QW20}
proposed the generalized valley Delta conjecture and proved
its specializations at $q=0$ and $t=0$.
Iraci and Vanden Wyngaerd~\cite{IW22}
showed that the refinement of the valley Delta conjecture
by the number of touches implies its generalized version.

The Theta operators $\Theta_f$ were introduced by D'Adderio et al.~\cite{DIVW21}.
D'Adderio et al.~\cite{DIVW21} conjectured a
formula for $\Theta_{e_\ell}\Theta_{e_k}\nabla e_d$ at $q=1$ in terms of parking
functions with both kinds of decorations, and
D'Adderio, Interdonato, Iraci and Pagaria~\cite{DIIP} proved this formula at $q=1$.
Interdonato and Iraci~\cite{Thetaconj} then defined a $\dinv$
statistic on parking functions carrying both kinds of decorations
and formulated the following conjecture.

\begin{conjecture}[Interdonato and Iraci~\cite{Thetaconj}]
\label{conjIntroTheta}
For integers $d\geq1$ and $k,\ell\geq0$, we have
\[
  \Theta_{e_\ell}\Theta_{e_k}\nabla e_d
  =
  \sum_{\pi\in\PF_{d+k+\ell}^{\ast k,\bullet\ell}}
  q^{\dinv(\pi)}t^{\area(\pi)}x^\pi.
\]
\end{conjecture}

Setting $\ell=0$ or $k=0$ gives the rise or valley version
of Conjecture~\ref{conjIntroDelta}, respectively.
Interdonato and Iraci~\cite{Thetaconj}
also proposed a refinement by $\touch(\pi)$, the number of
undecorated north steps starting on the main diagonal.

\begin{conjecture}[Interdonato and Iraci~\cite{Thetaconj}]
\label{Thetaconjecture}
For integers $d\geq j\geq1$ and $k,\ell\geq0$, we have
\[
  \Theta_{e_\ell}\Theta_{e_k}\nabla E_{d,j}
  =
  \sum_{\substack{
    \pi\in\PF_{d+k+\ell}^{\ast k,\bullet\ell}\\
    \touch(\pi)=j
  }}
  q^{\dinv(\pi)}t^{\area(\pi)}x^\pi.
\]
\end{conjecture}

By~\eqref{eqRefinedS13ESum}, summing
Conjecture~\ref{Thetaconjecture} over $j$ gives
Conjecture~\ref{conjIntroTheta}.
Our main result proves the touching refinement.

\begin{theorem}[{= Theorem~\ref{thmTouchingTheta}}]
\label{thmIntroMain}
For integers $d\geq j\geq1$ and $k,\ell\geq0$, we have
\[
  \Theta_{e_\ell}\Theta_{e_k}\nabla E_{d,j}
  =
  \sum_{\substack{
    \pi\in\PF_{d+k+\ell}^{\ast k,\bullet\ell}\\
    \touch(\pi)=j
  }}
  q^{\dinv(\pi)}t^{\area(\pi)}x^\pi.
\]
\end{theorem}

Summing Theorem~\ref{thmIntroMain} over $j$ gives the Theta conjecture, and then setting
$k=0$ gives the valley Delta conjecture.
The case $k=0$ of Theorem~\ref{thmIntroMain} also gives
the touching refinement required by Iraci and
Vanden Wyngaerd~\cite{IW22}, yielding the
generalized valley Delta conjecture.

\begin{corollary}\label{corThetaDeltaConsequences}
\leavevmode
\begin{enumerate}
  \item[\textup{(1)}]
  Conjecture \ref{conjIntroTheta} holds.

  \item[\textup{(2)}]
  The valley Delta conjecture of Haglund et al.~\cite{HRW18} holds.
  That is, for $d\geq1$ and $\ell\geq0$,
  \[
    \Delta'_{e_{d-1}}e_{d+\ell}
    =
    \sum_{\pi\in\PF_{d+\ell}^{\ast0,\bullet\ell}}
    q^{\dinv(\pi)}t^{\area(\pi)}x^\pi.
  \]

  \item[\textup{(3)}]
  The generalized valley Delta conjecture of Qiu and
  Wilson~\cite{QW20} holds.
  That is, for $d\geq1$ and $m,\ell\geq0$,
  \[
    \Delta_{h_m}\Delta'_{e_{d-1}}e_{d+\ell}
    =
    \sum_{\pi}q^{\dinv(\pi)}t^{\area(\pi)}x^\pi,
  \]
  where the sum ranges over partial parking functions with
  $m$ zero labels, $d+\ell$ positive labels, no decorated
  rises, and $\ell$ decorated contractible valleys.
  Here a partial parking function has nonnegative integer labels,
  strictly increasing from bottom to top in each column, with a positive label on the
  first north step; we set $x_0=1$.  Contractibility, area, and
  $\dinv$ use the same rules as above, now allowing label $0$,
  as in \cite[Definitions~3.3--3.10]{IW22}.
\end{enumerate}
\end{corollary}

We outline the proof of Theorem~\ref{thmIntroMain}.
Let $u,v,r,s$ be formal variables recording, respectively,
the number of north steps carrying neither decoration,
the touch number, the number of decorated rises, and the
number of decorated valleys.  Define
\begin{equation}\label{eqThetaalg}
  L[X;u,v,r,s]
  :=
  \sum_{j\geq1}\sum_{d\geq j}\sum_{k,\ell\geq0}
  u^dv^jr^ks^\ell
  \bigl(\Theta_{e_\ell}\Theta_{e_k}\nabla E_{d,j}\bigr)[X].
\end{equation}
The corresponding generating function for parking functions is
\begin{equation}\label{eqRefinedS13RIntroduction}
  R(X;u,v,r,s)
  :=
  1+
  \sum_{j\geq1}\sum_{d\geq j}\sum_{k,\ell\geq0}
  u^dv^jr^ks^\ell
  \sum_{\substack{
    \pi\in\PF_{d+k+\ell}^{\ast k,\bullet\ell}\\
    \touch(\pi)=j
  }}
  q^{\dinv(\pi)}t^{\area(\pi)}x^\pi,
\end{equation}
where the constant term accounts for the empty parking function.
Thus Theorem~\ref{thmIntroMain} is equivalent to $R=1+L$.

We use the operator $U_q$ defined in~\eqref{eqRefinedS13UqDefinition} 
from Haglund, Haiman and Loehr~\cite{HHL05}.
By Proposition~\ref{thmUqProperties}\textup{(1)}, it extends the plethystic substitution
$g[X]\mapsto g[(q-1)X]$ to quasisymmetric functions.
By Proposition~\ref{propUqSignedPaths}, $U_qR$ is the generating
function for signed parking functions, with labels $i^-$
and $i^+$ weighted by $-x_i$ and $qx_i$, respectively.

On the algebraic side, we use the endomorphism
$\boldsymbol{\theta}_u$ of D'Adderio et al.~\cite{DIIP} and its action on the
polynomials, recalled in
Theorem~\ref{DIIPthm}, to obtain the following expansion.

\begin{theorem}[{= Theorem~\ref{thmAlgMainRecursion}}]
\label{thmIntroAlgebraicExpansion}
For a finite alphabet $X_n$ and $j\geq1$, we have
\begin{equation}\label{eqintroalg}
  [v^j]L[(q-1)X_n;u,v,r,s]
  =
  (-u)^j(q;q)_j
  \sum_{\substack{I\subseteq X_n\\|I|=j}}
  \left(
    \prod_{x\in I}x
    \prod_{\substack{x\in I\\h\in X_n\setminus I}}
    \frac{x-qh}{x-h}
  \right)
  \mathcal L_I(X_n;ut,r,s).
\end{equation}
The functions $\mathcal L_I$ are uniquely determined by
the recurrence~\eqref{eqAlgRecursion} and the initial condition
$\mathcal L_\varnothing=1$, as formal power series in $u,r,s$
with coefficients rational in the letters of $X_n$.
\end{theorem}

On the combinatorial side, we delete the north steps starting
on the main diagonal and record the data needed to reconstruct
the original path.
The changes in $\dinv$ and $\area$ are given by
Theorem~\ref{thmMainDiagonalStatistics}.
Summing over the possible labels and decorations of the
deleted steps leads to a corresponding expansion for $U_qR$.

\begin{theorem}[{= Theorem~\ref{thmCombinatorialMarkedExpansion}}]
\label{thmIntroCombinatorialExpansion}
For a finite alphabet $X_n$ and $0\leq j\leq n$, we have
\[
  [v^j](U_qR)(X_n;u,v,r,s)
  =
  (-u)^j(q;q)_j
  \sum_{\substack{I\subseteq X_n\\|I|=j}}
  \left(
    \prod_{x\in I}x
    \prod_{\substack{x\in I\\h\in X_n\setminus I}}
    \frac{x-qh}{x-h}
  \right)
  \mathcal R_I(X_n;ut,r,s),
\]
where $\mathcal R_I$ is defined in~\eqref{eqRefinedS7RDefinition} and
$\mathcal R_\varnothing=1$.
For $j>n$, the coefficient is zero.
\end{theorem}

Theorem~\ref{thmMainCombinatorialRecursion} shows that
$\mathcal R_I$ satisfies~\eqref{eqAlgRecursion} with
$\mathcal L$ replaced by $\mathcal R$.
The common recurrence and initial condition therefore give
$\mathcal L_I=\mathcal R_I$.
Comparing Theorems~\ref{thmIntroAlgebraicExpansion}
and~\ref{thmIntroCombinatorialExpansion} yields
\[
  (U_qR)(X_n;u,v,r,s)
  =
  1+L[(q-1)X_n;u,v,r,s].
\]
Since this holds for every finite alphabet and $U_q$ is
invertible by Proposition~\ref{thmUqProperties}\textup{(2)},
we obtain $R=1+L$, proving Theorem~\ref{thmIntroMain}.

The paper is organized as follows.
Section~\ref{secPreliminaries} recalls symmetric functions
and introduces parking functions and their signed
versions.
Section~\ref{secThetaAlgebra} proves
Theorem~\ref{thmIntroAlgebraicExpansion} using the Theta
operators and the Dyck path algebra.
Section~\ref{secDiagonalDeletion} describes deletion of the
main diagonal and calculates its effect on the statistics.
Section~\ref{secInsertionsWithoutRise} evaluates the sums
arising from reconstruction without rise decorations.
Section~\ref{secInsertionsWithRise} determines how these sums
change when rise decorations are allowed.
Finally, Section~\ref{secCombinatorialRecurrence} proves
Theorem~\ref{thmIntroCombinatorialExpansion} and the recurrence
for $\mathcal R_I$, completing the proofs of
Theorem~\ref{thmIntroMain} and
Corollary~\ref{corThetaDeltaConsequences}.

A Lean 4 formalization of
Theorem \ref{thmIntroMain} is available at \href{https://github.com/qiudun123/theta-conjecture}{https://github.com/qiudun123/theta-conjecture}.

\section{Preliminaries}\label{secPreliminaries}
In Subsection~\ref{subsecSymmetricFunctions}, we recall symmetric functions
and plethystic notation needed below; see Macdonald~\cite{Macdonald95} for more details.
Subsection~\ref{subsecCombinatorialObjects} recalls the parking functions
and their related statistics, following the notation of Interdonato and Iraci~\cite{Thetaconj}.

\subsection{Symmetric functions}\label{subsecSymmetricFunctions}

A \emph{partition} $\mu=(\mu_1,\ldots,\mu_\ell)$ is a weakly decreasing
sequence of positive integers.  We write $|\mu|=\mu_1+\cdots+\mu_\ell$
for its size, $\ell(\mu)=\ell$ for its length, and $\mu\vdash n$
when $|\mu|=n$.  We identify $\mu$ with its Ferrers diagram in French notation:
rows are numbered from bottom to top, and a cell $(i,j)$ lies in column $i$ and row $j$, with $1\leq i\leq\mu_j$.

Let $\Lambda=\bigoplus_{n\geq0}\Lambda^{(n)}$ be the graded
algebra of symmetric functions over $\mathbb{Q}(q,t)$.
We use the standard bases of $\Lambda$: $e_\mu$, $h_\mu$, and $p_\mu$ for the elementary,
complete homogeneous, and power-sum symmetric functions, respectively.
The involution $\omega$ is determined by $\omega(e_m)=h_m$.

Throughout, $X=x_1+x_2+\cdots$ denotes the infinite alphabet,
while $X_n=x_1+\cdots+x_n$ denotes the finite alphabet of $n$
distinct variables, with $X_0=\varnothing$.

For $f\in\Lambda$, the \emph{plethystic substitution} $f[A]$ is obtained
by expressing $f$ in the power sums and replacing each $p_m$
by $p_m[A]$, while leaving the coefficients of $f$ unchanged.
Here $p_m[A]$ replaces every variable and parameter in $A$
by its $m$th power, leaving rational numerical coefficients unchanged.
The \emph{plethystic exponential} $\Omega[X]:=\sum_{m\geq0}h_m[X]$,
interpreted as a formal series, satisfies
\begin{equation}\label{eqOmega}
  \Omega[X+Y]=\Omega[X]\Omega[Y],
  \qquad
  \Omega[-X]=\sum_{m\geq0}(-1)^m e_m[X].
\end{equation}

Set $(a;q)_m=\prod_{i=0}^{m-1}(1-aq^i)$, and
$[m]_q=(1-q^m)/(1-q)$, with empty products equal to $1$.

Let $\widetilde H_\mu=\widetilde H_\mu[X;q,t]$ be the modified Macdonald polynomials, which form a homogeneous
basis of $\Lambda$, with $\widetilde H_\varnothing=1$.
For a single row, we recall the classical identity:
\begin{equation}\label{eqMacdonaldRow}
  \widetilde H_{(m)}[X;q,t]
  =
  (q;q)_m h_m\left[\frac{X}{1-q}\right],
  \qquad m\geq1.
\end{equation}
For a partition $\mu$, set
$B_\mu=\sum_{(i,j)\in\mu}q^{i-1}t^{j-1}$.
The operators $\nabla$, $\Delta_f$, and $\Delta'_f$, for
$f\in\Lambda$, act diagonally on this basis:
\begin{equation}\label{eqRefinedS13MacdonaldOperators}
  \Delta_f\widetilde H_\mu:=f[B_\mu]\widetilde H_\mu,
  \qquad
  \Delta'_f\widetilde H_\mu:=f[B_\mu-1]\widetilde H_\mu,
  \qquad
  \nabla\widetilde H_\mu := e_{|\mu|}[B_\mu] \widetilde H_\mu.
\end{equation}
By~\eqref{eqRefinedS13MacdonaldOperators}, we have the operator identity $\nabla=\Delta_{e_n}=\Delta'_{e_{n-1}}$ on $\Lambda^{(n)}$ for $n\geq1$.

Following Garsia and Haglund~\cite{GHag}, we define $E_{n,k}[X]$ for $1\leq k\leq n$ by
\begin{equation}\label{eqRefinedS13EDefinition}
  e_n\left[X\frac{1-z}{1-q}\right]
  =
  \sum_{k=1}^{n}\frac{(z;q)_k}{(q;q)_k}E_{n,k}[X],
  \qquad n\geq1,
\end{equation}
where $z$ is a formal variable.
Setting $z=q$ in~\eqref{eqRefinedS13EDefinition} gives \begin{equation}\label{eqRefinedS13ESum}\sum_{k=1}^{n}E_{n,k}=e_n.\end{equation}

Write $\mathrm{QSym} = \bigoplus_{m\geq 0} \mathrm{QSym}^m$ for the graded algebra of quasisymmetric
functions over $\mathbb{Q}(q,t)$.  We define Gessel's \emph{fundamental quasisymmetric function} by
\begin{equation}\label{eqRefinedS13Fundamental}
  F_{m,S}[X]
  :=
  \sum_{\substack{
    1\leq a_1\leq\cdots\leq a_m\\
    i\in S\Rightarrow a_i<a_{i+1}
  }}
  x_{a_1}\cdots x_{a_m},
  \qquad \text{ for }S\subseteq\{1,\ldots,m-1\}.
\end{equation}
\begin{definition}\label{defRefinedS13SignedWeights}
We order the \emph{signed alphabet} by $1^-<1^+<2^-<2^+<\cdots$ and set $\wt(i^-)=-x_i$ and $\wt(i^+)=qx_i$.
\end{definition}
Following Haglund et al.~\cite{HHL05},
define $U_q$ linearly by $U_q1=1$ and
\begin{equation}\label{eqRefinedS13UqDefinition}
  (U_qF_{m,S})[X]
  :=
  \sum_{1\leq a_1\leq\cdots\leq a_m}
  \ \sum_{\substack{
    c_i\in\{a_i^-,a_i^+\},\ c_1\leq\cdots\leq c_m\\
    c_i=c_{i+1}\text{ positive}\Rightarrow i\notin S\\
    c_i=c_{i+1}\text{ negative}\Rightarrow i\in S
  }}
  \prod_{i=1}^{m}\wt(c_i).
\end{equation}

\begin{proposition}\label{thmUqProperties}
The operator $U_q$ has the following properties.
\begin{enumerate}[(1)]
  \item
  For $g\in\Lambda$, we have $(U_qg)[X]=g[(q-1)X]$.

  \item
  The operator $U_q$ is invertible on $\mathrm{QSym}$.
  Moreover, $F\in\Lambda$ if and only if $U_qF\in\Lambda$.
\end{enumerate}
\end{proposition}

\begin{proof}
\textup{(1)}
By~\cite[Proposition~4.2]{HHL05}, 
the restriction to $\Lambda$ is obtained from $\omega_Yg[X+Y]$
by the ordinary substitutions $x_i\mapsto qx_i$ and
$y_i\mapsto-x_i$.
On the power-sum basis, this gives
\[
  p_m\longmapsto
  \bigl(q^m+(-1)^{m-1}(-1)^m\bigr)p_m
  =(q^m-1)p_m,
\]
which is just the plethystic substitution $X\mapsto(q-1)X$.

\textup{(2)}
By~\eqref{eqRefinedS13Fundamental}, \eqref{eqRefinedS13UqDefinition},
and~\cite[Definition~11]{GV22}, for $m\geq1$ and
$D\subseteq\{1,\ldots,m-1\}$, we have
\[
  (U_qF_{m,D})[X]
  = q^mL_{m,D}^{(-q^{-1})}[X],
\]
where $L_{m,D}^{(z)}$ is defined in~\cite[Definition~11]{GV22}.
Since $-q^{-1}$ is not a root of unity, it follows from
\cite[Theorem~4 and Remark~1]{GV22} that
$\{L_{m,D}^{(-q^{-1})}:D\subseteq\{1,\ldots,m-1\}\}$
is a basis of $\mathrm{QSym}^m$.
Together with $U_q1=1$, this proves that $U_q$ is invertible.
Meanwhile, for $g\in\Lambda$, Proposition~\ref{thmUqProperties}\textup{(1)} gives
$(U_q^{-1}g)[X]=g[X/(q-1)]$.
Thus both $U_q$ and $U_q^{-1}$ preserve $\Lambda$, and hence
$F\in\Lambda$ if and only if $U_qF\in\Lambda$.
\end{proof}

\subsection{Combinatorial objects}\label{subsecCombinatorialObjects}
Let $\pi$ be a \emph{Dyck path} of size $n$, that is, a lattice path from
$(0,0)$ to $(n,n)$ with north and east steps that stays weakly
above $y=x$.  Index its north steps from bottom to top.
The \emph{area sequence} $\mathbf a=(a_1,\ldots,a_n)$ records the
number of squares between each north step and the diagonal;
thus $a_1=0$ and $0\leq a_{i+1}\leq a_i+1$.

\begin{definition}\label{defRefinedS13Parent}
Associate a \emph{rooted forest} with this area sequence.
North step $i$ is a \emph{root} if $a_i=0$; otherwise its \emph{parent} is the
last north step $p<i$ with $a_p=a_i-1$.
Such a parent exists by the area-sequence inequalities.
Children are ordered by their indices.
\end{definition}

\begin{definition}\label{defRefinedS13DecoratedPaths}
A \emph{parking function} is a Dyck path with a positive integer label $w_i$
on each north step, with labels strictly increasing up in each column.
For $i\geq2$, north step $i$ is a \emph{rise} if $a_i=a_{i-1}+1$
and a \emph{valley} if $a_i\leq a_{i-1}$.
A valley is \emph{contractible} if $a_i<a_{i-1}$, or if
$a_i=a_{i-1}$ and $w_{i-1}<w_i$.
We decorate selected rises with $\ast$ and selected
contractible valleys with $\bullet$, which we call \emph{decorated rises} and \emph{decorated valleys}.
Here $\RR(\pi)$ is the set of
indices of decorated rises, and $\DD(\pi)$ is the set of indices
of decorated valleys.  Let $\PF_n^{\ast k,\bullet\ell}$ denote the family of parking
functions with $|\RR(\pi)|=k$
and $|\DD(\pi)|=\ell$.

A \emph{fall} is an east step followed by another east step.
Match each rise with the first subsequent fall crossing
the same diagonal.  Here the \emph{diagonal} of a rise or fall is the
line of slope $1$ through the common vertex of its two
consecutive steps.
We decorate a fall with a star from its matched decorated rise.
A \emph{star valley} is a decorated valley whose two immediately
preceding east steps form a fall decorated by a star.
Write $\DD^*(\pi)$ for the set of star-valley indices.

We define the \emph{monomial}, \emph{area}, and \emph{touch number} by
\begin{equation}\label{eqRefinedS13PathStatistics}
  x^\pi:=\prod_{i=1}^{n}x_{w_i},\qquad
  \area(\pi):=\sum_{i\notin\RR(\pi)}a_i,\qquad
  \touch(\pi):=\#\{i:a_i=0,\ i\notin\RR(\pi)\cup\DD(\pi)\}.
\end{equation}
An \emph{inversion pair} is a pair $i<j$ satisfying one of the following conditions:
\begin{equation}\label{eqRefinedS13Inversions}
\begin{alignedat}{2}
  \textup{(1)}\;& i\notin\DD(\pi), a_i=a_j, w_i<w_j;
  &\quad \textup{(2)}\;& i\notin\DD(\pi),a_i=a_j+1,w_j<w_i;\\
  \textup{(3)}\;& i\in\DD^*(\pi), a_i=a_j, j\in\DD(\pi), w_i<w_j;
  &\quad \textup{(4)}\;& i\in\DD^*(\pi),a_i=a_j, j\notin\DD(\pi), w_j<w_i.
\end{alignedat}
\end{equation}

Define
\begin{equation}\label{eqRefinedS13Dinv}
  \dinv(\pi):=
  \#\{\text{inversion pairs}\}
  -|\DD(\pi)\setminus\DD^*(\pi)|
  -\#\{i<j:i\in\RR(\pi),\ j\in\DD(\pi),\ a_i=a_j+1\}.
\end{equation}
\end{definition}

We extend $\PF_n^{\ast k,\bullet\ell}$ by allowing signed labels
on the north steps, ordered by $1^-<1^+<2^-<2^+<\cdots$.
\begin{definition}\label{defRefinedS13SignedComparison}
For signed letters $c,d$, write $c\prec d$ if $c<d$ or $c=d=i^-$ for some $i$.
\end{definition}
The labels must increase up each column with respect to $\prec$.
In the contractibility conditions of Definition~\ref{defRefinedS13DecoratedPaths}
and the four inversion conditions in~\eqref{eqRefinedS13Inversions},
we replace each comparison $<$ by $\prec$.
The definitions of area, touch number, star valleys, and the
correction terms in $\dinv$ in Definition~\ref{defRefinedS13DecoratedPaths} remain unchanged.

We denote the resulting family of \emph{signed parking functions} with
$k$ decorated rises and $\ell$ decorated contractible valleys by
$\PF_{n,\pm}^{\ast k,\bullet\ell}$.
For $\pi$ with signed labels $w_1,\ldots,w_n$, its \emph{label weight} is
$\prod_{i=1}^{n}\wt(w_i)$, where, as in Definition~\ref{defRefinedS13SignedWeights},
$\wt(i^-)=-x_i$ and $\wt(i^+)=qx_i$.

Following Interdonato and Iraci~\cite{Thetaconj},
the \emph{reading order} of the north steps is the lexicographic order on
$(a_i,\mathbf1_{i\in\DD(\pi)},i)$.

\begin{proposition}\label{propUqSignedPaths}
For integers $d\geq j\geq1$ and $k,\ell\geq0$, set
$n=d+k+\ell$.  We have
\[
  U_q\left(\sum_{\substack{\pi\in\PF_n^{\ast k,\bullet\ell}\\
                       \touch(\pi)=j}}
       q^{\dinv(\pi)}t^{\area(\pi)}x^\pi \right)
  =
  \sum_{\substack{\pi\in\PF_{n,\pm}^{\ast k,\bullet\ell}\\
                   \touch(\pi)=j}}
       q^{\dinv(\pi)}t^{\area(\pi)} \prod_{i=1}^{n}\wt(w_i).
\]
\end{proposition}

\begin{proof}
For a parking function $\pi$ with labels $1,\ldots,n$, each used
once, let $\operatorname{ides}(\pi)$ be the inverse descent set of its reverse
reading word: $i\in\operatorname{ides}(\pi)$ if and only if $i+1$ precedes
$i$ in that word.
Grouping parking functions by their standardization, as in
\cite[Definition~3.13 and Lemma~3.14]{Thetaconj}, gives
\[
  \sum_{\substack{\pi\in\PF_n^{\ast k,\bullet\ell}\\
                   \touch(\pi)=j}}
  q^{\dinv(\pi)}t^{\area(\pi)}x^\pi
  =
  \sum_{\substack{\pi\in\PF_n^{\ast k,\bullet\ell}\\
                   \touch(\pi)=j,\ \{w_1,\ldots,w_n\}=[n]}}
  q^{\dinv(\pi)}t^{\area(\pi)}F_{n,\operatorname{ides}(\pi)}[X],
\]
so the left-hand side is quasisymmetric.

Fix $\pi$ with labels $1,\ldots,n$, each used once, and replace
label $i$ by $c_i$, where $c_1\leq\cdots\leq c_n$.
By the definition of $\operatorname{ides}(\pi)$ and
\eqref{eqRefinedS13UqDefinition}, within each block of equal
replacement letters, the original labels decrease in reading
order if the common letter is positive, and increase if it
is negative.

In the column conditions, the comparison is from the earlier
step to the later step in reading order.
The comparisons at selected decorated valleys and
in~\eqref{eqRefinedS13Inversions} have the same direction.
If the replacement letters are distinct, they have the same
relative order as the original labels.
For equal positive letters, both the original $<$ comparison
and the new $\prec$ comparison are false; for equal negative
letters, both are true.
Thus the replacement gives a signed parking function with
the same inversion pairs.
The path and decorations are unchanged, so area, touch number,
and the correction terms in~\eqref{eqRefinedS13Dinv}
are also preserved.

Conversely, number the steps of a signed parking function
$1,\ldots,n$ in increasing order of their signed labels,
using reverse reading order for equal positive labels and
reading order for equal negative labels.
This uniquely recovers $\pi$ and $(c_i)$; the same comparisons
show that $\pi\in\PF_n^{\ast k,\bullet\ell}$.
Therefore
\begin{align*}
  U_q\left(
    \sum_{\substack{\pi\in\PF_n^{\ast k,\bullet\ell}\\
                     \touch(\pi)=j}}
    q^{\dinv(\pi)}t^{\area(\pi)}x^\pi
  \right)
  &=
  \sum_{\substack{\pi\in\PF_n^{\ast k,\bullet\ell}\\
                   \touch(\pi)=j,\ \{w_1,\ldots,w_n\}=[n]}}
  q^{\dinv(\pi)}t^{\area(\pi)}
  (U_qF_{n,\operatorname{ides}(\pi)})[X]\\
  &=
  \sum_{\substack{\pi\in\PF_{n,\pm}^{\ast k,\bullet\ell}\\
                   \touch(\pi)=j}}
  q^{\dinv(\pi)}t^{\area(\pi)}
  \prod_{i=1}^{n}\wt(w_i).\qedhere
\end{align*}
\end{proof}

Recall~\eqref{eqRefinedS13RIntroduction}:
\begin{equation*}
    R(X;u,v,r,s):=1+
  \sum_{j\geq1}\sum_{d\geq j}\sum_{k,\ell\geq0}
  u^d v^j r^k s^\ell
  \sum_{\substack{\pi\in\PF_{d+k+\ell}^{\ast k,\bullet\ell}\\
                   \touch(\pi)=j}}
  q^{\dinv(\pi)}t^{\area(\pi)}x^\pi.
\end{equation*}

Here the constant term $1$ records the empty path.
By Proposition~\ref{propUqSignedPaths},
\(
  (U_qR)(X_n;u,v,r,s)-1
\)
enumerates the nonempty signed paths with labels
$1^-,1^+,\ldots,n^-,n^+$.
The set of such paths in $\PF_{d+k+\ell,\pm}^{\ast k,\bullet\ell}$ is
denoted by $\PF_{d+k+\ell,\pm}^{\ast k,\bullet\ell}(X_n)$.

\begin{example}\label{exSignedDecoratedPath}
  The parking functions in Figure~\ref{figSignedDecoratedPath}
have area sequence $(0,1,2,1,0,1,1,2,0,0)$ and decorations
$\RR(\pi)=\{3,6,8\}$, $\DD(\pi)=\{4,7,9\}$ and
$\DD^*(\pi)=\{4,9\}$.
Thus $k=\ell=3$, $d=10-k-\ell=4$, and $\touch(\pi)=3$,
giving the common factor $u^4v^3r^3s^3$.
Their roots are steps $1,5,9,10$, of which $1,5,10$ are
undecorated.
Root $1$ has children $2,4$ and descendants $2,3,4$;
root $5$ has children $6,7$ and descendants $6,7,8$.

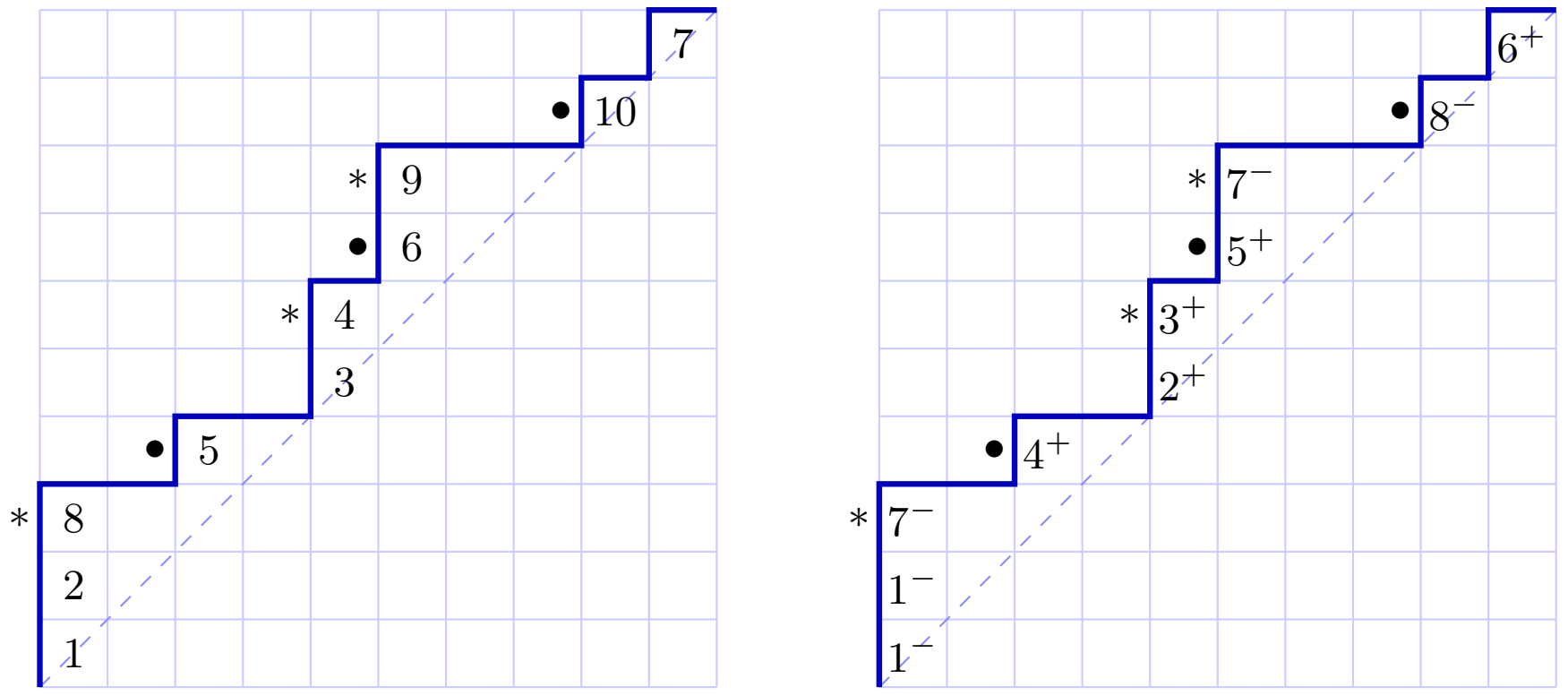
\begin{figure}[htbp]
  \centering
  \begin{tikzpicture}[x=.5cm,y=.5cm,font=\small]
    \foreach \shift in {0cm,6.2cm}{
      \begin{scope}[xshift=\shift]
        \draw[step=1,blue!20,thin] (0,0) grid (10,10);
        \draw[dashed,blue!45] (0,0)--(10,10);
        \draw[very thick,blue!75!black]
          (0,0)--(0,3)--(2,3)--(2,4)--(4,4)--(4,6)
          --(5,6)--(5,8)--(8,8)--(8,9)--(9,9)
          --(9,10)--(10,10);

        \node at (-.3,2.5) {$\ast$};
        \node at (3.7,5.5) {$\ast$};
        \node at (4.7,7.5) {$\ast$};
        \node at (1.7,3.5) {$\bullet$};
        \node at (4.7,6.5) {$\bullet$};
        \node at (7.7,8.5) {$\bullet$};
      \end{scope}
    }

    \foreach \x/\y/\w in {
      .5/.5/1,
      .5/1.5/2,
      .5/2.5/8,
      2.5/3.5/5,
      4.5/4.5/3,
      4.5/5.5/4,
      5.5/6.5/6,
      5.5/7.5/9,
      8.5/8.5/10,
      9.5/9.5/7
    }
      \node at (\x,\y) {$\w$};
    \node at (5,-.75) {(1) Positive labels};

    \begin{scope}[xshift=6.2cm]
      \foreach \x/\y/\w in {
        .5/.5/{1^-},
        .5/1.5/{1^-},
        .5/2.5/{7^-},
        2.5/3.5/{4^+},
        4.5/4.5/{2^+},
        4.5/5.5/{3^+},
        5.5/6.5/{5^+},
        5.5/7.5/{7^-},
        8.5/8.5/{8^-},
        9.5/9.5/{6^+}
      }
        \node at (\x,\y) {$\w$};
      \node at (5,-.75) {(2) Signed labels};
    \end{scope}
  \end{tikzpicture}
  \caption{Positive and signed labelings of the same decorated Dyck path.}
  \label{figSignedDecoratedPath}
\end{figure}

\begin{enumerate}[(1)]
  \item
  The positive labeling is $(1,2,8,5,3,4,6,9,10,7)$.
  The four inversion types in~\eqref{eqRefinedS13Inversions} contribute $10,3,1,2$ pairs,
  respectively.
  The only decorated valley outside $\DD^*(\pi)$ is step $7$,
  and the pairs in the last term of~\eqref{eqRefinedS13Dinv} are
  $(3,4),(3,7),(6,9)$.
  Hence
  \[
    \area(\pi)=a_2+a_4+a_7=3,\qquad
    \dinv(\pi)=10+3+1+2-1-3=12.
  \]
  \item
  The signed labeling is
  $(1^-,1^-,7^-,4^+,2^+,3^+,5^+,7^-,8^-,6^+)$.
  The column conditions hold since
  $1^-\prec1^-\prec7^-$, $2^+\prec3^+$, and
  $5^+\prec7^-$.
  The decorated valleys at steps $4$ and $9$ are contractible
  independently of their labels, while step $7$ is
  contractible because $3^+\prec5^+$.
  Using $\prec$, the four inversion counts are again
  $10,3,1,2$.
  The decorations and correction terms are unchanged, so
  $\area(\pi)=3$, $\dinv(\pi)=12$, and $\touch(\pi)=3$.
\end{enumerate}
\end{example}

\begin{lemma}\label{lemForestStarValleys}
\leavevmode
\begin{enumerate}[(1)]
  \item
  A north step is a rise if and only if it is the first
  child of its parent.

  \item
  Let $j\in\DD(\pi)$ satisfy $a_j<a_{j-1}$, and let $p$ be
  the last north step before $j$ with the same area as $j$.
  Then $j\in\DD^*(\pi)$ if and only if $p+1\in\RR(\pi)$.
\end{enumerate}
\end{lemma}

\begin{proof}
\textup{(1)}
If $i$ is a rise, its parent is $i-1$, and it is necessarily
the first child.
Conversely, let $i$ be the first child of $p$.
By Definition~\ref{defRefinedS13Parent}, $a_i=a_p+1$, and no entry
strictly between $p$ and $i$ equals $a_p$.
Since area entries increase by at most one at each step,
$a_{p+1}=a_p+1$.
Thus $p+1$ is a child of $p$, so $i=p+1$ and $i$ is a rise.

\textup{(2)}
Since $a_{j-1}>a_j$ and area entries increase by at most one,
$p$ exists.
Its choice gives $a_{p+1}=a_p+1$ and $a_i>a_p$ for $p<i<j$.
Therefore, the rise $p+1$ is matched with the fall formed by
the last two east steps before $j$.
By Definition~\ref{defRefinedS13DecoratedPaths},
$j$ is a star valley precisely when this fall carries a star,
that is, when $p+1\in\RR(\pi)$.
\end{proof}

\section{The Theta operator and the Dyck path algebra}\label{secThetaAlgebra}
In this section, we prove Theorem \ref{thmIntroAlgebraicExpansion}
and derive the recurrence~\eqref{eqAlgRecursion} for $\mathcal L_I$.
Subsection~\ref{subsecThetaOperators} recalls the Theta
operators of D'Adderio et al.~\cite{DIVW21}, and
Subsection~\ref{subsecDyckPathAlgebra} recalls the Dyck path
algebra of Carlsson and Mellit~\cite{CM18}.
In Subsection~\ref{subsecThetaDyckConnection}, we use the
endomorphism and its action on $V_\ell$ constructed by
D'Adderio et al.~\cite{DIIP} to prove the recurrence
for $f_j$ in Theorem~\ref{thmFjRecursion}.
Finally, Subsection~\ref{subsecAlgebraicRecurrence} specializes
this recurrence to finite alphabets and proves
Theorem~\ref{thmIntroAlgebraicExpansion}.

\subsection{The Theta operator}\label{subsecThetaOperators}
D'Adderio et al.~\cite{DIVW21} introduced the
Theta operators and used them to formulate
refinements of the Delta conjecture of Haglund et al.~\cite{HRW18}.

For a partition $\mu$, set
\begin{equation*}
  \Pi_\mu
  :=
  \prod_{\substack{(i,j)\in\mu\\(i,j)\neq(1,1)}}
  (1-q^{i-1}t^{j-1}),
  \qquad
  \Pi_\varnothing:=1,
\end{equation*}
and define the invertible operator $\mathbf{\Pi}$ by $
  \mathbf{\Pi}\widetilde H_\mu[X;q,t]
  :=
  \Pi_\mu\widetilde H_\mu[X;q,t]$.

Let $f$ and $F$ be homogeneous
symmetric functions of degrees $a$ and $b$, respectively.  Define
\begin{equation}\label{eqRefinedS13ThetaDefinition}
  \Theta_f F
  :=
  \begin{cases}
    0,
      & a>0 \text{ and } b=0,\\
    fF,
      & a=b=0,\\
    \mathbf{\Pi}
    \left(
      f\left[\dfrac{X}{(1-q)(1-t)}\right]
      \mathbf{\Pi}^{-1}F
    \right),
      & b>0.
  \end{cases}
\end{equation}
We extend~\eqref{eqRefinedS13ThetaDefinition} linearly in $f$ and $F$. 
In particular, the Theta operators commute, $\Theta_{e_0}$ is the
identity, and $\Theta_{e_m}1=0$ for $m>0$.

D'Adderio et al.~\cite[Theorem~3.1 and Corollary~3.2]{DIVW21} proved the following identities.

\begin{theorem}[D'Adderio et al.~\cite{DIVW21}]\label{thmThetaDeltaIdentity}
For integers $d\geq 1$ and $k\geq 0$, we have
\begin{equation*}
  \Theta_{e_k}\nabla e_d
  =\Delta'_{e_{d-1}}e_{d+k},\qquad \sum_{j=1}^{d}\Theta_{e_k}\nabla E_{d,j}
  =\Delta'_{e_{d-1}}e_{d+k}.
\end{equation*}
\end{theorem}

In our notation, Haglund~\cite[equation~(215)]{Hag04}
gives the following identity.
\begin{proposition}[Haglund~\cite{Hag04}]\label{propHag04}
For integers $1\leq j\leq d$, we have
\begin{equation*}
  \nabla E_{d,j}=t^{d-j}\Theta_{h_{d-j}}\widetilde H_{(j)}.
\end{equation*}
\end{proposition}

Following D'Adderio et al.~\cite{DIIP}, let $v$ be a formal variable and define
\begin{equation}\label{eqRefinedS13ThetaSeries}\Theta(v):=\sum_{m\geq 0}v^m\Theta_{e_m}.\end{equation}
The identity
\[
  \bigl(\sum_{m\geq 0}v^m e_m\bigr)
  \bigl(\sum_{m\geq 0}(-v)^m h_m\bigr)=1
\]
gives
\begin{equation}\label{eqRefinedS13ThetaSeriesInverse}\Theta(v)^{-1}=\sum_{m\geq 0}(-v)^m\Theta_{h_m},\qquad \Theta(-v)^{-1}=\sum_{m\geq 0}v^m\Theta_{h_m}.\end{equation}

By Proposition~\ref{propHag04} and~\eqref{eqRefinedS13ThetaSeries}--\eqref{eqRefinedS13ThetaSeriesInverse}, we can rewrite~\eqref{eqThetaalg} as
\begin{align*}
  L[X;u,v,r,s]
  &=
  \sum_{j\geq 1}\sum_{d\geq j}\sum_{k,\ell\geq 0}
   u^d v^j r^k s^\ell t^{d-j}
  \bigl(
    \Theta_{e_\ell}\Theta_{e_k}
    \Theta_{h_{d-j}}\widetilde H_{(j)}
  \bigr)[X]\\
  &=
  \sum_{j\geq 1}(uv)^j
  \sum_{d\geq j}\sum_{k,\ell\geq 0}
  (ut)^{d-j}r^k s^\ell
  \bigl(
    \Theta_{e_\ell}\Theta_{e_k}
    \Theta_{h_{d-j}}\widetilde H_{(j)}
  \bigr)[X]\\
  &=
  \sum_{j\geq 1}(uv)^j
  \bigl(
    \Theta(r)\Theta(s)\Theta(-ut)^{-1}
    \widetilde H_{(j)}
  \bigr)[X].
\end{align*}

\subsection{The Dyck path algebra}\label{subsecDyckPathAlgebra}
The \emph{Dyck path algebra} $\mathbb{A}_{q,t}$ was introduced by
Carlsson and Mellit~\cite{CM18}.  It is the quotient of the
path algebra over $\mathbb{Q}(q,t)$ of the quiver with vertices
$\ell\geq 0$, loops $T_i$ at vertex $\ell$ for $1\leq i<\ell$,
arrows $d_+,d_+^*$ from $\ell$ to $\ell+1$, and arrows $d_-$
from $\ell$ to $\ell-1$ for $\ell\geq 1$, by the relations
of Carlsson and Mellit~\cite[Definition~5.1 and Section~6]{CM18}.
Here we use the notation and normalization recalled
by D'Adderio et al.~\cite{DIIP}.

\begin{definition}\label{defd-d+}
For $\ell\geq0$ set
$V_\ell:=\mathbb{Q}(q,t)[y_1,\ldots,y_\ell]\otimes\Lambda[X]$
and $V_*:=\bigoplus_{\ell\geq0}V_\ell$.
We write an element of $V_\ell$ as $F[X;y_1,\ldots,y_\ell]$.

The actions of $\mathbb{A}_{q,t}$ on $V_*$ are defined as follows.
For $F\in V_\ell$, the generators act by
\begin{align}
  T_iF
  &:=
  s_iF+(1-q)y_i\frac{F-s_iF}{y_i-y_{i+1}},
  \qquad 1\leq i<\ell,\label{eqTi}\\
  d_-F
  &:=
  \langle y_\ell^0\rangle
  \left(
    F[X-(q-1)y_\ell;y_1,\ldots,y_\ell]\,
    \Omega[-y_\ell^{-1}X]
  \right),\label{eqdminus}\\
  d_+F
  &:=
  -T_1\cdots T_\ell
  \left(
    y_{\ell+1}
    F[X+(q-1)y_{\ell+1};y_1,\ldots,y_\ell]
  \right),\notag\\
  d_+^*F
  &:=
  F[X+(q-1)ty_1;y_2,\ldots,y_{\ell+1}].
  \label{eqdplusstar}
\end{align}
Here $T_i$ preserves $V_\ell$,
$d_-$ maps $V_\ell$ to $V_{\ell-1}$ for $\ell\geq1$,
and $d_+,d_+^*$ map $V_\ell$ to $V_{\ell+1}$.
The operator $s_i$ interchanges $y_i$ and $y_{i+1}$,
and $\langle y_\ell^0\rangle$ extracts the constant term
in $y_\ell$, with
$\Omega[-y_\ell^{-1}X]
=\sum_{m\geq0}(-1)^m y_\ell^{-m}e_m[X]$.
Empty products of the $T_i$ are understood to be the identity.

For $\ell\geq1$ we define the algebra elements
$y_1,\ldots,y_\ell$ and $z_1,\ldots,z_\ell$:
\begin{align}
  y_1
  &:=
  \frac{(d_+d_--d_-d_+)T_{\ell-1}\cdots T_1}
       {q^{\ell-1}(q-1)},
  \qquad
  y_{i+1}:=qT_i^{-1}y_iT_i^{-1}
  \quad(1\leq i<\ell),\label{eqyi}\\
  z_1
  &:=
  \frac{q^\ell}{1-q}
  (d_+^*d_--d_-d_+^*)T_{\ell-1}^{-1}\cdots T_1^{-1},
  \qquad
  z_{i+1}:=q^{-1}T_i z_iT_i
  \quad(1\leq i<\ell).\label{eqzi}
\end{align}
Under the action defined above, each algebra element $y_i$ acts
by multiplication by the variable $y_i$.
\end{definition}

For $j\geq 0$, let 
\begin{equation}\label{eqRefinedS13PDefinition}
  P_j:=(-1)^j y_1\cdots y_j,
\end{equation}
where the empty product gives $P_0=1$.

\begin{proposition}\label{propSeed}
For every integer $j\geq 1$, we have
$P_j=-y_1d_+^*P_{j-1}$ and $d_-^jP_j=\widetilde H_{(j)}$.
\end{proposition}

\begin{proof}
Since $P_{j-1}$ in~\eqref{eqRefinedS13PDefinition} is independent of $X$, \eqref{eqdplusstar}
gives
\[
  P_j
  =-y_1\left((-1)^{j-1}y_2\cdots y_j\right)
  =-y_1d_+^*P_{j-1}.
\]

To prove the second identity, let $\ell\geq 1$ and $m\geq 0$.
Taking the coefficient of $y_\ell^{m+1}$, noting $\frac{X}{1-q}-X=\frac{qX}{1-q}$ and
using~\eqref{eqOmega}, we obtain
\begin{align*}
  \sum_{a=1}^{m+1}
  (-1)^{a-1}e_a[X]
  h_{m+1-a}\left[\frac{X}{1-q}\right]
  =
  h_{m+1}\left[\frac{X}{1-q}\right]
  -\langle y_\ell^{m+1}\rangle
  \left(
    \Omega\left[\frac{y_\ell X}{1-q}\right]
    \Omega[-y_\ell X]
  \right)
  =
  (1-q^{m+1})h_{m+1}\left[\frac{X}{1-q}\right].
\end{align*}

Equation~\eqref{eqdminus} and the identity
\begin{equation*}
  h_m\left[\frac{X}{1-q}+y_\ell\right]
  =
  \sum_{b=0}^{m}
  y_\ell^b h_{m-b}\left[\frac{X}{1-q}\right]
\end{equation*}
yield
\begin{equation}\label{eqRefinedS13SeedStep}
\begin{aligned}
  d_-\left(
    P_\ell(q;q)_m h_m\left[\frac{X}{1-q}\right]
  \right)
  &=
  (-1)^\ell(q;q)_m
  \left(\prod_{i=1}^{\ell-1}y_i\right)
  \langle y_\ell^0\rangle
  \left(
    y_\ell h_m\left[\frac{X}{1-q}+y_\ell\right]
    \Omega[-y_\ell^{-1}X]
  \right)\\
  &=
  (-1)^{\ell-1}(q;q)_m
  \left(\prod_{i=1}^{\ell-1}y_i\right)
  \sum_{a=1}^{m+1}
  (-1)^{a-1}e_a[X]
  h_{m+1-a}\left[\frac{X}{1-q}\right]\\
  &=
  P_{\ell-1}(q;q)_{m+1}
  h_{m+1}\left[\frac{X}{1-q}\right].
\end{aligned}
\end{equation}
Starting with
$P_j=P_j(q;q)_0h_0\left[\frac{X}{1-q}\right]$,
repeated application of~\eqref{eqRefinedS13SeedStep} gives
\begin{equation*}
  d_-^mP_j
  =
  P_{j-m}(q;q)_m h_m\left[\frac{X}{1-q}\right],
  \qquad 0\leq m\leq j.
\end{equation*}
Therefore, by~\eqref{eqMacdonaldRow},
\[
  d_-^jP_j = P_0(q;q)_j h_j\left[\frac{X}{1-q}\right] = (q;q)_j h_j\left[\frac{X}{1-q}\right] = \widetilde H_{(j)}[X;q,t].\qedhere
\]
\end{proof}

Let $n\geq 0$ and consider the finite alphabet $X_n$.
For $x\in X_n$ and $I\subseteq X_n$, define
\begin{equation}\label{eqRefinedS13Eta}\eta_x(X_n):=\prod_{x'\in X_n\setminus\{x\}}\frac{x-qx'}{x-x'}
,\qquad \eta_I(X_n):=\prod_{x\in I}\prod_{x'\in X_n\setminus I}
\frac{x-qx'}{x-x'}.\end{equation}
Empty products are understood to be $1$.

\begin{proposition}\label{propDminusFiniteAlphabet}
Let $j\geq 1$ and let $F\in V_j$ be divisible by $y_j$.
Then
\begin{equation*}
  (d_-F)[(q-1)X_n;y_1,\ldots,y_{j-1}]
  =
  (1-q)\sum_{x\in X_n}\eta_x(X_n)
  F[(q-1)(X_n\setminus\{x\});y_1,\ldots,y_{j-1},x].
\end{equation*}
\end{proposition}

\begin{proof}
We first verify the identity
\begin{equation}\label{eqEtaPartialFractions}
  \prod_{x\in X_n}\frac{y_j-qx}{y_j-x}
  =
  q^n
  +(1-q)\sum_{x\in X_n}
  \eta_x(X_n)\frac{y_j}{y_j-x}.
\end{equation}
We regard both sides of \eqref{eqEtaPartialFractions} as rational functions of \(y_j\), 
with \(q,x_1,\ldots,x_n\) as parameters. 
After multiplying both sides by $\prod_{x\in X_n}(y_j-x)$, we obtain two polynomials in \(y_j\) of degree at most \(n\). 
At \(y_j=0\), both polynomials have value $(-q)^n\prod_{x\in X_n}x$.
Now let \(x\in X_n\) and set \(y_j=x\) in these polynomials. 
On the right, the term containing \(q^n\) and all summands except the one indexed by \(x\) vanish. 
By \eqref{eqRefinedS13Eta}, the remaining summand is
$$(1-q)x \eta_x(X_n)
  \prod_{x'\in X_n\setminus\{x\}}(x-x')
=(1-q)x
  \prod_{x'\in X_n\setminus\{x\}}(x-qx')$$
which is also the value of the polynomial on the left. 
Thus the two polynomials agree at the \(n+1\) distinct values \(0,x_1,\ldots,x_n\), 
proving \eqref{eqEtaPartialFractions}.

We use the formal expansion
$\frac{y_j}{y_j-x}=\sum_{a\geq0}x^a y_j^{-a}$.
Multiplying a polynomial in $y_j$ by this series and extracting
the constant term therefore evaluates the polynomial at $y_j=x$.
By the definition of \(\Omega\) and \eqref{eqOmega}, we have
$$\Omega[(1-q)X_n/y_j]
=\Omega[X_n/y_j]\,\Omega[-qX_n/y_j]
=\prod_{x\in X_n}\frac{1-qx/y_j}{1-x/y_j}
=\prod_{x\in X_n}\frac{y_j-qx}{y_j-x}.$$
Using~\eqref{eqdminus} and~\eqref{eqEtaPartialFractions}, we obtain
\begin{align*}
  (d_-F)[(q-1)X_n;y_1,\ldots,y_{j-1}]
  ={}&
  \langle y_j^0\rangle
  \left(
    F[(q-1)(X_n-y_j);y_1,\ldots,y_j]\,
    \Omega[(1-q)X_n/y_j]
  \right)\\
  ={}&
  \langle y_j^0\rangle
  \left(
    \prod_{x\in X_n}\frac{y_j-qx}{y_j-x}\,
    F[(q-1)(X_n-y_j);y_1,\ldots,y_j]
  \right)\\
  ={}&
  (1-q)\sum_{x\in X_n}\eta_x(X_n)
  \langle y_j^0\rangle
  \left(
    \frac{y_j}{y_j-x}
    F[(q-1)(X_n-y_j);y_1,\ldots,y_j]
  \right)\\
  ={}&
  (1-q)\sum_{x\in X_n}\eta_x(X_n)
  F[(q-1)(X_n\setminus\{x\});y_1,\ldots,y_{j-1},x].
\end{align*}
The term $q^n$ contributes zero since
$F[(q-1)(X_n-y_j);y_1,\ldots,y_j]$
is a polynomial in $y_j$ divisible by $y_j$.
\end{proof}

Repeated application of Proposition~\ref{propDminusFiniteAlphabet}
gives the following corollary.

\begin{corollary}\label{corDminusFiniteAlphabet}
Let $j\geq 1$ and let $F\in V_j$ be symmetric in
$y_1,\ldots,y_j$.  Then
\begin{equation*}
  \bigl(d_-^j(P_j F)\bigr)[(q-1)X_n]
  =
  (-1)^j(q;q)_j
  \sum_{\substack{I\subseteq X_n\\|I|=j}}
  \left(\prod_{x\in I}x\right)
  \eta_I(X_n)
  F[(q-1)(X_n\setminus I);I].
\end{equation*}
For $I=\{x_{i_1},\ldots,x_{i_j}\}$, we use the notation
\(
  F[(q-1)(X_n\setminus I);I]
  :=
  F[(q-1)(X_n\setminus I);x_{i_1},\ldots,x_{i_j}].
\)
Since $F$ is symmetric in $y_1,\ldots,y_j$, this evaluation is independent of the ordering of $I$.
\end{corollary}

\begin{proof}
Applying~\eqref{eqEtaPartialFractions} to $I\subseteq X_n$ of size $j$ and letting $y_j\to\infty$ gives
\begin{equation}\label{eqEtaSubsetSum}
  (q;q)_{j-1}(1-q)\sum_{x\in I}\eta_x(I)=(q;q)_{j-1}(1-q^j)=(q;q)_{j}.
\end{equation}
Equation~\eqref{eqRefinedS13Eta} also gives
\begin{equation}\label{eqEtaRemoval}
  \eta_x(X_n)
  \eta_{I\setminus\{x\}}(X_n\setminus\{x\})
  =
  \eta_I(X_n)\eta_x(I),
  \qquad x\in I.
\end{equation}

We now prove the identity by induction on $j$.
For $j=1$, it follows from
Proposition~\ref{propDminusFiniteAlphabet} applied to $-y_1F$.
Suppose $j\geq 2$ and assume the identity holds for $j-1$.
During the first $j-1$ applications of $d_-$, the variable $y_1$
is treated as a parameter.  Hence multiplication by $y_1$ and
substitution of $y_1=x$ commute with these applications.
Applying Proposition~\ref{propDminusFiniteAlphabet} to the final
application of $d_-$, followed by the induction hypothesis in
$y_2,\ldots,y_j$, gives
\begin{alignat*}{2}
  &\quad\bigl(d_-^j(P_jF)\bigr)[(q-1)X_n]\\*
  &={}
  -(1-q)\sum_{x\in X_n}x\eta_x(X_n)
  \bigl(
    d_-^{j-1}((-1)^{j-1}y_2\cdots y_jF)
  \bigr)[(q-1)(X_n\setminus\{x\});x]\\
  &={}
  (-1)^j(1-q)(q;q)_{j-1}
  \sum_{x\in X_n}
  \sum_{\substack{I\subseteq X_n\setminus\{x\}\\|I|=j-1}}
  x\left(\prod_{x'\in I}x'\right)
  \eta_x(X_n)\eta_I(X_n\setminus\{x\})
   F[(q-1)(X_n\setminus(I\cup\{x\}));I\cup\{x\}]\\
  &={}
  (-1)^j(1-q)(q;q)_{j-1}
  \sum_{\substack{I\subseteq X_n\\|I|=j}}
  \left(\prod_{x\in I}x\right)\eta_I(X_n)
  F[(q-1)(X_n\setminus I);I]
  \sum_{x\in I}\eta_x(I)
  &\text{(by~\eqref{eqEtaRemoval})}&\\
  &={}
  (-1)^j(q;q)_{j}
  \sum_{\substack{I\subseteq X_n\\|I|=j}}
  \left(\prod_{x\in I}x\right)\eta_I(X_n)
  F[(q-1)(X_n\setminus I);I].
  &\text{(by~\eqref{eqEtaSubsetSum})}&
\end{alignat*}
\end{proof}

\subsection{The connection between the Theta operators and the Dyck path algebra}\label{subsecThetaDyckConnection}
D'Adderio et al.~\cite{DIIP} constructed the following
endomorphism of $\mathbb{A}_{q,t}[[u]]$ and related it to the Theta operators.

\begin{theorem}[D'Adderio et al.~\cite{DIIP}]\label{DIIPthm}
Let $u$ be a formal parameter.
\begin{enumerate}
  \item
  There is an algebra endomorphism
  $\boldsymbol{\theta}_u:\mathbb{A}_{q,t}[[u]]
    \to\mathbb{A}_{q,t}[[u]]$ determined by
  \begin{align*}
    \boldsymbol{\theta}_u(T_i)=T_i,\quad
    \boldsymbol{\theta}_u(d_+)
      =\left(1-\frac{u}{qt}(1+uy_1)^{-1}y_1z_1\right)d_+,\quad
    \boldsymbol{\theta}_u(d_-)=d_-,\quad
    \boldsymbol{\theta}_u(d_+^*)
      =\left(1+u\frac{qt-z_1}{qt}y_1\right)^{-1}d_+^*.
  \end{align*}
  All inverses involving $u$ are expanded as formal geometric
series in $u$;
  \item
  There is a
  $\mathbb{Q}(q,t)[[u]]$-linear map
  $\boldsymbol{\theta}'_u:
    \bigoplus_{j\geq 0}V_j[[u]]
   \to
    \bigoplus_{j\geq 0}V_j[[u]]$
  satisfying $\boldsymbol{\theta}'_u(1)=1$ and
  $\boldsymbol{\theta}'_u(LF)
    =
    \boldsymbol{\theta}_u(L)\,
    \boldsymbol{\theta}'_u(F)$
  for $L\in\mathbb{A}_{q,t}$ and $F\in\bigoplus_{j\geq 0}V_j[[u]]$.

  \item
  On $V_0=\Lambda$, the map $\boldsymbol{\theta}'_u$ is
  $\Theta(u)=\sum_{m\geq 0}u^m\Theta_{e_m}$.
\end{enumerate}
\end{theorem}
Whenever $y_1^{-1}$ occurs, we use the extension of D'Adderio et al.
\cite[Definition~3.10 and Proposition~5.9]{DIIP}, in which
$\boldsymbol{\theta}'_u(y_1^{-1}F) = \boldsymbol{\theta}_u(y_1)^{-1}\boldsymbol{\theta}'_u(F)$.

The map $\boldsymbol{\theta}'_u$ and its formal inverse preserve $V_\ell[[u]]$.
Indeed, for $F\in\Lambda$,
$(d_+^*)^\ell F=F[X+(q-1)t(y_1+\cdots+y_\ell)]$.
Induction on the degree in $X$ shows that $V_\ell$ is generated
from $V_0$ by $d_+^*$ and multiplication by the $y_i$.
Indeed, the difference between the displayed expression and $F[X]$
has smaller degree in $X$, with polynomial coefficients in the $y_i$;
the degree-zero case consists of polynomials in the $y_i$.
Their images under $\boldsymbol{\theta}_u$ preserve polynomials
coefficientwise. Hence Theorem~\ref{DIIPthm}(2)--(3) implies that
$\boldsymbol{\theta}'_u$ maps $V_\ell[[u]]$ into $V_\ell[[u]]$ for every
$\ell\geq0$.
Since $\boldsymbol{\theta}'_u$ is the identity modulo $u$,
its formal inverse also maps $V_\ell[[u]]$ into $V_\ell[[u]]$.

\begin{corollary}\label{corThetaIdentities}
We have
\begin{align}
  \boldsymbol{\theta}_u(y_1)
  &=
  \frac{1}{qt}(1+uy_1)^{-1}y_1
  \bigl(qt+u(qt-z_1)y_1\bigr),
  \label{eqThetaY}\\
  \boldsymbol{\theta}_u
  \left(\frac{z_1}{qt}\right)
  &=
  \bigl(qt+u(qt-z_1)y_1\bigr)^{-1}z_1,
  \label{eqRefinedS13ThetaZ}\\
  \boldsymbol{\theta}_u
  \left(\frac{y_1z_1}{qt}\right)
  &=
  \frac{1}{qt}(1+uy_1)^{-1}y_1z_1,
  \label{eqRefinedS13ThetaYZ}\\
  \boldsymbol{\theta}_u
  \left(\frac{qt-z_1}{qt}y_1\right)
  &=
  \frac{qt-z_1}{qt}y_1,
  \label{eqThetaInvariant}\\
  \boldsymbol{\theta}_u(y_1^{-1})
  &=
  qt\bigl(qt+u(qt-z_1)y_1\bigr)^{-1}
  (y_1^{-1}+u).
  \label{eqRefinedS13ThetaYInverse}
\end{align}
\end{corollary}

\begin{proof}
By~\cite[Lemma~5.5]{CM18}, we have $d_-y_1=y_1d_-$.
Apply the same relation to the conjugate action of $\mathbb{A}_{q^{-1}}$ in \cite[Theorem 6.1]{CM18}. 
In this action, \(y_1\) corresponds to \(z_1\), and \(d_-^*=d_-\). 
Hence \(d_-z_1=z_1d_-\).  Both relations are used as maps
$V_{\ell+1}\to V_\ell$ for $\ell\geq1$.
Using these relations,
Theorem~\ref{DIIPthm}\textup{(1)}, and~\eqref{eqyi}, we obtain
\begin{align*}
  \boldsymbol{\theta}_u(y_1)
  =
  \frac{
    \bigl(
      \boldsymbol{\theta}_u(d_+)d_-
      -d_-\boldsymbol{\theta}_u(d_+)
    \bigr)T_{\ell-1}\cdots T_1
  }{q^{\ell-1}(q-1)}
  =
  \left(
    1-\frac{u}{qt}(1+uy_1)^{-1}y_1z_1
  \right)y_1
  =
  \frac{1}{qt}(1+uy_1)^{-1}y_1
  \bigl(qt+u(qt-z_1)y_1\bigr).
\end{align*}
Similarly, \eqref{eqzi} gives
\begin{align*}
  \boldsymbol{\theta}_u
  \left(\frac{z_1}{qt}\right)
  =
  \frac{q^{\ell-1}}{t(1-q)}
  \bigl(
    \boldsymbol{\theta}_u(d_+^*)d_-
    -d_-\boldsymbol{\theta}_u(d_+^*)
  \bigr)
  T_{\ell-1}^{-1}\cdots T_1^{-1}=
  \bigl(qt+u(qt-z_1)y_1\bigr)^{-1}z_1.
\end{align*}
Multiplying~\eqref{eqThetaY} and~\eqref{eqRefinedS13ThetaZ} yields~\eqref{eqRefinedS13ThetaYZ}:
\[
  \boldsymbol{\theta}_u
  \left(\frac{y_1z_1}{qt}\right)
  =
  \boldsymbol{\theta}_u(y_1)
  \boldsymbol{\theta}_u\left(\frac{z_1}{qt}\right)
  =
  \frac{1}{qt}(1+uy_1)^{-1}y_1z_1.
\]

For~\eqref{eqThetaInvariant}, we have
\begin{align*}
  \boldsymbol{\theta}_u
  \left(\frac{qt-z_1}{qt}y_1\right)
  ={}&
  \frac{1}{qt}
  \left(
    1-\bigl(qt+u(qt-z_1)y_1\bigr)^{-1}z_1
  \right)
  (1+uy_1)^{-1}y_1
  \bigl(qt+u(qt-z_1)y_1\bigr)\\
  ={}&
  \frac{1}{qt}
  \bigl(qt+u(qt-z_1)y_1\bigr)^{-1}
  (qt-z_1)(1+uy_1)
  (1+uy_1)^{-1}y_1
  \bigl(qt+u(qt-z_1)y_1\bigr)\\
  ={}&
  \frac{1}{qt}
  \bigl(qt+u(qt-z_1)y_1\bigr)^{-1}
  (qt-z_1)y_1
  \bigl(qt+u(qt-z_1)y_1\bigr)=
  \frac{qt-z_1}{qt}y_1.
\end{align*}
The last equality holds because $qt+u(qt-z_1)y_1$ and its
formal inverse commute with $(qt-z_1)y_1$.

Finally, inverting~\eqref{eqThetaY} gives
\[
  \boldsymbol{\theta}_u(y_1^{-1})
  =
  \boldsymbol{\theta}_u(y_1)^{-1}
  =
  qt\bigl(qt+u(qt-z_1)y_1\bigr)^{-1}
  (y_1^{-1}+u).\qedhere
\]
\end{proof}

Define
\begin{equation}\label{eqdefvartheta}
  \vartheta
  :=
  \boldsymbol{\theta}_r
  \circ\boldsymbol{\theta}_s
  \circ\boldsymbol{\theta}_{-ut}^{-1},\qquad \vartheta'
  :=
  \boldsymbol{\theta}'_r
  \circ\boldsymbol{\theta}'_s
  \circ\bigl(\boldsymbol{\theta}'_{-ut}\bigr)^{-1}.
\end{equation}
Theorem~\ref{DIIPthm}(2) gives
\begin{align*}
  \vartheta'(LF)
  =
  \boldsymbol{\theta}'_r\boldsymbol{\theta}'_s
  \left(
    \boldsymbol{\theta}_{-ut}^{-1}(L)
    \bigl(\boldsymbol{\theta}'_{-ut}\bigr)^{-1}(F)
  \right)=
  \bigl(
    \boldsymbol{\theta}_r
    \boldsymbol{\theta}_s
    \boldsymbol{\theta}_{-ut}^{-1}
  \bigr)(L)\,
  \bigl(
    \boldsymbol{\theta}'_r
    \boldsymbol{\theta}'_s
    \bigl(\boldsymbol{\theta}'_{-ut}\bigr)^{-1}
  \bigr)(F)=\vartheta(L)\vartheta'(F),
\end{align*}
for $L\in \mathbb{A}_{q,t}$ and $F\in \bigoplus_{j\geq 0}V_j[[u,r,s]]$. 

\begin{corollary}\label{corThetaComposition}
For $\vartheta$ defined in~\eqref{eqdefvartheta}, we have
\begin{equation}\label{eqRefinedS13VarthetaZ}
  \vartheta\left(\frac{z_1}{qt}\right)
  =
  \bigl(qt+r(qt-z_1)y_1\bigr)^{-1}
  \bigl(qt+s(qt-z_1)y_1\bigr)^{-1}
  \bigl(qt-ut(qt-z_1)y_1\bigr)z_1.
\end{equation}
\begin{equation}\label{eqRefinedS13VarthetaDPlus}
  \vartheta(d_+^*)
  =
  qt\bigl(qt+r(qt-z_1)y_1\bigr)^{-1}
  \bigl(qt+s(qt-z_1)y_1\bigr)^{-1}
  \bigl(qt-ut(qt-z_1)y_1\bigr)d_+^*.
\end{equation}
\begin{equation}\label{eqRefinedS13VarthetaYInverse}
  \vartheta(y_1^{-1})
  =
  qt\bigl(qt+r(qt-z_1)y_1\bigr)^{-1}
  \bigl(qt+s(qt-z_1)y_1\bigr)^{-1}
  \bigl(qt(y_1^{-1}+r+s)+rs(qt-z_1)y_1+ut z_1\bigr).
\end{equation}
\end{corollary}

\begin{proof}
By~\eqref{eqThetaInvariant}, every $\boldsymbol{\theta}_u$
and its formal inverse fix $(qt-z_1)y_1$.
The factors involving only this element commute with one another,
as do their formal inverses. This justifies their reordering below;
no such factor is moved past $z_1$ or $y_1^{-1}$.

For~\eqref{eqRefinedS13VarthetaZ}, equation~\eqref{eqRefinedS13ThetaZ} gives
\[
  \boldsymbol{\theta}_{-ut}^{-1}\left(\frac{z_1}{qt}\right)
  =
  \left(1-ut\frac{qt-z_1}{qt}y_1\right)\frac{z_1}{qt}.
\]
Applying $\boldsymbol{\theta}_s$ and then
$\boldsymbol{\theta}_r$ yields
\begin{align*}
  \vartheta\left(\frac{z_1}{qt}\right)
  &=
  \bigl(qt-ut(qt-z_1)y_1\bigr)
  \bigl(qt+s(qt-z_1)y_1\bigr)^{-1}
  \bigl(qt+r(qt-z_1)y_1\bigr)^{-1}z_1\\
  &=
  \bigl(qt+r(qt-z_1)y_1\bigr)^{-1}
  \bigl(qt+s(qt-z_1)y_1\bigr)^{-1}
  \bigl(qt-ut(qt-z_1)y_1\bigr)z_1.
\end{align*}

For~\eqref{eqRefinedS13VarthetaDPlus}, Theorem~\ref{DIIPthm}\textup{(1)} and~\eqref{eqThetaInvariant} give
\begin{align*}
  \vartheta(d_+^*)
  =
  \boldsymbol{\theta}_r\boldsymbol{\theta}_s
  \left(
    \left(1-ut\frac{qt-z_1}{qt}y_1\right)d_+^*
  \right)
  =
  qt\bigl(qt+r(qt-z_1)y_1\bigr)^{-1}
  \bigl(qt+s(qt-z_1)y_1\bigr)^{-1}
  \bigl(qt-ut(qt-z_1)y_1\bigr)d_+^*.
\end{align*}

For~\eqref{eqRefinedS13VarthetaYInverse}, equations~\eqref{eqRefinedS13ThetaZ} and~\eqref{eqRefinedS13ThetaYInverse} give
\begin{align*}
  \boldsymbol{\theta}_{-ut}
  \left(y_1^{-1}+ut\frac{z_1}{qt}\right)
  &=
  \bigl(qt-ut(qt-z_1)y_1\bigr)^{-1}
  \bigl(qt(y_1^{-1}-ut)+ut z_1\bigr)\\
  &=
  \bigl(qt-ut(qt-z_1)y_1\bigr)^{-1}
  \bigl(qt-ut(qt-z_1)y_1\bigr)y_1^{-1}
  =
  y_1^{-1}.
\end{align*}
Thus
$\boldsymbol{\theta}_{-ut}^{-1}(y_1^{-1})
=y_1^{-1}+ut\,z_1/(qt)$.
Applying $\boldsymbol{\theta}_s$ and then
$\boldsymbol{\theta}_r$, we obtain
\begin{align*}
  \vartheta(y_1^{-1})
  ={}&
  \boldsymbol{\theta}_r
  \left(
    \bigl(qt+s(qt-z_1)y_1\bigr)^{-1}
    \bigl(qt(y_1^{-1}+s)+ut z_1\bigr)
  \right)\\
  ={}&
  qt\bigl(qt+s(qt-z_1)y_1\bigr)^{-1}
  \bigl(qt+r(qt-z_1)y_1\bigr)^{-1}
  \bigl(qt(y_1^{-1}+r)+s\bigl(qt+r(qt-z_1)y_1\bigr)+ut z_1\bigr)\\
  ={}&
  qt\bigl(qt+r(qt-z_1)y_1\bigr)^{-1}
  \bigl(qt+s(qt-z_1)y_1\bigr)^{-1}
  \bigl(qt(y_1^{-1}+r+s)+rs(qt-z_1)y_1+ut z_1\bigr).\qedhere
\end{align*}
\end{proof}

We obtain the following proposition from
Corollary~\ref{corThetaComposition}.

\begin{proposition}\label{propSeedDivisibility}
For every $j\geq 1$, each coefficient of
$\vartheta'(P_j)$, viewed as a formal power series
in $u,r,s$, is a polynomial divisible by $y_1\cdots y_j$.
Moreover, it is symmetric in $y_1,\ldots,y_j$.
\end{proposition}

\begin{proof}
By Proposition~\ref{propSeed}, $y_1^{-1}P_j=-d_+^*P_{j-1}$.
Applying $\vartheta'$ and using
Corollary~\ref{corThetaComposition} yields
\begin{align}
  \left(
    y_1^{-1}+r+s
    +rs\frac{qt-z_1}{qt}y_1
    +ut\frac{z_1}{qt}
  \right)\vartheta'(P_j) =
  -\left(
    1-ut\frac{qt-z_1}{qt}y_1
  \right)d_+^*\vartheta'(P_{j-1}).
  \label{eqTransportedSeed}
\end{align}
Multiplying by $y_1$ on the left gives
\begin{align*}
  \vartheta'(P_j)
  =
  -y_1\left(
    1-ut\frac{qt-z_1}{qt}y_1
  \right)d_+^*\vartheta'(P_{j-1})
  -y_1\left(
    r+s+rs\frac{qt-z_1}{qt}y_1
    +ut\frac{z_1}{qt}
  \right)\vartheta'(P_j).
\end{align*}
$\vartheta'(P_j)$ is coefficientwise divisible by $y_1$
because the operators on the right preserve polynomials.

Since $\vartheta(T_i)=T_i$ and $T_iP_j=P_j$, we also have
$T_i\vartheta'(P_j)
  =
  \vartheta'(T_iP_j)=
  \vartheta'(P_j)$.
Equation~\eqref{eqTi} therefore implies
\begin{align*}
  0
  =
  (y_{i+1}-y_i)
  \bigl(T_i\vartheta'(P_j)-\vartheta'(P_j)\bigr)=
  (y_{i+1}-qy_i)
  \bigl(s_i\vartheta'(P_j)-\vartheta'(P_j)\bigr).
\end{align*}
Consequently, $\vartheta'(P_j)$ is symmetric in
$y_1,\ldots,y_j$.  Divisibility by $y_1$ then implies divisibility
by every $y_i$, and hence by $y_1\cdots y_j$.
\end{proof}

\begin{definition}\label{defRefinedS13FjDefinition}
Define
$f_j[X;y_1,\ldots,y_j]:=\vartheta'(P_j)/P_j$, with $f_0:=1$.
\end{definition}
Each coefficient of $f_j$ is a polynomial symmetric in $y_1,\ldots,y_j$. 
Definition~\ref{defRefinedS13FjDefinition}, Proposition~\ref{propSeed},
Theorem~\ref{DIIPthm}\textup{(2)--(3)}, and~\eqref{eqdefvartheta} give
\begin{align}\label{eqFjThetaIdentity}
  d_-^j\bigl(P_j f_j\bigr)
  =
  d_-^j\vartheta'(P_j)
  =
  \vartheta'(d_-^jP_j)
  =
  \Theta(r)\Theta(s)\Theta(-ut)^{-1}
  \widetilde H_{(j)}.
\end{align}

The functions $f_j$ satisfy the following recurrence.
\begin{theorem}\label{thmFjRecursion}
For every $j\geq 1$, we have
\begin{equation}\label{eqFjRecursion}
\begin{aligned}
  &(1+ry_1)(1+sy_1)f_j[X;y_1,\ldots,y_j] =
  f_{j-1}[X+(q-1)ty_1;y_2,\ldots,y_j]\\
  &\qquad{}
  +y_1\langle y_{j+1}^0\rangle
  \left(
    \frac{\Omega[-X/y_{j+1}]}{1-ty_1/y_{j+1}}
    (rs\,y_{j+1}-ut)\,
    f_j\!\left[
      X+(1-q)(y_{j+1}-ty_1);
      y_2,\ldots,y_j,y_{j+1}
    \right]
  \right).
\end{aligned}
\end{equation}
\end{theorem}

\begin{proof}
Since $P_jf_j$ is symmetric in $y_1,\ldots,y_j$, equation~\eqref{eqTi}
gives $T_i(P_jf_j)=P_jf_j$ and
$T_i(y_{i+1}P_jf_j)=qy_iP_jf_j$.  Hence
\begin{equation}\label{eqRefinedS13HeckeSymmetry}
  T_{j-1}^{-1}\cdots T_1^{-1}(P_jf_j)=P_jf_j,
  \qquad
  T_{j-1}^{-1}\cdots T_1^{-1}(y_1P_jf_j)
  =q^{1-j}y_jP_jf_j.
\end{equation}
We first compute the two compositions in~\eqref{eqzi}.
\begin{alignat}{2}
  &\quad d_+^*d_-(P_jf_j)\notag\\*
  &=d_+^*\Bigl(\langle y_j^0\rangle
    \bigl(P_j f_j[X+(1-q)y_j;y_1,\ldots,y_j]
      \Omega[-X/y_j]\bigr)\Bigr)
  &\qquad\text{(by~\eqref{eqdminus})}&\notag\\
  &={}\frac{P_j}{y_1}\langle y_{j+1}^0\rangle
    \Bigl(y_{j+1}\Omega[-(X+(q-1)ty_1)/y_{j+1}]\,f_j[X+(q-1)ty_1+(1-q)y_{j+1};y_2,\ldots,y_j,y_{j+1}]\Bigr)
  &\qquad\text{(by~\eqref{eqdplusstar}, \eqref{eqRefinedS13PDefinition})}&\notag\\
  &={}\frac{P_j}{y_1}\langle y_{j+1}^0\rangle
    \Bigl(y_{j+1}\Omega[-X/y_{j+1}]
      \Omega[(1-q)ty_1/y_{j+1}]\,f_j[X+(1-q)(y_{j+1}-ty_1);y_2,\ldots,y_j,y_{j+1}]\Bigr)
  &\qquad\text{(by~\eqref{eqOmega})}&\rlap{,}\notag\\
  &\quad d_-d_+^*(P_jf_j)\notag\\*
  &=d_-\Bigl(\frac{P_j}{y_1}y_{j+1}
    f_j[X+(q-1)ty_1;y_2,\ldots,y_j,y_{j+1}]\Bigr)
  &\qquad\text{(by~\eqref{eqdplusstar}, \eqref{eqRefinedS13PDefinition})}&\notag\\
  &={}\frac{P_j}{y_1}\langle y_{j+1}^0\rangle
    \Bigl(y_{j+1}\Omega[-X/y_{j+1}]\,f_j[X+(1-q)(y_{j+1}-ty_1);y_2,\ldots,y_j,y_{j+1}]\Bigr)
  &\qquad\text{(by~\eqref{eqdminus})}&\rlap{.}
  \label{eqSeedCompositions}
\end{alignat}
Here $P_j/y_1$ is a polynomial.  In the second composition,
$d_-$ replaces $X$ by $X+(1-q)y_{j+1}$ after applying
$d_+^*$ and then multiplies by $\Omega[-X/y_{j+1}]$;
there is no subsequent substitution in this factor.
By~\eqref{eqOmega},
\begin{equation}\label{eqSeedOmegaDifference}
  \Omega[(1-q)ty_1/y_{j+1}]-1
  =\frac{1-qty_1/y_{j+1}}{1-ty_1/y_{j+1}}-1
  =\frac{(1-q)ty_1/y_{j+1}}{1-ty_1/y_{j+1}}.
\end{equation}
Subtracting the two compositions and using this difference gives
\begin{alignat}{2}
  &\quad q^{1-j}\frac{z_1}{qt}(P_jf_j)\notag\\*
  &=\frac{1}{t(1-q)}(d_+^*d_--d_-d_+^*)(P_jf_j)
  &\qquad\makebox[0pt][r]{\text{(by~\eqref{eqzi}, \eqref{eqRefinedS13HeckeSymmetry})}}&\notag\\
  &=\frac{P_j}{t(1-q)y_1}\langle y_{j+1}^0\rangle
    \Bigl(y_{j+1}\Omega[-X/y_{j+1}]
      \bigl(\Omega[(1-q)ty_1/y_{j+1}]-1\bigr)
      f_j[X+(1-q)(y_{j+1}-ty_1);y_2,\ldots,y_j,y_{j+1}]\Bigr)
  &\qquad\text{(by~\eqref{eqSeedCompositions})}&\notag\\
  &=P_j\langle y_{j+1}^0\rangle
    \Biggl(\frac{\Omega[-X/y_{j+1}]}{1-ty_1/y_{j+1}}\,
      f_j[X+(1-q)(y_{j+1}-ty_1);y_2,\ldots,y_j,y_{j+1}]\Biggr)
  &\qquad\text{(by~\eqref{eqSeedOmegaDifference})}&\rlap{.}
  \label{eqZSeedAction}
\end{alignat}
The same calculation gives
\begin{alignat}{2}
  \frac{z_1}{qt}(y_1P_jf_j)&=\frac{1}{t(1-q)}(d_+^*d_--d_-d_+^*)(y_jP_jf_j)
  &\qquad\text{(by~\eqref{eqzi}, \eqref{eqRefinedS13HeckeSymmetry})}&\notag\\
  &=P_j\langle y_{j+1}^0\rangle
    \Biggl(\frac{y_{j+1}\Omega[-X/y_{j+1}]}{1-ty_1/y_{j+1}}\,f_j[X+(1-q)(y_{j+1}-ty_1);y_2,\ldots,y_j,y_{j+1}]\Biggr)
  &\qquad\text{(by~\eqref{eqSeedCompositions}, \eqref{eqOmega})}&\rlap{.}
  \label{eqZYSeedAction}
\end{alignat}
By Definition~\ref{defRefinedS13FjDefinition},
$\vartheta'(P_j)=P_jf_j$.
The factors involving $r$ and $s$ commute, since they are
polynomials in $(qt-z_1)y_1$.
Thus
\begin{alignat}{2}
  q^{1-j}
    \left(1-ut\frac{qt-z_1}{qt}y_1\right)
    \frac{z_1}{qt}P_jf_j
  &=q^{1-j}
    \left(1+r\frac{qt-z_1}{qt}y_1\right)
    \left(1+s\frac{qt-z_1}{qt}y_1\right)
    \vartheta\left(\frac{z_1}{qt}\right)(P_jf_j)
  &\quad\text{(by~\eqref{eqRefinedS13VarthetaZ})}&\notag\\
  &=q^{1-j}
    \left(1+r\frac{qt-z_1}{qt}y_1\right)
    \left(1+s\frac{qt-z_1}{qt}y_1\right)
    \vartheta'\left(\frac{z_1}{qt}P_j\right)
  &\quad\text{(by Theorem~\ref{DIIPthm}\textup{(2)},
                   \eqref{eqdefvartheta})}&\notag\\
  &=\left(1+r\frac{qt-z_1}{qt}y_1\right)
    \left(1+s\frac{qt-z_1}{qt}y_1\right)P_jf_j
  &\quad\text{(by~\eqref{eqZSeedAction} with $f_j=1$)}&
    \rlap{.}
    \label{eqTransportedSeedEigenvalue}
\end{alignat}

Next, expanding the product and cancelling the adjacent
factors $y_1y_1^{-1}$ gives
\begin{align}
  \frac{qt-z_1}{qt}y_1
    \left(
      y_1^{-1}+r+s+rs\frac{qt-z_1}{qt}y_1
    \right)
  &=1-\frac{z_1}{qt}
    +(r+s)\frac{qt-z_1}{qt}y_1
    +rs\left(\frac{qt-z_1}{qt}y_1\right)^2
    \notag\\
  &=\left(1+r\frac{qt-z_1}{qt}y_1\right)
    \left(1+s\frac{qt-z_1}{qt}y_1\right)
    -\frac{z_1}{qt}.
    \label{eqRefinedS13SeedProductExpansion}
\end{align}
Keeping $z_1$ on the right in the correction
term, we then have
\begin{alignat}{2}
  &\quad \left(1-ut\frac{qt-z_1}{qt}y_1\right)
    \left(y_1^{-1}+r+s+rs\frac{qt-z_1}{qt}y_1
          +utq^{1-j}\frac{z_1}{qt}\right)P_jf_j\notag\\
  &=\left(y_1^{-1}+r+s+rs\frac{qt-z_1}{qt}y_1
          +ut\frac{z_1}{qt}\right)P_jf_j\notag\\
  &\quad+ut\Biggl(q^{1-j}
       \left(1-ut\frac{qt-z_1}{qt}y_1\right)\frac{z_1}{qt}
       -\left(1+r\frac{qt-z_1}{qt}y_1\right)
        \left(1+s\frac{qt-z_1}{qt}y_1\right)\Biggr)P_jf_j
  &\quad\text{(by~\eqref{eqRefinedS13SeedProductExpansion})}&\notag\\
  &=\left(y_1^{-1}+r+s+rs\frac{qt-z_1}{qt}y_1
          +ut\frac{z_1}{qt}\right)P_jf_j
  &\quad\text{(by~\eqref{eqTransportedSeedEigenvalue})}&\notag\\
  &=-\left(1-ut\frac{qt-z_1}{qt}y_1\right)
       d_+^*(P_{j-1}f_{j-1})
  &\quad\text{(by~\eqref{eqTransportedSeed})}&\rlap{.}
  \label{eqSeedReduction}
\end{alignat}

The common left factor
$1-ut\frac{qt-z_1}{qt}y_1$ in~\eqref{eqSeedReduction}
has constant term $1$ in $u$ and is therefore invertible.
Comparing the first and last expressions in
\eqref{eqSeedReduction} and cancelling this factor gives
\[
  \left(
    y_1^{-1}+r+s+rs\frac{qt-z_1}{qt}y_1
    +utq^{1-j}\frac{z_1}{qt}
  \right)P_jf_j
  =-d_+^*(P_{j-1}f_{j-1}).
\]
Multiply by $y_1$ on the left and use
$y_1\frac{qt-z_1}{qt}y_1
 =y_1^2-y_1\frac{z_1}{qt}y_1$.
Moving the two terms involving $z_1$ to the right, we obtain
\[
  (1+ry_1)(1+sy_1)P_jf_j
  =-y_1d_+^*(P_{j-1}f_{j-1})
    +rs\,y_1\frac{z_1}{qt}(y_1P_jf_j)
    -utq^{1-j}y_1\frac{z_1}{qt}(P_jf_j).
\]
By~\eqref{eqdplusstar} and~\eqref{eqRefinedS13PDefinition},
$-y_1d_+^*(P_{j-1}f_{j-1})
=P_jf_{j-1}[X+(q-1)ty_1;y_2,\ldots,y_j]$.
For the remaining terms,
\begin{alignat*}{2}
  rs\,y_1\frac{z_1}{qt}(y_1P_jf_j)
  &=rs\,y_1P_j\langle y_{j+1}^0\rangle
    \Biggl(
      \frac{y_{j+1}\Omega[-X/y_{j+1}]}{1-ty_1/y_{j+1}}\,
      f_j[X+(1-q)(y_{j+1}-ty_1);y_2,\ldots,y_j,y_{j+1}]
    \Biggr)
  &\quad\text{(by~\eqref{eqZYSeedAction})}&\rlap{,}\\[4pt]
  -utq^{1-j}y_1\frac{z_1}{qt}(P_jf_j)
  &=-ut\,y_1P_j\langle y_{j+1}^0\rangle
    \Biggl(
      \frac{\Omega[-X/y_{j+1}]}{1-ty_1/y_{j+1}}\,
      f_j[X+(1-q)(y_{j+1}-ty_1);y_2,\ldots,y_j,y_{j+1}]
    \Biggr)
  &\quad\text{(by~\eqref{eqZSeedAction})}&\rlap{.}
\end{alignat*}
Substituting these expressions and using
$P_j=(-1)^jy_1\cdots y_j\neq0$, we obtain
\begin{align*}
  &(1+ry_1)(1+sy_1)f_j[X;y_1,\ldots,y_j] =
  f_{j-1}[X+(q-1)ty_1;y_2,\ldots,y_j]\\
  &\qquad{}
  +y_1\langle y_{j+1}^0\rangle
  \left(
    \frac{\Omega[-X/y_{j+1}]}{1-ty_1/y_{j+1}}
    (rs\,y_{j+1}-ut)\,
    f_j\!\left[
      X+(1-q)(y_{j+1}-ty_1);
      y_2,\ldots,y_j,y_{j+1}
    \right]
  \right),
\end{align*}
as required.
\end{proof}

\subsection{The algebraic recurrence}\label{subsecAlgebraicRecurrence}
Let $n\geq 0$ and consider the finite alphabet $X_n$.
For $I\subseteq X_n$, define
\begin{equation}\label{eqDefGAlg}
  \mathcal L_I(X_n;ut,r,s)
  :=
  f_{|I|}[(q-1)(X_n\setminus I);I].
\end{equation}
In particular, $\mathcal L_{\varnothing}(X_n;ut,r,s)=1$
since $f_0=1$ by Definition~\ref{defRefinedS13FjDefinition}.
This evaluation is independent of the order of the letters of $I$
because $f_{|I|}$ is symmetric in $y_1,\ldots,y_{|I|}$.

For $x\in X_n$, write $X_n^{x\to tx}$ for the alphabet obtained
by replacing $x$ with $tx$.
Its letters are again algebraically independent over $\mathbb{Q}(q,t)$,
so $X_n^{x\to tx}$ is also a finite alphabet.

\begin{proposition}\label{propAlgRecursion}
For every $I\subseteq X_n$ and $x\in I$, we have
\begin{equation}\label{eqAlgRecursion}
\begin{aligned}
  &(1+rx)(1+sx)\mathcal L_I(X_n;ut,r,s)
  =
  \mathcal L_{I\setminus\{x\}}
  (X_n^{x\to tx};ut,r,s)\\[-2pt]
  &\quad{}
  +tx(rsx-u)
  \prod_{x'\in X_n\setminus I}\frac{tx-qx'}{tx-x'}\,
  \mathcal L_{(I\setminus\{x\})\cup\{tx\}}
  (X_n^{x\to tx};ut,r,s)\\[-2pt]
  &\quad{}
  +(1-q)x
  \sum_{x'\in X_n\setminus I}
  \frac{x'\eta_{x'}(X_n\setminus I)}{x'-tx}
  (rsx'-ut)\,
  \mathcal L_{(I\setminus\{x\})\cup\{x'\}}
  (X_n^{x\to tx};ut,r,s).
\end{aligned}
\end{equation}
\end{proposition}

\begin{proof}
Set $j=|I|$.  Specializing~\eqref{eqFjRecursion} at
$X\mapsto(q-1)(X_n\setminus I)$ and
$(y_1,\ldots,y_j)=(x,I\setminus\{x\})$, and
using~\eqref{eqDefGAlg}, gives
\begin{multline*}
  (1+rx)(1+sx)\mathcal L_I(X_n;ut,r,s)
  =
  \mathcal L_{I\setminus\{x\}}(X_n^{x\to tx};ut,r,s)\\
  +
  x\langle y_{j+1}^0\rangle
  \Biggl(
    \frac{y_{j+1}(rs\,y_{j+1}-ut)}{y_{j+1}-tx}
    \prod_{x'\in X_n\setminus I}
    \frac{y_{j+1}-qx'}{y_{j+1}-x'}
    f_j[(q-1)((X_n\setminus I)-y_{j+1}+tx);
        I\setminus\{x\},y_{j+1}]
  \Biggr).
\end{multline*}
We use the following identity to separate the factors in the
denominator:
\begin{equation}\label{eqRefinedS13PartialFractionShift}
  \frac{y_{j+1}}{y_{j+1}-tx}
  \prod_{x'\in X_n\setminus I}
  \frac{y_{j+1}-qx'}{y_{j+1}-x'}
  =
  \prod_{x'\in X_n\setminus I}\frac{tx-qx'}{tx-x'}
  \frac{y_{j+1}}{y_{j+1}-tx}
  +
  (1-q)\sum_{x'\in X_n\setminus I}
  \frac{x'\eta_{x'}(X_n\setminus I)}{x'-tx}
  \frac{y_{j+1}}{y_{j+1}-x'}.
\end{equation}
To verify~\eqref{eqRefinedS13PartialFractionShift}, multiply both sides by
$(y_{j+1}-tx)\prod_{x'\in X_n\setminus I}(y_{j+1}-x')$.
The resulting polynomials have degree at most $n-j+1$ in
$y_{j+1}$ and agree at the $n-j+2$ distinct values
$0$, $tx$, and $x'\in X_n\setminus I$.
They are therefore equal.

Expanding $y_{j+1}/(y_{j+1}-tx)$ as
$\sum_{h\geq0}(tx)^h y_{j+1}^{-h}$, we see that
$y_{j+1}^m$ contributes to the constant term only through
$(tx)^m y_{j+1}^{-m}$.  Hence
\(
  \langle y_{j+1}^0\rangle
  \left(
    \frac{y_{j+1}}{y_{j+1}-tx}\,y_{j+1}^m
  \right)
  =(tx)^m,
\)
for $m\geq0$.
Thus, when the fraction multiplies a polynomial, taking the
constant term amounts to evaluating that polynomial at $tx$.

For the fraction $y_{j+1}/(y_{j+1}-tx)$, substitute
$y_{j+1}=tx$.  The symmetric-function alphabet in $f_j$ becomes
$(q-1)(X_n\setminus I)$, and its polynomial variable becomes
$tx$.  Thus~\eqref{eqDefGAlg} gives
\[
  f_j[(q-1)(X_n\setminus I);I\setminus\{x\},tx]
  =
  \mathcal L_{(I\setminus\{x\})\cup\{tx\}}
  (X_n^{x\to tx};ut,r,s).
\]
The factor $rs\,y_{j+1}-ut$ becomes $rs\,tx-ut$.
Multiplying by the outer factor $x$ and the coefficient
$\prod_{x'\in X_n\setminus I}(tx-qx')/(tx-x')$ gives the second term
of~\eqref{eqAlgRecursion}.

For each $x'\in X_n\setminus I$, the fraction
$y_{j+1}/(y_{j+1}-x')$ instead substitutes $y_{j+1}=x'$.
The symmetric-function alphabet in $f_j$ becomes
$(q-1)((X_n\setminus I)-x'+tx)$, so~\eqref{eqDefGAlg} gives
\[
  f_j[(q-1)((X_n\setminus I)-x'+tx);
      I\setminus\{x\},x']
  =
  \mathcal L_{(I\setminus\{x\})\cup\{x'\}}
  (X_n^{x\to tx};ut,r,s).
\]
Here $rs\,y_{j+1}-ut$ becomes $rsx'-ut$.
Multiplying by $x$ and the coefficient
$(1-q)x'\eta_{x'}(X_n\setminus I)/(x'-tx)$,
then summing over $x'\in X_n\setminus I$, gives the third term
of~\eqref{eqAlgRecursion}.
\end{proof}

We now state the main recursion on the algebraic side.
\begin{theorem}\label{thmAlgMainRecursion}
For every $n\geq 0$, we have
\begin{equation*}
  L[(q-1)X_n;u,v,r,s]
  =
  \sum_{j=1}^{n}(-uv)^j(q;q)_j
  \sum_{\substack{I\subseteq X_n\\|I|=j}}
  \left(\prod_{x\in I}x\right)
  \eta_I(X_n)\,
  \mathcal L_I(X_n;ut,r,s),
\end{equation*}
where the functions $\mathcal L_I$ satisfy
the recurrence~\eqref{eqAlgRecursion} and the initial condition
$\mathcal L_{\varnothing}(X_n;ut,r,s)=1$.
These conditions uniquely determine the family
$\mathcal L_I$ as formal power series
in $u,r,s$ with coefficients rational in the letters of $X_n$.
\end{theorem}

\begin{proof}
By~\eqref{eqThetaalg}, Proposition~\ref{propHag04},
\eqref{eqFjThetaIdentity}, Corollary~\ref{corDminusFiniteAlphabet},
and~\eqref{eqDefGAlg}, we obtain
\begin{align*}
  L[(q-1)X_n;u,v,r,s]
  &=
  \sum_{j\geq 1}(uv)^j
  \bigl(
    \Theta(r)\Theta(s)\Theta(-ut)^{-1}
    \widetilde H_{(j)}
  \bigr)[(q-1)X_n]\\
  &=
  \sum_{j\geq 1}(uv)^j
  \bigl(
    d_-^j(P_j f_j)
  \bigr)[(q-1)X_n]\\
  &=
  \sum_{j=1}^{n}(-uv)^j(q;q)_j
  \sum_{\substack{I\subseteq X_n\\|I|=j}}
  \left(\prod_{x\in I}x\right)
  \eta_I(X_n)\,
  f_j[(q-1)(X_n\setminus I);I]\\
  &=
  \sum_{j=1}^{n}(-uv)^j(q;q)_j
  \sum_{\substack{I\subseteq X_n\\|I|=j}}
  \left(\prod_{x\in I}x\right)
  \eta_I(X_n)\,
  \mathcal L_I(X_n;ut,r,s).
\end{align*}

The initial condition $\mathcal L_\varnothing=1$ determines
all coefficients for $I=\varnothing$.
We determine the remaining coefficients by induction on their
total degree in $u,r,s$, simultaneously for all finite alphabets.

In degree zero, \eqref{eqAlgRecursion} reduces to
$\mathcal L_I(X_n;0,0,0)
  =
  \mathcal L_{I\setminus\{x\}}(X_n^{x\to tx};0,0,0)$.
Repeatedly removing a letter from $I$ therefore gives
$\mathcal L_I(X_n;0,0,0)=1$.

Now fix $d\geq1$ and suppose that all coefficients of total
degree less than $d$ have been determined.
On the left-hand side of~\eqref{eqAlgRecursion}, multiplication
by $(1+rx)(1+sx)=1+(r+s)x+rsx^2$ leaves only the degree-$d$
coefficients of $\mathcal L_I$ to be determined.
On the right-hand side, the last two terms involve only
coefficients of degree less than $d$, because their prefactors
have positive degree.
Thus the degree-$d$ coefficients of $\mathcal L_I$ are determined
by those of
$\mathcal L_{I\setminus\{x\}}(X_n^{x\to tx};ut,r,s)$.

Starting from $I=\varnothing$, induction on $|I|$ determines
all coefficients of degree $d$.
Since this induction applies to every finite alphabet, it also
applies to $X_n^{x\to tx}$.
Induction on $d$ proves uniqueness.
\end{proof}

\section{Statistics under main-diagonal deletion}\label{secDiagonalDeletion}
In this section, we study deletion of the main diagonal of a
signed parking function $\pi$.
We group the roots of the resulting path $\pi'$ into blocks
and record the deletion data $(\pi';W,\delta,h,\chi)$.
We characterize the admissible deletion data and prove that
they determine $\pi$ uniquely in
Proposition~\ref{propMainDiagonalDeletion}.
We then compute the changes in $\area$ and $\dinv$ in
Theorem~\ref{thmMainDiagonalStatistics}, obtaining the
insertion weights used in
Sections~\ref{secInsertionsWithoutRise}
and~\ref{secInsertionsWithRise}.

\subsection{Deleting the main diagonal}\label{subsecDeletingMainDiagonal}
\begin{definition}\label{defRefinedS45DeletionData}
Let $\pi\in\PF_{d+k+\ell,\pm}^{\ast k,\bullet\ell}(X_n)$ be a path with $L\geq1$ roots.
Delete the roots and subtract $1$ from each remaining area entry.
Remove the rise decoration from the first child of each deleted root, if present.
Retain the other labels and decorations, and reindex the north steps.
Write $\pi'$ for the resulting path.

Write $W=(c_1,\ldots,c_L)$ for the root labels of $\pi$, in order.
If a root is decorated, its decoration must be a valley decoration.
Let $\delta_j=1$ if the $j$th root is a decorated valley, and $\delta_j=0$ otherwise.

A \emph{block} consists of a root outside $\DD(\pi')$ together
with all consecutive roots in $\DD(\pi')$ immediately following
it in path order.
The first root of a block is its \emph{head}.

Suppose that there are $H$ blocks, with head labels
$\beta=(\beta_1,\ldots,\beta_H)$.
Let $h_j$ be the number of blocks \emph{assigned} to the $j$th deleted root, that is, whose heads are children of this root in $\pi$.  Then
\[
  h_j\geq0,\qquad h_1+\cdots+h_L=H.
\]
Finally, let $\chi_j=1$ if the first child of the $j$th deleted root carries a rise decoration in $\pi$, and let $\chi_j=0$ otherwise.
In particular, $\chi_j=1$ implies $h_j>0$.
The tuple $(\pi';W,\delta,h,\chi)$ is the \emph{deletion data} of $\pi$.
\end{definition}

\begin{example}\label{exDeletingMainDiagonal}
Recall the signed path $\pi$ from Example~\ref{exSignedDecoratedPath}.
Delete its roots at steps $1,5,9,10$ and subtract $1$ from each remaining area entry.
Also remove the decorated rise on step $6$, the first child of the root at step $5$.
The resulting path $\pi'$ is shown in Figure~\ref{figMainDiagonalDeletion}.

\begin{figure}[!htbp]
  \centering
  \begin{tikzpicture}[x=.5cm,y=.5cm,font=\small]
    \begin{scope}
      \draw[step=1,blue!20,thin] (0,0) grid (10,10);
      \draw[dashed,blue!45] (0,0)--(10,10);
      \draw[very thick,blue!75!black]
        (0,0)--(0,3)--(2,3)--(2,4)--(4,4)--(4,6)
        --(5,6)--(5,8)--(8,8)--(8,9)--(9,9)
        --(9,10)--(10,10);

      \foreach \x/\y/\w in {
        .5/.5/{1^-},
        .5/1.5/{1^-},
        .5/2.5/{7^-},
        2.5/3.5/{4^+},
        4.5/4.5/{2^+},
        4.5/5.5/{3^+},
        5.5/6.5/{5^+},
        5.5/7.5/{7^-},
        8.5/8.5/{8^-},
        9.5/9.5/{6^+}
      }
        \node at (\x,\y) {$\w$};

      \node at (-.3,2.5) {$\ast$};
      \node at (3.7,5.5) {$\ast$};
      \node at (4.7,7.5) {$\ast$};
      \node at (1.7,3.5) {$\bullet$};
      \node at (4.7,6.5) {$\bullet$};
      \node at (7.7,8.5) {$\bullet$};
      \node at (5,-.7) {$\pi$};
    \end{scope}

    \begin{scope}[xshift=6.2cm]
      \draw[step=1,blue!20,thin] (0,0) grid (6,6);
      \draw[dashed,blue!45] (0,0)--(6,6);
      \draw[very thick,blue!75!black]
        (0,0)--(0,2)--(2,2)--(2,3)--(3,3)
        --(3,4)--(4,4)--(4,6)--(6,6);

      \foreach \x/\y/\w in {
        .5/.5/{1^-},
        .5/1.5/{7^-},
        2.5/2.5/{4^+},
        3.5/3.5/{3^+},
        4.5/4.5/{5^+},
        4.5/5.5/{7^-}
      }
        \node at (\x,\y) {$\w$};

      \node at (-.3,1.5) {$\ast$};
      \node at (3.7,5.5) {$\ast$};
      \node at (1.7,2.5) {$\bullet$};
      \node at (3.7,4.5) {$\bullet$};
      \node at (3,-.7) {$\pi'$};
    \end{scope}
  \end{tikzpicture}
  \caption{Main-diagonal deletion of $\pi$ (left), giving $\pi'$ (right).}
  \label{figMainDiagonalDeletion}
\end{figure}
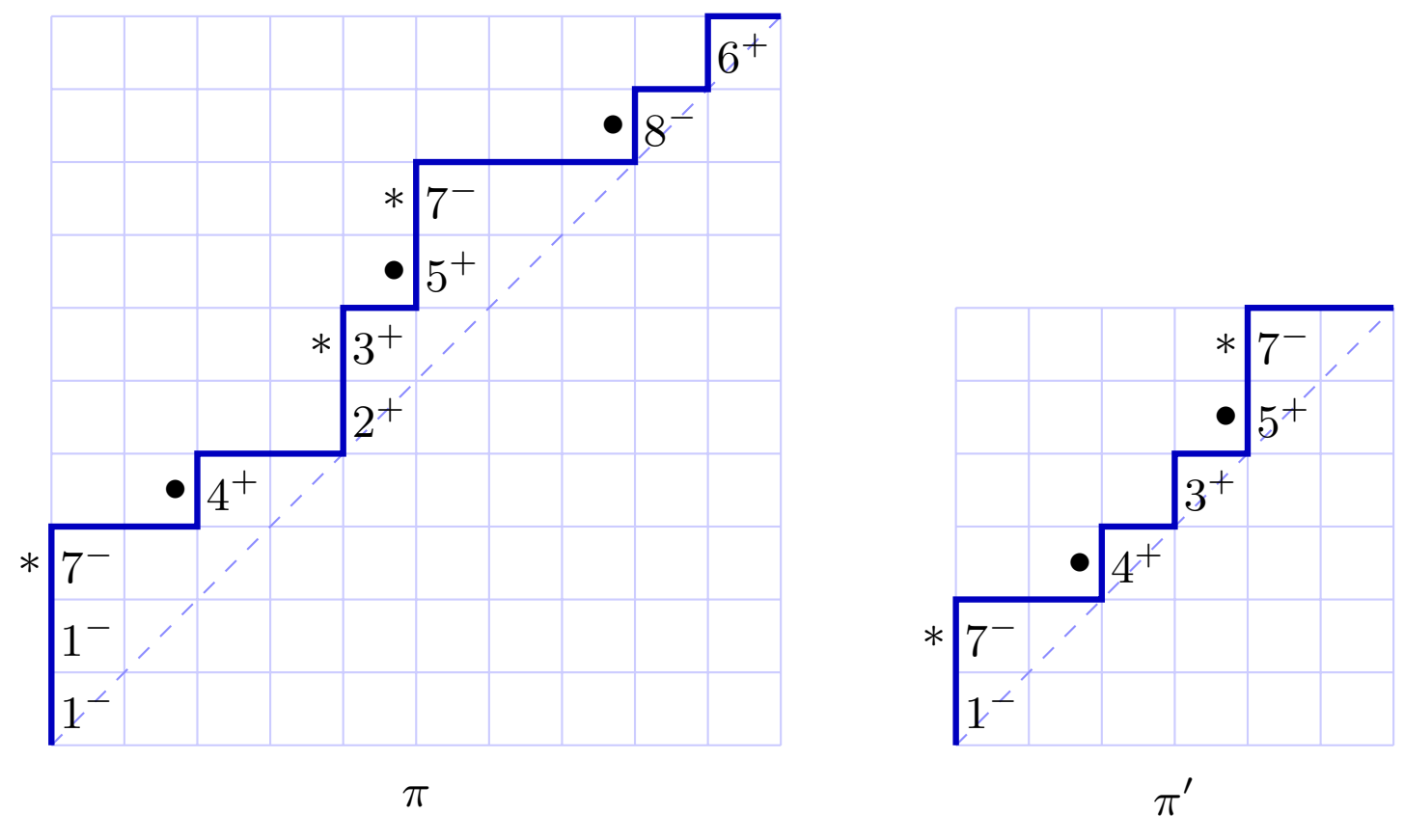

The path $\pi'$ has area sequence $(0,1,0,0,0,1)$, with
$\RR(\pi')=\{2,6\}$ and $\DD(\pi')=\{3,5\}$.
Its two blocks have root indices $(1,3)$ and $(4,5)$, respectively.
Thus the deletion data are
\begin{gather*}
  W=(1^-,2^+,8^-,6^+),\qquad
  \delta=(0,0,1,0),\qquad
  h=(1,1,0,0),\qquad
  \chi=(0,1,0,0).
\end{gather*}
Here $\delta_3=1$ records the deleted decorated valley at step $9$,
and $\chi_2=1$ records the removed decorated rise on step $6$.
\end{example}

The following proposition characterizes the deletion data.
\begin{proposition}\label{propMainDiagonalDeletion}
Let $\pi'$ be a signed parking function on $X_n$, possibly empty, whose $H$ blocks have head labels $(\beta_1,\ldots,\beta_H)$.
Fix $L\geq1$, a sequence $W=(c_1,\ldots,c_L)$ of signed labels on $X_n$, and vectors
$\delta,\chi\in\{0,1\}^L$ and $h\in\mathbb{Z}_{\geq0}^L$ with $h_1+\cdots+h_L=H$.
There is a unique signed parking function $\pi$ with deletion data $(\pi';W,\delta,h,\chi)$ if and only if
    $\textup{(1)} \delta_1=0$;
    $\textup{(2)}
    \text{If }\delta_j=1,\text{ then }
      h_{j-1}>0\text{ or }c_{j-1}\prec c_j$;
    $\textup{(3)}
    \text{If }h_j>0,\text{ then }
      c_j\prec\beta_{h_1+\cdots+h_{j-1}+1}$;
    $\textup{(4)}
    \text{If }\chi_j=1,\text{ then }h_j>0$.

Here condition~\textup{(2)} applies for $2\leq j\leq L$, and
conditions~\textup{(3)}--\textup{(4)} apply for $1\leq j\leq L$.
\end{proposition}

\begin{proof}
Suppose first that the data are obtained from $\pi$.
The first root is not a valley, which proves \textup{(1)}.
For $j\geq2$, if $h_{j-1}>0$, the $j$th root follows at least
two east steps and is therefore contractible.
If $h_{j-1}=0$, the two successive roots are separated by one
east step, so contractibility requires $c_{j-1}\prec c_j$.
This proves \textup{(2)}.
If $h_j>0$, the first child of the $j$th root has label
$\beta_{h_1+\cdots+h_{j-1}+1}$.
The column condition therefore gives \textup{(3)}.
Finally, $\chi_j=1$ means that the first child is a decorated
rise, so $h_j>0$, proving \textup{(4)}.

Conversely, for $j=1,\ldots,L$ in order, insert a root labeled
$c_j$, followed by the next $h_j$ blocks of $\pi'$ together
with the descendants of their roots.
Increase the area entries of all these appended north steps by $1$.
If $h_j=0$, the inserted root has no children.
Keep the decorations of $\pi'$, decorate the inserted root
as a valley when $\delta_j=1$, and decorate its first child
as a rise when $\chi_j=1$.

The resulting area sequence defines a Dyck path.
Each nonempty group of assigned blocks begins at a root outside
$\DD(\pi')$.
Consequently, every retained decorated valley has its preceding
north step in the same group; their area difference and labels
are unchanged, so the valley remains contractible.
Condition~\textup{(3)} gives the column comparison between each
inserted root and its first child.
Conditions~\textup{(1)}--\textup{(2)} ensure that every newly
decorated valley is contractible.
Condition~\textup{(4)} ensures that a first child exists whenever
a rise decoration is prescribed, and this child is a rise by
Lemma~\ref{lemForestStarValleys}(1).

Every step of this construction is determined by the data.
Deleting the main diagonal recovers $\pi'$ and the same
$W,\delta,h,\chi$, so the construction is the unique inverse.
\end{proof}
\begin{definition}\label{defAdmissibleDeletionData}
The deletion data are \emph{admissible} if they satisfy
conditions~\textup{(1)}--\textup{(4)} of
Proposition~\ref{propMainDiagonalDeletion}.
\end{definition}

\subsection{Changes in the statistics}\label{subsecStatisticsChanges}
Fix deletion data $(\pi';W,\delta,h,\chi)$ satisfying the
conditions of Proposition~\ref{propMainDiagonalDeletion},
and let $\pi$ be the corresponding parking function.
We define the following statistics on the deletion data.

\begin{definition}\label{defRefinedS45Statistics}
For $1\leq b\leq H$, let $o(b)$ be the index of the deleted root
to which the $b$th block is assigned.
Equivalently,
$h_1+\cdots+h_{o(b)-1}<b\leq h_1+\cdots+h_{o(b)}$.
Define
\begin{equation}\label{eqRefinedS45RootStatistics}
  \dinv_0(W)
  :=
  \sum_{\substack{1\leq j<k\leq L\\ \delta_j=0}}
  \mathbf{1}_{\{c_j\prec c_k\}}
  -\sum_{j=1}^{L}\delta_j,
  \qquad
  C(W;\delta,h;\beta)
  :=
  \sum_{b=1}^{H}\sum_{k=o(b)+1}^{L}
  \mathbf{1}_{\{c_k\prec\beta_b\}}.
\end{equation}
For $1\leq j<L$, define
\begin{equation}\label{eqRefinedS45RiseCorrection}
  \nu_j
  :=
  \delta_{j+1}\sum_{k=j+2}^{L}
  \left(
    \delta_k\mathbf{1}_{\{c_{j+1}\prec c_k\}}
    +(1-\delta_k)\mathbf{1}_{\{c_k\prec c_{j+1}\}}
  \right)
  -\sum_{k=j+2}^{L}\delta_k,
\end{equation}
and set $\nu_L=0$.
We write $|\pi'|$ for the number of north steps of $\pi'$.
\end{definition}

\begin{example}\label{exMainDiagonalCorrections}
Recall the deletion data in
Example~\ref{exDeletingMainDiagonal}.
The two blocks belong to the first and second deleted roots,
so $o(1)=1$ and $o(2)=2$.
The first two deleted roots contribute $3$ and $2$ inversion
pairs of type \textup{(1)}, respectively.
Consequently,
\[
  \dinv_0(W)=3+2-1=4,
  \qquad
  C(W;\delta,h;\beta)=(0+0+0)+(0+0)=0.
\]
We also have
$\nu_1=-1$,
$\nu_2=\mathbf{1}_{\{6^+\prec8^-\}}=1$, and
$\nu_3=\nu_4=0$.
Since $\chi=(0,1,0,0)$, we obtain
$\sum_{j=1}^{4}\chi_j\nu_j=1$,
$|\pi'|=6$,
$|\RR(\pi')|=2$, and
$\sum_{j=1}^{4}\chi_j=1$.
\end{example}

\begin{theorem}\label{thmMainDiagonalStatistics}
The parking function $\pi$ reconstructed from the deletion data
$(\pi';W,\delta,h,\chi)$ satisfies
\begin{align*}
  \dinv(\pi)-\dinv(\pi')
  =
  \dinv_0(W)+C(W;\delta,h;\beta)
  +\sum_{j=1}^{L}\chi_j\nu_j,\quad
  \area(\pi)-\area(\pi')
  =
  |\pi'|-|\RR(\pi')|-\sum_{j=1}^{L}\chi_j.
\end{align*}
\end{theorem}

\begin{proof}
First reconstruct $\pi$ without restoring the rise decorations
recorded by $\chi$.
By Lemma~\ref{lemForestStarValleys}\textup{(2)},
a retained decorated valley is a star valley
before deletion if and only if it is a
star valley after deletion.
Thus inversion pairs involving only retained steps, together
with their correction terms, contribute equally to
$\dinv(\pi)$ and $\dinv(\pi')$.
No retained decorated rise forms a correction pair with an
inserted root, since its area is at least two.

The inversion pairs of type \textup{(1)} in~\eqref{eqRefinedS13Inversions} between inserted roots
contribute
$\sum_{j<k,\,\delta_j=0}\mathbf{1}_{\{c_j\prec c_k\}}$.
At this stage, every inserted decorated valley lies outside
$\DD^*(\pi)$, so the term
$-|\DD(\pi)\setminus\DD^*(\pi)|$ contributes
$-\sum_j\delta_j$.
Their sum is $\dinv_0(W)$ by~\eqref{eqRefinedS45RootStatistics}.
An inversion pair involving one retained step and one inserted
root must be of type \textup{(2)}.
The retained step must have area $1$ and lie outside $\DD(\pi)$,
so it is a block head.
These pairs contribute $C(W;\delta,h;\beta)$ by~\eqref{eqRefinedS45RootStatistics}.

Now restore a rise decoration recorded by $\chi_j$, where $j<L$.
If $\delta_{j+1}=1$, the next inserted root now belongs to
$\DD^*(\pi)$; no other decorated valley enters or leaves
$\DD^*(\pi)$.
Consequently, the term
$-|\DD(\pi)\setminus\DD^*(\pi)|$ increases by $\delta_{j+1}$,
and the new inversion pairs of types \textup{(3)} and
\textup{(4)} contribute
\[
  \delta_{j+1}\sum_{k=j+2}^{L}
  \left(
    \delta_k\mathbf{1}_{\{c_{j+1}\prec c_k\}}
    +(1-\delta_k)\mathbf{1}_{\{c_k\prec c_{j+1}\}}
  \right).
\]
The restored decorated rise also forms a correction pair
with every later root that is a decorated valley, contributing
$-\sum_{k=j+1}^{L}\delta_k$.
The term $-\delta_{j+1}$ cancels the contribution
$+\delta_{j+1}$ above, leaving $\nu_j$.
For $j=L$, there is no later root, so the change is zero.
This proves the first identity.

Finally, each retained north step outside $\RR(\pi')$
contributes one additional unit to the area, except for the
$\sum_j\chi_j$ first children whose rise decorations are restored.
These first children have area zero in $\pi'$ and area one
before restoring their decorations, so each restoration removes
exactly one unit.
The inserted roots contribute no area.
This proves the second identity.
\end{proof}

\begin{example}
For the parking functions in Example~\ref{exDeletingMainDiagonal},
counting the four types of inversion pairs and the two correction
terms gives
\[
  \dinv(\pi)=10+3+1+2-1-3=12,
  \qquad
  \dinv(\pi')=5+3+1+1-1-2=7.
\]
The positive area entries of $\pi$ at steps outside $\RR(\pi)$
occur at $2,4,7$, each with value $1$.
Every positive area entry of $\pi'$ occurs at a step in
$\RR(\pi')$.
Hence $\area(\pi)=3$ and $\area(\pi')=0$.

Together with the values calculated from the deletion data in
Example~\ref{exMainDiagonalCorrections}, this gives
\[
  \dinv(\pi)-\dinv(\pi')=12-7=5=4+0+1,
  \qquad
  \area(\pi)-\area(\pi')=3-0=3=6-2-1,
\]
verifying both identities of
Theorem~\ref{thmMainDiagonalStatistics}.
\end{example}

\section{Insertion without rise decorations}\label{secInsertionsWithoutRise}
In this section, we consider insertions without rise decorations.
We fix the labels of the undecorated roots in their order
and the multiset of labels of the roots carrying a valley
decoration.
We then sum over all admissible ways of inserting the latter
roots and the heads among the fixed undecorated roots.
The sum includes the row weights of the heads, the label
weights of the inserted roots, and the changes in $\dinv$,
with the weights and mutual inversions of the fixed
undecorated roots removed.
Using Theorem~\ref{thmMainDiagonalStatistics}, we evaluate
this sum in Theorem~\ref{thmRootHeadInsertion}.

\subsection{Head labels and linear functionals}\label{subsecHeadLinearFunctional}
Fix admissible deletion data $(\pi';W,\delta,h,\chi)$ as in
Definition~\ref{defAdmissibleDeletionData}, and let
$\beta=(\beta_1,\ldots,\beta_H)$ be the head labels, in their order along the path.
We isolate their label weights and the inversion pairs between heads.

Recall the signed alphabet $1^-<1^+<\cdots<n^-<n^+$ from Definition~\ref{defRefinedS13SignedWeights} and the comparison rule $\prec$ from Definition~\ref{defRefinedS13SignedComparison}.
For a word $\beta$, set
\begin{equation}\label{eqRefinedS45HeadInversions}
  I(\beta)
  :=
  \sum_{b<b'}\mathbf{1}_{\{\beta_b\prec\beta_{b'}\}}
  =
  \#\{b<b':\beta_b<\beta_{b'}\}
  +\sum_{i=1}^{n}\binom{m_i}{2},
  \text{ where }
  m_i:=\#\{b:\beta_b=i^-\}.
\end{equation}

Let the letters $Z_c$, for $c\in\{1^-,1^+,\ldots,n^-,n^+\}$, generate the associative algebra over
$\mathbb{Q}(q,t)$ subject to the relations
\begin{equation}\label{eqRefinedS45ZRelations}
  Z_cZ_{c'}=qZ_{c'}Z_c\qquad\text{for }c<c'.
\end{equation}
Using~\eqref{eqRefinedS45ZRelations}, every word reduces to
the expression below, independently of the order of the exchanges.
The resulting linear reduction respects the defining relations
and fixes every weakly decreasing word, so these words form a basis:
\begin{equation}\label{eqRefinedS45ZOrdering}
  Z_{\beta_1}\cdots Z_{\beta_H}
  =
  q^{\#\{b<b':\beta_b<\beta_{b'}\}}
  Z_{\beta_{i_1}}\cdots Z_{\beta_{i_H}},
  \qquad
  \beta_{i_1}\geq\cdots\geq\beta_{i_H}.
\end{equation}

For $A\subseteq[n]$, define
\begin{equation}\label{eqRefinedS45Evaluation}
  e_A(\beta)
  :=
  \prod_{b=1}^{H}q^{\#\{a\in A:\beta_b\prec a^-\}},
  \qquad
  \operatorname{ev}_A
  \bigl(Z_{\beta_1}\cdots Z_{\beta_H}\bigr)
  :=
  q^{I(\beta)}e_A(\beta)\prod_{b=1}^{H}\wt(\beta_b),
\end{equation}
where $\wt(i^-)=-x_i$ and $\wt(i^+)=qx_i$ as in
Definition~\ref{defRefinedS13SignedWeights}, and $\operatorname{ev}_A(1)=1$.
Extend $\operatorname{ev}_A$ linearly.

For $c<c'$, replacing adjacent letters $c',c$ by $c,c'$
increases $I(\beta)$ by $1$ and leaves the label weights
and $e_A(\beta)$ unchanged.  Hence
\(
  \operatorname{ev}_A(\cdots Z_cZ_{c'}\cdots)
  =
  q\,\operatorname{ev}_A(\cdots Z_{c'}Z_c\cdots),
\)
in agreement with~\eqref{eqRefinedS45ZRelations}.
We set
\begin{equation}\label{eqRefinedS45QGenerators}
  Q_i:=Z_{i^-}+Z_{i^+}
\end{equation}
and have the following proposition.

\begin{proposition}\label{propHeadLinearFunctional}
\leavevmode
\begin{enumerate}
  \item For $i<j$, we have $Q_iQ_j=qQ_jQ_i$.
  \item For $A\subseteq[n]$ and $m\geq2$, we have
  $\operatorname{ev}_A(Q_i^m)=0$.
\end{enumerate}
\end{proposition}

\begin{proof}
\textup{(1)}
For $i<j$, every letter $i^\pm$ precedes every letter $j^\pm$.
Equation~\eqref{eqRefinedS45ZRelations} gives
\[
  Q_iQ_j
  =(Z_{i^-}+Z_{i^+})(Z_{j^-}+Z_{j^+})
  =q(Z_{j^-}+Z_{j^+})(Z_{i^-}+Z_{i^+})
  =qQ_jQ_i.
\]

\textup{(2)}
We use induction on $m$.
For a word of length $m$ in $i^-,i^+$, adding $i^-$ at the
beginning creates $m$ pairs counted by $I$, whereas adding
$i^+$ creates none, by~\eqref{eqRefinedS45HeadInversions}.
The combined factors contributed by $i^-$ and $i^+$ to the
label weight and $e_A$ are, respectively,
$-x_iq^{\#\{a\in A:i\leq a\}}$ and
$qx_iq^{\#\{a\in A:i<a\}}$.
Thus, for $m\geq0$,
\begin{align}
  \operatorname{ev}_A(Q_i^{m+1})
  &=
  \operatorname{ev}_A(Z_{i^-}Q_i^m)
  +\operatorname{ev}_A(Z_{i^+}Q_i^m)\notag\\
  &=
  \left(
    -x_iq^{m+\#\{a\in A:i\leq a\}}
    +qx_iq^{\#\{a\in A:i<a\}}
  \right)\operatorname{ev}_A(Q_i^m)\notag\\
  &=
  x_iq^{\#\{a\in A:i<a\}}
  \bigl(q-q^{m+\mathbf{1}_{\{i\in A\}}}\bigr)
  \operatorname{ev}_A(Q_i^m).\label{eqRefinedS45RepeatedIndexEvaluation}
\end{align}
If $i\in A$, taking $m=0$ gives
$\operatorname{ev}_A(Q_i)=0$, and then $m=1$ gives
$\operatorname{ev}_A(Q_i^2)=0$.
If $i\notin A$, taking $m=1$ makes the factor in parentheses
equal to $q-q=0$, so again $\operatorname{ev}_A(Q_i^2)=0$.
Equation~\eqref{eqRefinedS45RepeatedIndexEvaluation} now proves
$\operatorname{ev}_A(Q_i^m)=0$ for every $m\geq2$ by induction.
\end{proof}

\begin{definition}\label{defRefinedS45Quotient}
Let $\mathcal Q_n$ be the quotient of the subalgebra generated
by $Q_1,\ldots,Q_n$ by the two-sided ideal generated by
$Q_i^2$ for $1\leq i\leq n$.
We use the same symbols $Q_i$ for their images in $\mathcal Q_n$.
\end{definition}

The functionals $\operatorname{ev}_A$ are well defined on
$\mathcal Q_n$.  In a decreasing product of powers of the
$Q_i$, no inversion joins letters with different absolute indices.
The remaining factors in~\eqref{eqRefinedS45Evaluation}
also separate over these indices.  Thus, by~\eqref{eqRefinedS45HeadInversions}
and~\eqref{eqRefinedS45Evaluation},
\[
  \operatorname{ev}_A(Q_n^{m_n}\cdots Q_1^{m_1})
  =
  \prod_{i=1}^{n}\operatorname{ev}_A(Q_i^{m_i}),
\]
which vanishes whenever some $m_i\geq2$, by
Proposition~\ref{propHeadLinearFunctional}(2).
Together with Proposition~\ref{propHeadLinearFunctional}(1),
this shows that $\operatorname{ev}_A$ vanishes on the ideal
generated by the $Q_i^2$.

For $I=\{i_1<\cdots<i_m\}\subseteq[n]$, define
\begin{equation}\label{eqRefinedS45OrderedQ}
  Q_I:=Q_{i_m}\cdots Q_{i_1}\qquad\text{and}\qquad Q_\varnothing:=1.
\end{equation}
Before passing to $\mathcal Q_n$, the monomials
$Q_n^{m_n}\cdots Q_1^{m_1}$ are linearly independent:
their terms $Z_{n^-}^{m_n}\cdots Z_{1^-}^{m_1}$ are distinct
and have coefficient $1$.
By Proposition~\ref{propHeadLinearFunctional}\textup{(1)},
the ideal generated by the $Q_i^2$ is spanned by the
ordered monomials with some $m_i\geq2$.
Consequently, the $\{Q_I: I\subseteq [n]\}$ form a basis of $\mathcal Q_n$
and $\dim\mathcal Q_n=2^n$.
By~\eqref{eqRefinedS45OrderedQ} and~\eqref{eqRefinedS45Evaluation},
the evaluations of these basis elements are
\begin{equation}\label{lemEvaluationOrderedMonomials}
  \operatorname{ev}_A(Q_I)
  =
  \begin{cases}
    (q-1)^{|I|}
    q^{\#\{(i,a)\in I\times A:i<a\}}
    \displaystyle\prod_{i\in I}x_i,
      & I\cap A=\varnothing,\\[4pt]
    0,
      & I\cap A\neq\varnothing.
  \end{cases}
\end{equation}

In $\mathcal Q_n$, define
\begin{equation}\label{eqRefinedS45HeadElements}
  S:=\sum_{i=1}^{n}Q_i,\qquad
  T_c:=\sum_{\substack{1\leq i\leq n\\c\prec i^-}}Q_i,
  \qquad
  \sigma_c(Q_i):=q^{\mathbf{1}_{\{c\prec i^-\}}}Q_i.
\end{equation}
The element $S$ allows any signed head label, whereas $T_c$
allows only labels satisfying the comparison with $c$.
Here $c\prec i^-$ is equivalent to $c\prec i^+$.
The map $\sigma_c$ defines an automorphism of $\mathcal Q_n$,
since it preserves $Q_i^2=0$ and $Q_iQ_j=qQ_jQ_i$ for $i<j$
in Definition~\ref{defRefinedS45Quotient} and
Proposition~\ref{propHeadLinearFunctional}\textup{(1)}.

\begin{example}\label{exHeadLinearFunctional}
Recall the path in Example~\ref{exDeletingMainDiagonal}, whose head word is
$\beta=(1^-,3^+)$.
Take $A=\{2,8\}$.
Then $I(\beta)=1$ and $e_A(\beta)=q^3$, so
$\operatorname{ev}_A(Z_{1^-}Z_{3^+})
  =
  q(-x_1)(qx_3)q^3
  =
  -q^5x_1x_3$.
Summing the signs at these two absolute indices gives
$\operatorname{ev}_A(Q_1Q_3)
  =
  q\,\operatorname{ev}_A(Q_3Q_1)
  =
  q^4(q-1)^2x_1x_3$.
\end{example}

\begin{proposition}\label{propFiniteHeadFactor}
Recall the notation for deletion data from
Definition~\ref{defRefinedS45DeletionData}, and fix $W,\delta$,
an index $1\leq j\leq L$, and $A\subseteq[n]$.
For each $m\geq1$, sum over all signed head words
$\beta=(\beta_1,\ldots,\beta_m)$ satisfying $c_j\prec\beta_1$.
The first head is either undecorated, with $\chi_j=0$,
or a decorated rise, with $\chi_j=1$; all remaining heads
are undecorated.
Then
\begin{align*}
  \sum_{m\geq1}
  \sum_{\substack{\beta=(\beta_1,\ldots,\beta_m)\\
                  c_j\prec\beta_1}}
  \sum_{\chi_j \in \{0,1\}}
  (ut)^{m-\chi_j}r^{\chi_j}
  q^{I(\beta)+\chi_j\nu_j}
  e_A(\beta)\prod_{b=1}^{m}\wt(\beta_b)
  = \operatorname{ev}_A
  \left((ut+rq^{\nu_j})T_{c_j}(1-utS)^{-1}\right).
\end{align*}
If the empty head word is allowed by
Proposition~\ref{propMainDiagonalDeletion}\textup{(2)},
its contribution is $1$.
\end{proposition}

\begin{proof}
Fix the number of heads $m\geq1$.
Equations~\eqref{eqRefinedS45HeadElements} and~\eqref{eqRefinedS45Evaluation} give
\[
  \sum_{\substack{\beta=(\beta_1,\ldots,\beta_m)\\
                  c_j\prec\beta_1}}
  q^{I(\beta)}e_A(\beta)\prod_{b=1}^{m}\wt(\beta_b)
  =
  \operatorname{ev}_A(T_{c_j}S^{m-1}).
\]

All heads have area entry $1$ in the reconstructed path.
If $\chi_j=0$, the first head is undecorated, and the row
weights of the $m$ heads multiply to $(ut)^m$.
If $\chi_j=1$, the first head is a decorated rise, and
their row weights multiply to $r(ut)^{m-1}$.
In this case, Theorem~\ref{thmMainDiagonalStatistics}
gives the additional factor $q^{\nu_j}$.
Thus summing over $\chi_j\in\{0,1\}$ gives
\[
  \sum_{\chi_j\in\{0,1\}}
  (ut)^{m-\chi_j}r^{\chi_j}q^{\chi_j\nu_j}
  =(ut)^m+r(ut)^{m-1}q^{\nu_j}
  =(ut)^{m-1}(ut+rq^{\nu_j}).
\]

By Definition~\ref{defRefinedS45Quotient} and
Proposition~\ref{propHeadLinearFunctional}\textup{(1)},
every product of $n+1$ generators vanishes in $\mathcal Q_n$.
Since $T_{c_j}$ and $S$ are linear combinations of the
generators, $T_{c_j}S^{m-1}=0$ for $m>n$.
Summing over $m$ therefore gives
\begin{align*}
  \sum_{m=1}^{n}
  (ut)^{m-1}(ut+rq^{\nu_j})
  \operatorname{ev}_A(T_{c_j}S^{m-1})
  &=
  \operatorname{ev}_A
  \left(
    (ut+rq^{\nu_j})T_{c_j}
    \sum_{m=0}^{n-1}(utS)^m
  \right)\\
  &=
  \operatorname{ev}_A
  \left(
    (ut+rq^{\nu_j})T_{c_j}(1-utS)^{-1}
  \right),
\end{align*}
where the last equality follows by expanding the geometric
series and using $T_{c_j}S^n=0$.
\end{proof}

\subsection{Evaluation of the insertion sum}\label{subsecNoRiseInsertionSum}
Throughout this subsection, set $r=0$.
Let $B=(b_1,\ldots,b_{|B|})$ list the labels of the undecorated
deleted roots in order, and let $E$ be the multiset of labels
of the deleted roots that are decorated valleys.
For $d\in\operatorname{supp}(E)$, let $\mu_d$ be its multiplicity,
where $\operatorname{supp}(E)$ is the set of distinct labels in $E$.

Insert a formal letter $\xi_i$ between $i^-$ and $i^+$, having
the order
$1^-<\xi_1<1^+<2^-<\xi_2<2^+<\cdots<n^-<\xi_n<n^+$.
Recall the comparison rule $\prec$ from Definition~\ref{defRefinedS13SignedComparison}, which we use on this alphabet.
Thus, for a signed label $c$, the conditions
$c\prec i^-$, $c<\xi_i$, and $c\prec i^+$ are equivalent.

For $I\subseteq[n]$ and a letter $z$ in this alphabet, define
\begin{equation}\label{eqRefinedS45InsertionCounts}
  \lambda_z(B):=\#\{1\leq i\leq|B|:b_i\prec z\},
  \qquad
  n_{>z}:=\#\{i\in I:\xi_i>z\},
  \qquad
  \rho_z:=\lambda_z(B)+n_{>z}.
\end{equation}
Set $\epsilon_d=1$ if $d$ is negative and $\epsilon_d=0$
if $d$ is positive.
For nonnegative integers $\mu,\rho$ and $\epsilon\in\{0,1\}$,
define
\begin{equation}\label{eqRefinedS45PhiFactors}
  \varphi_{\mu,\epsilon}(\rho):=
  \begin{cases}
    h_\mu(1,q^{-1},\ldots,q^{-(\rho-1)}),&\epsilon=1,\\
    e_\mu(1,q^{-1},\ldots,q^{-(\rho-1)}),&\epsilon=0.
  \end{cases}
\end{equation}
The alphabet in~\eqref{eqRefinedS45PhiFactors} is empty when $\rho=0$.
Finally, set
\begin{equation}\label{eqRefinedS45InsertionProduct}
  W_I(B,E):=
  \prod_{i\in I}[\rho_{\xi_i}]_q
  \prod_{d\in\operatorname{supp}(E)}
  q^{\mu_dn_{>d}}\,
  \varphi_{\mu_d,\epsilon_d}(\rho_d).
\end{equation}

Define
\begin{equation}\label{eqRefinedS45NoRiseInsertion}
  K_{B,E}(ut,0)
  :=
  \sum_{W,\delta,h,\beta}
  (ut)^H
  q^{\dinv_0(W)-I(B)+C(W;\delta,h;\beta)}
  \left(\prod_{d\in\operatorname{supp}(E)}
    \wt(d)^{\mu_d}\right)
  Q_{i_1}\cdots Q_{i_H}.
\end{equation}
For $L\geq1$, the sum runs over $W,\delta,h,\beta$
satisfying the conditions of Proposition~\ref{propMainDiagonalDeletion}
with $\chi=0$, where comparisons are interpreted in the
extended alphabet containing the letters $\xi_i$.
We require $B$ to be the sequence of labels $c_j$ with
$\delta_j=0$, taken in increasing order of $j$, and $E$
to be the multiset of labels $c_j$ with $\delta_j=1$.
Here $\beta=(\xi_{i_1},\ldots,\xi_{i_H})$ ranges over finite
head words, with $H=\sum_{j=1}^{L}h_j$.
The factors $Q_{i_1},\ldots,Q_{i_H}$ follow this order.

The exponent $\dinv_0(W)-I(B)+C(W;\delta,h;\beta)$ follows
from Theorem~\ref{thmMainDiagonalStatistics}, after removing
the comparisons between letters of $B$.
We also allow the unique empty insertion with $L=H=0$ and
$W,\delta,h,\beta$ all empty, and give it weight $1$.
Thus $K_{\varnothing,\varnothing}(ut,0)=1$, whereas
$K_{\varnothing,E}(ut,0)=0$ for $E\neq\varnothing$.

\begin{theorem}\label{thmRootHeadInsertion}
For fixed $B,E$, we have
\[
  K_{B,E}(ut,0)
  =
  \sum_{I\subseteq[n]}
  (ut)^{|I|}W_I(B,E)
  \left(
    \prod_{d\in\operatorname{supp}(E)}
    \wt(d)^{\mu_d}q^{\mu_d(\lambda_d(B)-1)}
  \right)Q_I.
\]
\end{theorem}

\begin{proof}
In~\eqref{eqRefinedS45NoRiseInsertion}, we have $\chi=0$.
If $B=\varnothing$, every nonempty choice of heads and labels
in $E$ has a largest letter $z$ with $\rho_z=0$, giving a zero
factor on the right-hand side.
Only $E=I=\varnothing$ contributes $1$, in agreement
with~\eqref{eqRefinedS45NoRiseInsertion}.
Hence assume $B\neq\varnothing$.

A product containing a repeated $Q_i$ vanishes by
Definition~\ref{defRefinedS45Quotient} and
Proposition~\ref{propHeadLinearFunctional}\textup{(1)}.
Fix $I\subseteq[n]$ and calculate the coefficient of $Q_I$.
Represent $(W,\delta,h,\beta)$ by placing the assigned heads,
in order, after each root of $W$.
By Proposition~\ref{propMainDiagonalDeletion} with $\chi=0$,
the resulting word begins with a letter of $B$ and satisfies
$c\prec z$ whenever a root $c$ is immediately followed by
a head $z$ or a letter $z\in E$.
The occurrences from $B$ retain their prescribed order.

Starting with $B$, insert the heads $\xi_i$, for $i\in I$,
and the labels in $E$ in decreasing order, inserting all copies
of each repeated label at the same step.
For a letter $z$, the possible positions are immediately after
a letter $b\in B$ with $b\prec z$ or an already inserted head
greater than $z$.
There are $\lambda_z(B)+n_{>z}=\rho_z$ such positions
by~\eqref{eqRefinedS45InsertionCounts}.
The roots already inserted from $E$ are larger than $z$
and cannot immediately precede it.
The following letter is either in $B$, with no comparison
required, or is an already inserted letter greater than $z$,
so insertion preserves the required comparisons.

Delete the inserted letters in increasing order, taking
all copies of each label together.
Every new required comparison is with a letter larger than
the deleted label, so the comparisons remain valid and
the insertion positions and multiplicities are recovered uniquely.

The comparisons between letters of $B$ contribute $I(B)$,
which is cancelled in~\eqref{eqRefinedS45NoRiseInsertion}.
Suppose a head $\xi_i$ is inserted in the $k$th possible
position from the left.
Of the letters after which $\xi_i$ may be inserted,
exactly $\rho_{\xi_i}-k$ lie to its right.
Each such root contributes $q$ through $C$, and each such
head contributes $q$ when the $Q$ factors are put in decreasing
index order, by~\eqref{eqRefinedS45RootStatistics} and
Proposition~\ref{propHeadLinearFunctional}\textup{(1)}.
Thus summing over the positions gives
\[
  \sum_{k=1}^{\rho_{\xi_i}}q^{\rho_{\xi_i}-k}
  =[\rho_{\xi_i}]_q.
\]

Next insert a root with label $d\in E$ in the $k$th possible
position from the left.
Exactly $k$ preceding letters are roots $b\in B$ with
$b\prec d$ or heads larger than $d$.
Each such root contributes $1$ to $\dinv_0(W)$, and each such
head contributes $1$ to $C$.
The decoration of $d$ contributes $-1$ through
$-\sum_j\delta_j$ in~\eqref{eqRefinedS45RootStatistics}.
Hence the insertion has weight $\wt(d)q^{k-1}$.

Negative copies may occupy the same position, since $d\prec d$,
whereas positive copies must occupy distinct positions.
There are no additional inversion pairs between copies of $d$,
since a decorated root cannot be the first root in a pair
counted by $\dinv_0(W)$.
Consequently, summing over all $\mu_d$ copies gives,
by~\eqref{eqRefinedS45PhiFactors},
\[
  \wt(d)^{\mu_d}
  \begin{cases}
    h_{\mu_d}(1,q,\ldots,q^{\rho_d-1}),&\epsilon_d=1,\\
    e_{\mu_d}(1,q,\ldots,q^{\rho_d-1}),&\epsilon_d=0
  \end{cases}
  =
  \wt(d)^{\mu_d}q^{\mu_d(\rho_d-1)}
  \varphi_{\mu_d,\epsilon_d}(\rho_d).
\]
The alphabets are empty when $\rho_d=0$.

At each step, the list of position weights depends only on
$B$ and the labels already inserted, not on their chosen
positions.
The successive sums therefore multiply.
Including the factor $(ut)^{|I|}$, the coefficient of $Q_I$ is
\[
  (ut)^{|I|}
  \prod_{i\in I}[\rho_{\xi_i}]_q
  \prod_{d\in\operatorname{supp}(E)}
  \left(
    \wt(d)^{\mu_d}q^{\mu_d(\rho_d-1)}
    \varphi_{\mu_d,\epsilon_d}(\rho_d)
  \right).
\]
Since $\rho_d=\lambda_d(B)+n_{>d}$,
\eqref{eqRefinedS45InsertionProduct} gives the asserted formula.
\end{proof}

\begin{example}
Take $B=(1^-,2^+,6^+)$, $E=\{8^-\}$, and $I=\{1,3\}$,
as in Example~\ref{exDeletingMainDiagonal}.
We calculate the coefficient of $Q_I$ in $K_{B,E}(ut,0)$
directly and compare it with Theorem~\ref{thmRootHeadInsertion}.

In~\eqref{eqRefinedS45NoRiseInsertion}, there are six admissible insertions for
each head order.
Using $Q_1Q_3=qQ_I$ and $Q_3Q_1=Q_I$, their contributions give
\[
  [Q_I]K_{B,E}(ut,0)
  =
  -(ut)^2x_8
  \left(
    \underbrace{q+2q^2+2q^3+q^4}_{(\xi_1,\xi_3)}
    +
    \underbrace{1+2q+2q^2+q^3}_{(\xi_3,\xi_1)}
  \right).
\]

Meanwhile, $\rho_{\xi_1}=\rho_{\xi_3}=2$,
$\rho_{8^-}=\lambda_{8^-}(B)=3$, and $\wt(8^-)=-x_8$.
Theorem~\ref{thmRootHeadInsertion} gives
\[
  [Q_I]K_{B,E}(ut,0)
  =
  -(ut)^2x_8q^2(1+q)^2(1+q^{-1}+q^{-2})
  =
  -(ut)^2x_8(1+3q+4q^2+3q^3+q^4),
\]
in agreement with the direct calculation from~\eqref{eqRefinedS45NoRiseInsertion}.
\end{example}

\section{Insertions with rise decorations}\label{secInsertionsWithRise}
In this section, we extend the insertion sum in \eqref{eqRefinedS45NoRiseInsertion}
to allow decorated rises.
Using the $q$-series identities in Proposition~\ref{propGIdentities}, we prove in
Theorem~\ref{thmRiseDifference} that the dependence on $r$
is determined by the case without rise decorations.
We then sum over the labels of roots carrying a valley
decoration, first for $r=0$ and then for general $r$.
This calculation gives $M_{X_n}(ut,r,s)$
in Theorem~\ref{thmRiseTransferMatrix}, which we use in
Section~\ref{secCombinatorialRecurrence} to derive the combinatorial recurrence.

\subsection{\texorpdfstring{$q$}{q}-series identities}\label{subsecQSeries}
For $j\geq0$, set
\begin{equation}\label{eqRefinedS6gDefinition}g_j(ut,r):=\prod_{h=0}^{j-1}(r+utq^h),\qquad g_0(ut,r):=1.\end{equation}

If $U$ is a linear combination of $Q_1,\ldots,Q_n$ with scalar coefficients, define
\begin{equation}\label{eqdefGr}
  G_r(U):=
  \sum_{j\geq0}
  \frac{(-1)^j g_j(ut,r)}{(q;q)_j}U^j.
\end{equation}
  
Since $U^{n+1}=0$, the sum in~\eqref{eqdefGr} is finite.
Its constant term in $U$ is $1$, so its inverse is also given
by a finite geometric series.

\begin{proposition}\label{propGIdentities}
The following identities hold.
\begin{enumerate}[(1)]
  \item $G_{qr}(U)=(1+rU)G_r(U)$.

  \item $G_r(qU)=(1-utU)^{-1}(1+rU)G_r(U)$.

  \item For $c\in\{1^-,1^+,\ldots,n^-,n^+\}$, we have
  \[
    G_R(\sigma_c(S))G_R(S)^{-1}
    =
    1+(ut+R)T_c(1-utS)^{-1},
  \]
  where $G_R$ is obtained by replacing $r$ with a formal
  parameter $R$ in \eqref{eqdefGr}.
\end{enumerate}
\end{proposition}

\begin{proof}
\textup{(1)}
For $j\geq1$, \eqref{eqRefinedS6gDefinition} gives
$g_j(ut,qr)-g_j(ut,r)=r(q^j-1)g_{j-1}(ut,r)$.
Substituting into~\eqref{eqdefGr}, we obtain
\[
  G_{qr}(U)-G_r(U)
  =
  rU\sum_{j\geq0}
  \frac{(-1)^jg_j(ut,r)}{(q;q)_j}U^j
  =
  rUG_r(U).
\]

\textup{(2)}
Using $g_{j+1}(ut,r)=(r+utq^j)g_j(ut,r)$ from~\eqref{eqRefinedS6gDefinition} in~\eqref{eqdefGr} gives
\[
  G_r(qU)-G_r(U)
  =
  U\bigl(rG_r(U)+utG_r(qU)\bigr).
\]
Rearranging and multiplying by $(1-utU)^{-1}$ proves \textup{(2)}.

\textup{(3)}
By~\eqref{eqRefinedS45HeadElements} and
Definition~\ref{defRefinedS13SignedComparison},
every index occurring in $S-T_c$ is smaller than every
index occurring in $T_c$.
Thus Proposition~\ref{propHeadLinearFunctional}\textup{(1)}
gives
\begin{align*}
  (S-T_c)T_c
  =\sum_{\substack{1\leq i,j\leq n\\
       c\not\prec i^-,\,c\prec j^-}}Q_iQ_j
  =q\sum_{\substack{1\leq i,j\leq n\\
       c\not\prec i^-,\,c\prec j^-}}Q_jQ_i
   =qT_c(S-T_c).
\end{align*}
Also, \eqref{eqRefinedS45HeadElements} gives
$\sigma_c(S)=(S-T_c)+qT_c$, so
\begin{align*}
  \sigma_c(S)T_c
  =(S-T_c)T_c+qT_c^2
  =qT_c(S-T_c)+qT_c^2
   =qT_cS.
\end{align*}

We now prove by induction that, for $j\geq1$,
$(\sigma_c(S))^j-S^j=(q^j-1)T_cS^{j-1}$.
The case $j=1$ is
$\sigma_c(S)-S=(q-1)T_c$.
For the induction step, the preceding identities give
\begin{align*}
  (\sigma_c(S))^{j+1}-S^{j+1}
  =\sigma_c(S)\bigl((\sigma_c(S))^j-S^j\bigr)
    +\bigl(\sigma_c(S)-S\bigr)S^j
  =(q^j-1)\sigma_c(S)T_cS^{j-1}
    +(q-1)T_cS^j
  =(q^{j+1}-1)T_cS^j.
\end{align*}

Substituting into~\eqref{eqdefGr}, using
$(q;q)_j=(1-q^j)(q;q)_{j-1}$, and then using
$g_{j+1}(ut,R)=(R+utq^j)g_j(ut,R)$
from~\eqref{eqRefinedS6gDefinition}, we obtain
\begin{align*}
  G_R(\sigma_c(S))-G_R(S)
  =\sum_{j\geq1}
    \frac{(-1)^jg_j(ut,R)(q^j-1)}{(q;q)_j}
    T_cS^{j-1}
  =T_c\sum_{j\geq0}
    \frac{(-1)^jg_{j+1}(ut,R)}{(q;q)_j}S^j
  =T_c\bigl(RG_R(S)+utG_R(qS)\bigr).
\end{align*}
By \textup{(2)} with $r$ replaced by $R$,
\begin{align*}
  T_c\bigl(RG_R(S)+utG_R(qS)\bigr)
  =T_c\bigl(R+ut(1-utS)^{-1}(1+RS)\bigr)G_R(S)
  =(ut+R)T_c(1-utS)^{-1}G_R(S),
\end{align*}
where the last equality uses
$R(1-utS)+ut(1+RS)=R+ut$.
Multiplying on the right by $G_R(S)^{-1}$ proves \textup{(3)}.
\end{proof}

\subsection{The contribution of rise decorations}\label{subsecRiseContribution}
Fix $B,E$.  Recall the summation indices and notation
in~\eqref{eqRefinedS45NoRiseInsertion}.
When $h_j>0$, we now allow the first head assigned to root $j$
to be a decorated rise.
By Theorem~\ref{thmMainDiagonalStatistics}, these heads
contribute $(ut)^{h_j}$ if the first head is undecorated,
and $r(ut)^{h_j-1}q^{\nu_j}$ if it is a decorated rise.
We then define
\begin{equation}\label{eqRefinedS6KDefinition}
  K_{B,E}(ut,r)
  :=
  \sum_{W,\delta,h,\beta}
  q^{\dinv_0(W)-I(B)+C(W;\delta,h;\beta)}
  \left(
    \prod_{d\in\operatorname{supp}(E)}\wt(d)^{\mu_d}
  \right)
  \left(
    \prod_{j:h_j>0}
    (ut)^{h_j-1}(ut+rq^{\nu_j})
  \right)
  Q_{i_1}\cdots Q_{i_H}.
\end{equation}
Setting $r=0$ recovers~\eqref{eqRefinedS45NoRiseInsertion}.

For $H\in\mathcal Q_n$, let $[H]_j$ denote the sum of its
$Q_I$-terms with $|I|=j$.
Define the automorphism $\tau$ by
\begin{equation}\label{eqRefinedS6TauDefinition}
\tau(Q_I)=q^{|I|}Q_I,
\end{equation}
fixing all scalar parameters.

\begin{theorem}\label{thmRiseDifference}
For fixed root data $B,E$ and $0\leq j\leq n$, we have
\[
  q^{\binom{j}{2}}
  [K_{B,E}(ut,r)G_r(S)]_j
  =
  g_j(ut,r)
  \left.[K_{B,E}(ut,r)G_r(S)]_j\right|_{ut=1,\ r=0}.
\]
Equivalently,
\[
  [K_{B,E}(ut,r)G_r(S)]_j
  =
  \prod_{h=0}^{j-1}(ut+rq^{-h})
  \left.[K_{B,E}(ut,r)G_r(S)]_j\right|_{ut=1,\ r=0}.
\]
\end{theorem}

Write $B=cA$, where $c$ is the first letter of $B$.
Recall $\lambda_d(A)$ from~\eqref{eqRefinedS45InsertionCounts}.
For $d\in\operatorname{supp}(E)$, let $E-d$ denote the multiset
obtained by removing one occurrence of $d$, and set
\begin{equation}\label{eqRefinedS6RemainingCounts}
v_d:=\#\{e\in E-d:d\prec e\},
\end{equation}
counting multiplicities.

Using the commuting automorphisms from
\eqref{eqRefinedS45HeadElements}, define
\begin{equation}\label{eqRefinedS6SigmaWords}
  \sigma_{A+E}
  :=
  \left(\prod_{a\in A}\sigma_a\right)
  \left(\prod_{d\in\operatorname{supp}(E)}
    \sigma_d^{\mu_d}\right),
  \qquad
  \sigma_A:=\prod_{a\in A}\sigma_a.
\end{equation}
In $\sum_{A=c'A'}$, $c'$ is the first letter of $A$;
the sum is zero if $A=\varnothing$.

Define $\Phi_c(A,E;r)$ by
\begin{equation}\label{eqPhiNormalization}
  K_{cA,E}(ut,r)G_r(S)
  =
  \left(
    \prod_{d\in\operatorname{supp}(E)}
    \wt(d)^{\mu_d}q^{\mu_d(\lambda_d(cA)-1)}
  \right)
  G_{r/q^{|E|}}(\sigma_c\sigma_{A+E}S)
  \Phi_c(A,E;r).
\end{equation}
The scalar product is a common factor of $K_{cA,E}(ut,r)$.
After canceling it, the invertibility of the $G$-factor
determines $\Phi_c(A,E;r)$ uniquely.

\begin{proposition}\label{propPhiEmpty}
For every word $A$ and signed label $c$, we have
$K_{cA,\varnothing}(ut,r)G_r(S)=G_r(\sigma_c\sigma_A S)$,
and consequently
\begin{equation}\label{eqRefinedS6PhiEmpty}
 \Phi_c(A,\varnothing;r)=1.
\end{equation}
\end{proposition}
\begin{proof}
When $E=\varnothing$, we have $W=(c_1,\ldots,c_L)=cA$
and $\delta_j=0$ for every $j$.
Thus~\eqref{eqRefinedS45RootStatistics} and
\eqref{eqRefinedS45RiseCorrection} give
$\dinv_0(W)=I(cA)$ and $\nu_j=0$.

For the heads assigned to $c_j$, the empty block contributes $1$.
A nonempty block of $h$ heads has weight
$(ut)^h+r(ut)^{h-1}=(ut)^{h-1}(ut+r)$,
and only its first label is restricted, by $c_j\prec i_1^-$.
By~\eqref{eqRefinedS45RootStatistics}, its contribution to $C$
counts the comparisons of these heads with $c_{j+1},\ldots,c_L$.
Hence~\eqref{eqRefinedS45HeadElements} gives
\begin{align*}
  \sum_{\substack{1\leq i_1,\ldots,i_h\leq n\\
                   c_j\prec i_1^-}}
    q^{\sum_{b=1}^{h}\sum_{k=j+1}^{L}
       \mathbf{1}_{\{c_k\prec i_b^-\}}}
    Q_{i_1}\cdots Q_{i_h}
  =\sigma_{c_{j+1}}\cdots\sigma_{c_L}
    \left(
      \left(\sum_{\substack{1\leq i\leq n\\c_j\prec i^-}}Q_i\right)
      \left(\sum_{i=1}^{n}Q_i\right)^{h-1}
    \right)
  =\sigma_{c_{j+1}}\cdots\sigma_{c_L}(T_{c_j}S^{h-1}).
\end{align*}
An empty composition of automorphisms acts as the identity.

Since every root is undecorated,
Proposition~\ref{propMainDiagonalDeletion}\textup{(2)}
imposes no condition between successive blocks.
We may therefore sum the heads assigned to each root
independently, keeping their $Q$ factors in root order.
Expanding~\eqref{eqRefinedS6KDefinition} gives
\begin{align*}
  K_{cA,\varnothing}(ut,r)
  &=\prod_{j=1}^{L}
    \left(
      1+\sum_{h\geq1}(ut)^{h-1}(ut+r)
      \sigma_{c_{j+1}}\cdots\sigma_{c_L}
      (T_{c_j}S^{h-1})
    \right)\\
  &=\prod_{j=1}^{L}
    \sigma_{c_{j+1}}\cdots\sigma_{c_L}
    \left(
      1+(ut+r)T_{c_j}
      \sum_{h\geq1}(ut)^{h-1}S^{h-1}
    \right)\\
  &=\prod_{j=1}^{L}
    \sigma_{c_{j+1}}\cdots\sigma_{c_L}
    \left(1+(ut+r)T_{c_j}(1-utS)^{-1}\right).
\end{align*}
All products are ordered by increasing $j$.
The geometric sum is finite because $S^{n+1}=0$.

Apply $\sigma_{c_{j+1}}\cdots\sigma_{c_L}$ to
Proposition~\ref{propGIdentities}\textup{(3)} with $R=r$.
These automorphisms commute, fix the scalar coefficients
of $G_r$, and preserve inverses. Thus
\begin{align*}
  \sigma_{c_{j+1}}\cdots\sigma_{c_L}
    \left(1+(ut+r)T_{c_j}(1-utS)^{-1}\right)
  =\sigma_{c_{j+1}}\cdots\sigma_{c_L}
    \left(G_r(\sigma_{c_j}S)G_r(S)^{-1}\right)
  =G_r(\sigma_{c_j}\cdots\sigma_{c_L}S)
    G_r(\sigma_{c_{j+1}}\cdots\sigma_{c_L}S)^{-1}.
\end{align*}
Substituting into the preceding product and successively cancelling adjacent inverse factors gives
\begin{equation}\label{eqRiseEmptyTelescoping}
  K_{cA,\varnothing}(ut,r)G_r(S)
  =\left(
      \prod_{j=1}^{L}
      G_r(\sigma_{c_j}\cdots\sigma_{c_L}S)
      G_r(\sigma_{c_{j+1}}\cdots\sigma_{c_L}S)^{-1}
    \right)G_r(S)
    =G_r(\sigma_{c_1}\cdots\sigma_{c_L}S)
  =G_r(\sigma_c\sigma_A S).
\end{equation}

The scalar product in~\eqref{eqPhiNormalization} is $1$
when $E=\varnothing$.
Substituting~\eqref{eqRiseEmptyTelescoping} there and cancelling
the invertible left factor $G_r(\sigma_c\sigma_A S)$ gives
$\Phi_c(A,\varnothing;r)=1$, proving~\eqref{eqRefinedS6PhiEmpty}.
\end{proof}

\begin{example}
Take $B=(1^-,2^+,6^+)$ and $E=\{8^-\}$, as in
Example~\ref{exDeletingMainDiagonal}.
We now check the coefficient of $Q_1$ in
Theorem~\ref{thmRiseDifference}.

With no heads, the three possible positions of $8^-$ give
$[K_{B,E}(ut,r)]_0=-x_8(1+q+q^2)$.
A single head $\xi_1$ must immediately follow $1^-$.
Placing $8^-$ before $2^+$, between $2^+$ and $6^+$,
or after $6^+$ gives
$\dinv_0(W)-I(B)=0,1,2$ and $\nu_1=2,-1,-1$,
respectively; in each case, $C=0$.
Thus~\eqref{eqRefinedS6KDefinition} gives
\begin{align*}
  [Q_1]K_{B,E}(ut,r)
  =-x_8\bigl((ut+rq^2)+q(ut+rq^{-1})
      +q^2(ut+rq^{-1})\bigr)
  =-x_8(1+q+q^2)(ut+r).
\end{align*}
Since $[G_r(S)]_0=1$ and
$[Q_1]G_r(S)=(ut+r)/(q-1)$, the left-hand side gives
\[
  [Q_1]\bigl(K_{B,E}(ut,r)G_r(S)\bigr)
  =
  -\frac{qx_8(1+q+q^2)}{q-1}(ut+r).
\]
On the right-hand side, $g_1(ut,r)=ut+r$, and specialization
of the same calculation gives
\[
  g_1(ut,r)
  \left.
    [Q_1]\bigl(K_{B,E}(ut,r)G_r(S)\bigr)
  \right|_{ut=1,\,r=0}
  =
  (ut+r)\left(-\frac{qx_8(1+q+q^2)}{q-1}\right).
\]
The two coefficients agree.
\end{example}

We use the following two lemmas to prove Theorem~\ref{thmRiseDifference}.
Their technical proofs are given in Appendix~\ref{appRiseCalculations}.

\begin{lemma}\label{lemHomogeneousFactor}
Let $H\in\mathcal Q_n[ut,r]$, with $[H]_j$ homogeneous
of degree $j$ in $ut,r$.
Then
\[
  (qr+ut\tau)H(qr)=(qr+ut)\tau H(r) \quad \Longleftrightarrow  \quad
[H]_j=g_j(ut,r)C_j \text{ for every }j,
\]
where $C_j$ is independent of $ut,r$.
Here $H(qr)$ denotes substitution of $qr$ for $r$.
\end{lemma}

\begin{lemma}\label{lemPhiOperatorCalculation}
Let $E\neq\varnothing$.
Suppose that, for every $j\geq0$,
$[G_{r/q^{|E|}}(\sigma_{A+E}S)\Phi_{c'}(A',E;r)]_j$
for $A=c'A'$, and
$[G_{rq^{1-|E|}}(\sigma_{A+E}S)\Phi_d(A,E-d;r)]_j$
for $d\in\operatorname{supp}(E)$,
are divisible by $g_j(ut,r)$, with quotients independent of $ut,r$.
Then
\[
  (qr+ut\tau)
  \left(
    G_{rq^{1-|E|}}(\sigma_c\sigma_{A+E}S)
    \Phi_c(A,E;qr)
  \right) =
  (qr+ut)\tau
  \left(
    G_{r/q^{|E|}}(\sigma_c\sigma_{A+E}S)
    \Phi_c(A,E;r)
  \right).
\]
\end{lemma}

\begin{proof}[Proof of Theorem~\ref{thmRiseDifference}]
We first prove that, for every $c$ and $j\geq0$,
\begin{equation}\label{eqRefinedS6RiseClaim}
  \frac{
    [G_{r/q^{|E|}}(\sigma_c\sigma_{A+E}S)
      \Phi_c(A,E;r)]_j
  }{g_j(ut,r)}
  \quad\text{is independent of }ut,r.
\end{equation}
By~\eqref{eqRefinedS6KDefinition}, \eqref{eqdefGr},
and~\eqref{eqPhiNormalization}, every monomial in the
numerator has total degree $j$ in $ut,r$.
Lemma~\ref{lemHomogeneousFactor} therefore reduces
\eqref{eqRefinedS6RiseClaim} to proving
\begin{align*}
  (qr+ut\tau)
    \left(
      G_{rq^{1-|E|}}(\sigma_c\sigma_{A+E}S)
      \Phi_c(A,E;qr)
    \right)
  =(qr+ut)\tau
    \left(
      G_{r/q^{|E|}}(\sigma_c\sigma_{A+E}S)
      \Phi_c(A,E;r)
    \right).
\end{align*}

We use induction on $|A|+|E|$.
If $E=\varnothing$, then
\eqref{eqRefinedS6PhiEmpty} and~\eqref{eqdefGr} give
\[
  \frac{
    [G_r(\sigma_c\sigma_A S)\Phi_c(A,\varnothing;r)]_j
  }{g_j(ut,r)}
  =
  \frac{(-1)^j}{(q;q)_j}(\sigma_c\sigma_A S)^j,
\]
which proves~\eqref{eqRefinedS6RiseClaim} directly.

For $E\neq\varnothing$, $A=c'A'$, and
$d\in\operatorname{supp}(E)$, we have
\begin{gather*}
  |A'|+|E|=|A|+|E-d|=|A|+|E|-1,\quad
  \sigma_{c'}\sigma_{A'+E}
  =\sigma_d\sigma_{A+(E-d)}
  =\sigma_{A+E},
  \quad
  \frac{r}{q^{|E-d|}}=rq^{1-|E|}.
\end{gather*}
The induction hypothesis gives the two conditions in
Lemma~\ref{lemPhiOperatorCalculation}.
When $A=\varnothing$, only the pairs $(A,E-d)$ are needed.
Lemma~\ref{lemPhiOperatorCalculation} proves the displayed
identity, and hence~\eqref{eqRefinedS6RiseClaim}.

Now let $B=cA\neq\varnothing$.
Using~\eqref{eqPhiNormalization},
\eqref{eqRefinedS6RiseClaim}, and
$g_j(1,0)=q^{\binom{j}{2}}$, we obtain
\begin{align*}
  \frac{[K_{cA,E}(ut,r)G_r(S)]_j}{g_j(ut,r)}
  &=\left(
      \prod_{d\in\operatorname{supp}(E)}
      \wt(d)^{\mu_d}q^{\mu_d(\lambda_d(cA)-1)}
    \right)
    \frac{
      [G_{r/q^{|E|}}(\sigma_c\sigma_{A+E}S)
        \Phi_c(A,E;r)]_j
    }{g_j(ut,r)}\\
  &=\left(
      \prod_{d\in\operatorname{supp}(E)}
      \wt(d)^{\mu_d}q^{\mu_d(\lambda_d(cA)-1)}
    \right)
    \left.
    \frac{
      [G_{r/q^{|E|}}(\sigma_c\sigma_{A+E}S)
        \Phi_c(A,E;r)]_j
    }{g_j(ut,r)}
    \right|_{ut=1,r=0}\\
  &=q^{-\binom{j}{2}}
    \left.[K_{cA,E}(ut,r)G_r(S)]_j\right|_{ut=1,r=0}.
\end{align*}
Here the product over $d$ is independent of $ut,r$.

For $B=\varnothing$, the same equality follows from
\[
  K_{\varnothing,E}(ut,r)
  =
  \begin{cases}
    1,&E=\varnothing,\\
    0,&E\neq\varnothing,
  \end{cases}
  \qquad
  \frac{[G_r(S)]_j}{g_j(ut,r)}
  =\frac{(-1)^j}{(q;q)_j}S^j.
\]
Thus, for all $B,E$,
\[
  q^{\binom{j}{2}}[K_{B,E}(ut,r)G_r(S)]_j
  =
  g_j(ut,r)
  \left.[K_{B,E}(ut,r)G_r(S)]_j\right|_{ut=1,r=0}.
\]
The second form follows from
$q^{-\binom{j}{2}}g_j(ut,r)
=\prod_{h=0}^{j-1}(ut+rq^{-h})$.
\end{proof}

\subsection{The summation without rise decorations}\label{subsecNoRiseTransferMatrix}
Fix $X_n=x_1+\cdots+x_n$, and write
$x_I:=\prod_{i\in I}x_i$ for $I\subseteq[n]$.
For a word $\alpha$ giving the signed labels of the
undecorated deleted roots in order, define
\begin{equation}\label{eqRefinedS6HDefinition}
  \mathcal H_\alpha(ut;r,s)
  :=
  \sum_E s^{|E|}K_{\alpha,E}(ut,r),
\end{equation}
where $E$ ranges over multisets of signed labels for the
deleted roots that are decorated valleys.
The sum is understood as a formal power series in $s$,
with each coefficient given by a finite sum.

Recall $\operatorname{ev}_A$ and $e_B$ from
\eqref{eqRefinedS45Evaluation}.
With the notation of~\eqref{eqRefinedS45InsertionCounts},
recall that
$e_B(\alpha)=\prod_{i\in B}q^{\lambda_{\xi_i}(\alpha)}$,
where $\lambda_{\xi_i}(\alpha)$ counts the occurrences
of labels $c$ in $\alpha$ satisfying $c\prec i^-$.

\begin{proposition}\label{propNoRiseTransferMatrix}
For every such word $\alpha$ and every $A\subseteq[n]$, we have
\[
  \operatorname{ev}_A\mathcal H_\alpha(ut;0,s)
  =
  \sum_{B\subseteq[n]}e_B(\alpha)V_{X_n}(ut,s)_{B,A},
\]
where
\begin{equation}\label{eqRefinedS6VDefinition}
  V_{X_n}(ut,s)_{B,A}
  =
  \prod_{i=1}^{n}\frac{v_i(B,A)}{1+sx_i},
  \qquad
  v_i(B,A)=
  \begin{cases}
    1-utx_iq^{|A\cap[i,n]|+|B\cap[1,i)|},
      &i\notin A,\ i\notin B,\\
    utx_iq^{|A\cap[i,n]|+|B\cap[1,i)|}+sx_i,
      &i\notin A,\ i\in B,\\
    1,&i\in A,\ i\notin B,\\
    sx_i,&i\in A,\ i\in B.
  \end{cases}
\end{equation}
\end{proposition}

\begin{proof}
Fix $I\subseteq[n]$.
Equation~\eqref{eqRefinedS45InsertionCounts} gives
$\rho_{i^+}=\lambda_{\xi_i}(\alpha)+|I\cap(i,n]|$ and
$\rho_{i^-}=\lambda_{\xi_i}(\alpha)+|I\cap[i,n]|$.
In \eqref{eqRefinedS45PhiFactors}, the elementary symmetric function
chooses each alphabet entry at most once, while the complete homogeneous
function allows arbitrary multiplicity. Their generating products,
applied to the multiplicities of $i^+$ and $i^-$, respectively, give
\begin{align*}
  \sum_{\mu\geq0}
    (sx_iq^{\rho_{i^+}})^\mu
    \varphi_{\mu,0}(\rho_{i^+})
  =
  \prod_{h=1}^{\rho_{i^+}}(1+sx_iq^h), \quad
  \sum_{\mu\geq0}
    (-sx_iq^{\rho_{i^-}-1})^\mu
    \varphi_{\mu,1}(\rho_{i^-})
  =
  \prod_{h=0}^{\rho_{i^-}-1}(1+sx_iq^h)^{-1}.
\end{align*}
Since $\rho_{i^-}=\rho_{i^+}+\mathbf{1}_{\{i\in I\}}$,
we have
\[
  \frac{\prod_{h=1}^{\rho_{i^+}}(1+sx_iq^h)}
       {\prod_{h=0}^{\rho_{i^-}-1}(1+sx_iq^h)}
  =
  \begin{cases}
    \dfrac{1+sx_iq^{\lambda_{\xi_i}(\alpha)+|I\cap(i,n]|}}
          {1+sx_i},
      &i\notin I,\\[6pt]
    \dfrac{1}{1+sx_i},
      &i\in I.
  \end{cases}
\]

We now substitute Theorem~\ref{thmRootHeadInsertion}
into~\eqref{eqRefinedS6HDefinition} and expand $W_I(\alpha,E)$
using~\eqref{eqRefinedS45InsertionProduct}.
For fixed $I$, the numbers $\rho_d$ are independent of $E$,
so the sum over $E$ factors into the sums over the multiplicities
of $i^+$ and $i^-$ computed above.
Using also
$\rho_{\xi_i}=\lambda_{\xi_i}(\alpha)+|I\cap(i,n]|$
from~\eqref{eqRefinedS45InsertionCounts}, we obtain
\begin{align}
  \mathcal H_\alpha(ut;0,s)
  &=\sum_E s^{|E|}K_{\alpha,E}(ut,0)
    \notag\\
  &=\sum_{I\subseteq[n]}(ut)^{|I|}
    \prod_{i\in I}[\rho_{\xi_i}]_q
    \left(
      \sum_E
      \prod_{d\in\operatorname{supp}(E)}
      \bigl(s\wt(d)q^{\rho_d-1}\bigr)^{\mu_d}
      \varphi_{\mu_d,\epsilon_d}(\rho_d)
    \right)Q_I
    \notag\\
  &=\sum_{I\subseteq[n]}(ut)^{|I|}
    \prod_{i\in I}[\rho_{\xi_i}]_q
    \prod_{i=1}^{n}
    \frac{\prod_{h=1}^{\rho_{i^+}}(1+sx_iq^h)}
         {\prod_{h=0}^{\rho_{i^-}-1}(1+sx_iq^h)}
    Q_I
    \notag\\
  &=\frac{1}{\prod_{i=1}^{n}(1+sx_i)}
    \sum_{I\subseteq[n]}(ut)^{|I|}
    \prod_{i\in I}
      [\lambda_{\xi_i}(\alpha)+|I\cap(i,n]|]_q
    \prod_{i\notin I}
      \bigl(1+sx_iq^{\lambda_{\xi_i}(\alpha)+|I\cap(i,n]|}\bigr)
    Q_I.
    \label{eqNoRiseHeadExpansion}
\end{align}

Using~\eqref{lemEvaluationOrderedMonomials},
$(q-1)[m]_q=q^m-1$, and
$q^{\lambda_{\xi_i}(\alpha)}=e_{\{i\}}(\alpha)$ gives
\begin{align*}
  \left(\prod_{i=1}^{n}(1+sx_i)\right)
    \operatorname{ev}_A\mathcal H_\alpha(ut;0,s)
  =
  \sum_{I\subseteq[n]\setminus A}
  \prod_{i\in I}
    utx_iq^{|A\cap(i,n]|}
    \bigl(e_{\{i\}}(\alpha)q^{|I\cap(i,n]|}-1\bigr)
  \prod_{i\notin I}
    \bigl(1+sx_i e_{\{i\}}(\alpha)q^{|I\cap(i,n]|}\bigr).
\end{align*}

Expand the products and group the terms obtained by choosing
$e_{\{i\}}(\alpha)$ exactly for $i\in B$.
The chosen factors satisfy
\begin{align*}
  \prod_{b\in B}e_{\{b\}}(\alpha)
  =e_B(\alpha),
  \prod_{b\in B}q^{|I\cap(b,n]|}
  =q^{\sum_{b\in B}\sum_{i\in I}\mathbf{1}_{\{b<i\}}}
   =\prod_{i\in I}q^{|B\cap[1,i)|}.
\end{align*}
For $i\in I$, choosing the term containing
$e_{\{i\}}(\alpha)$ gives the positive contribution,
whereas choosing $-1$ gives the negative contribution.
For $i\notin I$, the corresponding choices contribute
$sx_i$ and $1$, after the displayed powers of $q$
have been assigned to the factors indexed by $I$.
Since $I\cap A=\varnothing$, we also have
$|A\cap(i,n]|=|A\cap[i,n]|$ for $i\in I$.
Thus
\begin{align*}
  &\left(\prod_{i=1}^{n}(1+sx_i)\right)
    \operatorname{ev}_A\mathcal H_\alpha(ut;0,s)\\
  &=\sum_{B\subseteq[n]}e_B(\alpha)
    \sum_{I\subseteq[n]\setminus A}
    \prod_{i\in I\cap B}
      utx_iq^{|A\cap[i,n]|+|B\cap[1,i)|}
    \prod_{i\in I\setminus B}
      \left(-utx_iq^{|A\cap[i,n]|+|B\cap[1,i)|}\right)
    \prod_{i\in B\setminus I}sx_i.
\end{align*}

For fixed $B$, each factor now depends only on whether its
index belongs to $I$.
If $i\notin A$, both choices are allowed; if $i\in A$,
only $i\notin I$ is allowed.
The resulting sums at each index are precisely the four cases
in~\eqref{eqRefinedS6VDefinition}.
The choices at distinct indices are independent, so
the sum over $I$ is $\prod_{i=1}^{n}v_i(B,A)$.
Consequently,
\begin{align*}
  \operatorname{ev}_A\mathcal H_\alpha(ut;0,s)
  =\frac{1}{\prod_{i=1}^{n}(1+sx_i)}
    \sum_{B\subseteq[n]}e_B(\alpha)
    \prod_{i=1}^{n}v_i(B,A)=\sum_{B\subseteq[n]}
    e_B(\alpha)V_{X_n}(ut,s)_{B,A}.
\end{align*}
\end{proof}

\subsection{The summation with rise decorations}\label{subsecRiseTransferMatrix}
Recall $X_n$ from Subsection~\ref{subsecNoRiseTransferMatrix}.
We also recall $\mathcal H_\alpha$, $e_B$, and $V_{X_n}(ut,s)$ from
\eqref{eqRefinedS6HDefinition}, \eqref{eqRefinedS45Evaluation}, and
\eqref{eqRefinedS6VDefinition}.
Let $\mathcal V_n$ be the vector space over
$\mathbb Q(q,t,x_1,\ldots,x_n)$ with basis
$z_A:=\prod_{i\in A}z_i$,
indexed by $A\subseteq[n]$,
where $z_1,\ldots,z_n$ are commuting indeterminates
and $z_\varnothing=1$.
These indeterminates are distinct from the Dyck path algebra
operators $z_i$ and from the elements $Z_{i^\pm}$ and $Q_i$
in~\eqref{eqRefinedS45ZRelations} and~\eqref{eqRefinedS45QGenerators}.

Define \begin{equation}\label{eqRefinedS6EpsilonDefinition}
\epsilon f:=f(1,\ldots,1)
\end{equation} for $f\in \mathcal V_n$.
Note that the functionals $\operatorname{ev}_A$ act on $\mathcal Q_n$,
whereas $\epsilon$ acts on $\mathcal V_n$.
Matrices below act on column vectors in the basis
$\{z_A:A\subseteq[n]\}$.
Matrix rows and columns are indexed by subsets of $[n]$;
the entry $(B,A)$ lies in row $B$ and column $A$.

We extend scalars to include the parameters of the operators below.
The coefficients in~\eqref{eqRefinedS6VDefinition} define
the action of $V_{X_n}(ut,s)$ by
\[
  V_{X_n}(ut,s)(z_A)
  :=
  \sum_{B\subseteq[n]}V_{X_n}(ut,s)_{B,A}z_B,
  \qquad A\subseteq[n],
\]
extended linearly.
We use the same convention for the other operators
given by coefficients indexed by subsets of $[n]$.

\begin{definition}
Define the following operators on $\mathcal V_n$:
\begin{equation}\label{eqRefinedS6ABCD}
\begin{alignedat}{2}
  \mathcal A z_A
  &=q^{|A|}z_A,
  &\qquad \mathcal B z_A
  &=(q-1)\sum_{i\in A}
    q^{|A\cap[1,i)|}z_{A\setminus\{i\}},\\
  \mathcal C
  &=\mathcal A-\mathcal B,
  &\qquad \mathcal D z_A
  &=(q-1)\sum_{i\notin A}
    x_iq^{|A\cap(i,n]|}z_{A\cup\{i\}}.
\end{alignedat}
\end{equation}
These operators depend on the ordered alphabet $X_n$;
we omit the alphabet subscript when it is clear from context.
\end{definition}

\begin{lemma}\label{lemTransferMatrixRelations}
\leavevmode
\begin{enumerate}
  \item
  We have
  $V_{X_n}(qut,0)
    =
    \mathcal C V_{X_n}(ut,0)
    =
    \mathcal A V_{X_n}(ut,0)(1-ut\mathcal D)$.

  \item
  There is a matrix $\mathsf J_{X_n}(s)$, independent of $ut$,
  such that
  \[
    V_{X_n}(ut,s)
    =
    \mathsf J_{X_n}(s)V_{X_n}(ut,0),
    \qquad
    \mathcal A\mathsf J_{X_n}(s)\mathcal A^{-1}
    =
    (1+s\mathcal D)\mathsf J_{X_n}(s).
  \]
\end{enumerate}
\end{lemma}

\begin{proof}
\textup{(1)}
By Proposition~\ref{propGIdentities}\textup{(2)},
$G_0(qS)^{-1}=(1-utS)G_0(S)^{-1}$.
For $0\leq m\leq n$, normal ordering $S^m$ gives a nonzero
sum of powers of $q$ as the coefficient of any $Q_I$ with
$|I|=m$.  The elements $1,S,\ldots,S^n$ have different
$Q$-degrees and are therefore linearly independent.
Comparing coefficients of $S^m$, each coefficient is
$ut/(1-q^m)$ times the preceding one.
Since the constant term is $1$ and $S^{n+1}=0$, this gives
\begin{equation}\label{eqRefinedS6G0InverseSeries}
  G_0(S)^{-1}
  =
  \sum_{m=0}^{n}\frac{(ut)^m}{(q;q)_m}S^m.
\end{equation}
In the expansion of $S^m$, products with repeated indices
vanish.  For distinct indices, inserting the largest index
into a word of length $h-1$ contributes
$1,q,\ldots,q^{h-1}$ according to its position.
Thus Proposition~\ref{propHeadLinearFunctional}\textup{(1)} and
Definition~\ref{defRefinedS45Quotient} give
\begin{equation}\label{eqRefinedS6PowerSumQ}
  S^m
  =
  \prod_{h=1}^{m}(1+q+\cdots+q^{h-1})
  \sum_{\substack{I\subseteq[n]\\|I|=m}}Q_I
  =
  \frac{(q;q)_m}{(1-q)^m}
  \sum_{\substack{I\subseteq[n]\\|I|=m}}Q_I.
\end{equation}
Substituting~\eqref{eqRefinedS6PowerSumQ} into~\eqref{eqRefinedS6G0InverseSeries} yields
\begin{equation}\label{eqRefinedS6G0Inverse}
  G_0(S)^{-1}
  =
  \sum_{m=0}^{n}\frac{(ut)^m}{(1-q)^m}
  \sum_{\substack{I\subseteq[n]\\|I|=m}}Q_I
  =
  \sum_{I\subseteq[n]}
  \frac{(-ut)^{|I|}}{(q-1)^{|I|}}Q_I.
\end{equation}

For $C\cap A=\varnothing$, multiply~\eqref{eqRefinedS6G0Inverse} on the left by
$(ut)^{|C|}q^{\binom{|C|}{2}}Q_C/(q-1)^{|C|}$.
Applying $\operatorname{ev}_A$ and
using~\eqref{lemEvaluationOrderedMonomials}, we obtain
\begin{align*}
  \operatorname{ev}_A\left(
    \frac{(ut)^{|C|}q^{\binom{|C|}{2}}}{(q-1)^{|C|}}
    Q_CG_0(S)^{-1}
  \right)
  =
  (ut)^{|C|}
  q^{\binom{|C|}{2}+\#\{(c,a)\in C\times A:c<a\}}x_C
  \prod_{i\notin A\cup C}
  \left(1-utx_iq^{|A\cap(i,n]|+|C\cap[1,i)|}\right).
\end{align*}
This is $V_{X_n}(ut,0)_{C,A}$ by
Proposition~\ref{propNoRiseTransferMatrix}.
When $C\cap A\neq\varnothing$, both expressions vanish.
Hence
\begin{equation}\label{eqNoRiseMatrixEvaluation}
  V_{X_n}(ut,0)_{C,A}
  =
  \operatorname{ev}_A\left(
    \frac{(ut)^{|C|}q^{\binom{|C|}{2}}}{(q-1)^{|C|}}
    Q_CG_0(S)^{-1}
  \right).
\end{equation}

We use~\eqref{eqNoRiseMatrixEvaluation} together with
$Q_BS=\sum_{i\notin B}q^{|B\cap[1,i)|}Q_{B\cup\{i\}}$
and $\binom{|B|+1}{2}=\binom{|B|}{2}+|B|$.
These identities give
\begin{alignat*}{2}
\bigl(\mathcal C V_{X_n}(ut,0)\bigr)_{B,A} &={} q^{|B|}V_{X_n}(ut,0)_{B,A} -(q-1)\sum_{i\notin B} q^{|B\cap[1,i)|}V_{X_n}(ut,0)_{B\cup\{i\},A}  &\quad\makebox[41.6pt][r]{\text{(by~\eqref{eqRefinedS6ABCD})}} &\\
&= q^{|B|}\operatorname{ev}_A\left( \frac{(ut)^{|B|}q^{\binom{|B|}{2}}}{(q-1)^{|B|}} Q_B(1-utS)G_0(S)^{-1} \right)  &\quad\makebox[41.6pt][r]{\text{(by \eqref{eqNoRiseMatrixEvaluation})}} &\\
&= q^{|B|}\operatorname{ev}_A\left( \frac{(ut)^{|B|}q^{\binom{|B|}{2}}}{(q-1)^{|B|}} Q_BG_0(qS)^{-1} \right)  &\quad\makebox[41.6pt][r]{\text{(by Proposition~\ref{propGIdentities}\textup{(2)})}} &\\
&= V_{X_n}(qut,0)_{B,A}  &\quad\makebox[41.6pt][r]{\text{(by \eqref{eqNoRiseMatrixEvaluation})}}\rlap{.} &
\end{alignat*}

Meanwhile, \eqref{eqRefinedS45HeadElements} and
\eqref{lemEvaluationOrderedMonomials}, together with
Proposition~\ref{propHeadLinearFunctional}\textup{(1)}
and Definition~\ref{defRefinedS45Quotient}, give
\begin{equation}\label{aboveeq}
  \operatorname{ev}_A(Q_CS)
  =\sum_{i=1}^{n}\operatorname{ev}_A(Q_CQ_i)
  =\sum_{i\notin C}
    q^{|C\cap[1,i)|}\operatorname{ev}_A(Q_{C\cup\{i\}})
  =(q-1)\sum_{i\notin A}
    x_iq^{|A\cap(i,n]|}\operatorname{ev}_{A\cup\{i\}}(Q_C).
\end{equation}
By linearity, \eqref{aboveeq} holds with $Q_C$
replaced by any element of $\mathcal Q_n$.
Consequently,
\begin{alignat*}{2}
  &\quad \bigl(
    \mathcal A V_{X_n}(ut,0)(1-ut\mathcal D)
  \bigr)_{B,A}\\
  &=
  q^{|B|}
  \left(
    V_{X_n}(ut,0)_{B,A}
    -ut(q-1)\sum_{i\notin A}
    x_iq^{|A\cap(i,n]|}V_{X_n}(ut,0)_{B,A\cup\{i\}}
  \right)
  &\qquad\text{(by~\eqref{eqRefinedS6ABCD})}&\\
  &=
  q^{|B|}\operatorname{ev}_A\left(
    \frac{(ut)^{|B|}q^{\binom{|B|}{2}}}{(q-1)^{|B|}}
    Q_BG_0(S)^{-1}(1-utS)
  \right)
  &\qquad\text{(by~\eqref{eqNoRiseMatrixEvaluation},
                   \eqref{aboveeq})}&\\
  &=
  q^{|B|}\operatorname{ev}_A\left(
    \frac{(ut)^{|B|}q^{\binom{|B|}{2}}}{(q-1)^{|B|}}
    Q_BG_0(qS)^{-1}
  \right)
  &\qquad\text{(by Proposition~\ref{propGIdentities}\textup{(2)})}&\\
  &=
  V_{X_n}(qut,0)_{B,A}
  &\qquad\text{(by~\eqref{eqNoRiseMatrixEvaluation},
                   \eqref{eqRefinedS6G0InverseSeries})}&\rlap{.}
\end{alignat*}
This proves \textup{(1)}.

\textup{(2)}
Since $S^{n+1}=0$, \eqref{eqdefGr},
\eqref{eqRefinedS6gDefinition}, and
\eqref{eqRefinedS6PowerSumQ} give
\begin{equation}\label{eqRefinedS6G0Expansion}
  G_0(S)
  =
  \sum_{m=0}^{n}
  \frac{(-ut)^mq^{\binom{m}{2}}}{(q;q)_m}S^m
  =
  \sum_{J\subseteq[n]}
  \frac{(ut)^{|J|}q^{\binom{|J|}{2}}}{(q-1)^{|J|}}Q_J.
\end{equation}
By Proposition~\ref{propHeadLinearFunctional}\textup{(1)}
and Definition~\ref{defRefinedS45Quotient},
\begin{equation}\label{eqRefinedS6OrderedProduct}
  Q_IQ_J
  =
  \begin{cases}
    q^{\#\{(i,j)\in I\times J:i<j\}}Q_{I\cup J},
      &I\cap J=\varnothing,\\
    0,&I\cap J\neq\varnothing.
  \end{cases}
\end{equation}
Thus, for every $C\subseteq[n]$,
\begin{equation}\label{eqRefinedS6LiftCoefficient}
  [Q_C]\bigl(\mathcal H_\alpha(ut;0,s)G_0(S)\bigr)
  =
  \sum_{I\subseteq C}
  \frac{(ut)^{|C\setminus I|}
    q^{\binom{|C\setminus I|}{2}
      +\#\{(i,j)\in I\times(C\setminus I):i<j\}}}
       {(q-1)^{|C\setminus I|}}
  [Q_I]\mathcal H_\alpha(ut;0,s).
\end{equation}

Regard $e_{\{1\}}(\alpha),\ldots,e_{\{n\}}(\alpha)$
as independent variables.
By~\eqref{eqNoRiseHeadExpansion}, after factoring
$(ut)^{|I|}$ out of $[Q_I]\mathcal H_\alpha(ut;0,s)$,
the coefficients are independent of $ut$, and each
$e_{\{i\}}(\alpha)$ occurs at most once in each term.
Thus $(ut)^{|C|}$ is a common factor
in~\eqref{eqRefinedS6LiftCoefficient}.
After factoring out
$(ut)^{|C|}q^{\binom{|C|}{2}}/(q-1)^{|C|}$
from its right-hand side, let
$\mathsf J_{X_n}(s)_{B,C}$ be the coefficient of
$e_B(\alpha)=\prod_{i\in B}e_{\{i\}}(\alpha)$.
These coefficients depend only on $q,s,x_1,\ldots,x_n$,
so we obtain
\begin{equation}\label{eqNoRiseTransferLift}
  \mathcal H_\alpha(ut;0,s)G_0(S)
  =
  \sum_{B,C\subseteq[n]}
  e_B(\alpha)\mathsf J_{X_n}(s)_{B,C}
  \frac{(ut)^{|C|}q^{\binom{|C|}{2}}}{(q-1)^{|C|}}Q_C.
\end{equation}

Multiply~\eqref{eqNoRiseTransferLift} on the right by
$G_0(S)^{-1}$ and apply $\operatorname{ev}_A$.
Proposition~\ref{propNoRiseTransferMatrix}
and~\eqref{eqNoRiseMatrixEvaluation} give
\begin{align*}
  \sum_{B\subseteq[n]}e_B(\alpha)V_{X_n}(ut,s)_{B,A}
  =\operatorname{ev}_A\mathcal H_\alpha(ut;0,s)
  =\sum_{B,C\subseteq[n]}
    e_B(\alpha)\mathsf J_{X_n}(s)_{B,C}
    V_{X_n}(ut,0)_{C,A}.
\end{align*}

Comparing coefficients of $e_B(\alpha)$ therefore gives
\begin{equation}\label{eqRefinedS6VFactorization}
  V_{X_n}(ut,s)_{B,A}
  =
  \sum_{C\subseteq[n]}
  \mathsf J_{X_n}(s)_{B,C}V_{X_n}(ut,0)_{C,A}.
\end{equation}

For the second identity, the four cases
in~\eqref{eqRefinedS6VDefinition} give
\begin{equation}\label{eqLocalFactorShift}
\begin{aligned}
  q^{\mathbf1_{\{i\in B\}}}v_i(B,A)
  -\left.v_i(B,A)\right|_{ut\mapsto qut}
  &=(q-1)sx_i\mathbf1_{\{i\in B\}}
    \left.v_i(B\setminus\{i\},A)\right|_{ut\mapsto qut}\\
  &\quad+(q-1)utx_i
    q^{|A\cap(i,n]|+|B\cap[1,i]|}
    \mathbf1_{\{i\notin A\}}v_i(B,A\cup\{i\}).
\end{aligned}
\end{equation}
For $k\neq i$, \eqref{eqRefinedS6VDefinition} gives
\begin{equation}\label{theseidentities}
  \begin{aligned}
  \left.v_k(B\setminus\{i\},A)\right|_{ut\mapsto qut}
  &=
  \begin{cases}
    \left.v_k(B,A)\right|_{ut\mapsto qut},&k<i,\\
    v_k(B,A),&k>i,
  \end{cases}
  \qquad(i\in B),\\
  v_k(B,A\cup\{i\})
  &=
  \begin{cases}
    \left.v_k(B,A)\right|_{ut\mapsto qut},&k<i,\\
    v_k(B,A),&k>i,
  \end{cases}
  \qquad(i\notin A).
\end{aligned}
\end{equation}
Replacing the factors in increasing order of $i$, and then
using~\eqref{eqLocalFactorShift} and \eqref{theseidentities}, yields
\begin{align*}
  &\quad q^{|B|}\prod_{k=1}^{n}v_k(B,A)
   -\left.\prod_{k=1}^{n}v_k(B,A)\right|_{ut\mapsto qut}\\
  &=
  \sum_{i=1}^{n}q^{|B\cap(i,n]|}
  \left(\left.\prod_{k<i}v_k(B,A)\right|_{ut\mapsto qut}\right)
  \left(
    q^{\mathbf1_{\{i\in B\}}}v_i(B,A)
    -\left.v_i(B,A)\right|_{ut\mapsto qut}
  \right)
  \prod_{k>i}v_k(B,A)\\
  &=
  s(q-1)\sum_{i\in B}x_iq^{|B\cap(i,n]|}
  \left.\prod_{k=1}^{n}
    v_k(B\setminus\{i\},A)\right|_{ut\mapsto qut}
  +
  utq^{|B|}(q-1)\sum_{i\notin A}x_iq^{|A\cap(i,n]|}
  \prod_{k=1}^{n}v_k(B,A\cup\{i\}).
\end{align*}
In the first equality, consecutive intermediate products
cancel. In the last term, the powers of $q$ combine by
$|B\cap[1,i]|+|B\cap(i,n]|=|B|$.

Dividing by $\prod_{k=1}^{n}(1+sx_k)$ and
using~\eqref{eqRefinedS6VDefinition}, we obtain
\begin{align*}
  q^{|B|}V_{X_n}(ut,s)_{B,A}
  &=
  V_{X_n}(qut,s)_{B,A}
  +s(q-1)\sum_{i\in B}x_iq^{|B\cap(i,n]|}
    V_{X_n}(qut,s)_{B\setminus\{i\},A}\\
  &\quad+
  utq^{|B|}(q-1)\sum_{i\notin A}x_iq^{|A\cap(i,n]|}
    V_{X_n}(ut,s)_{B,A\cup\{i\}}.
\end{align*}
By~\eqref{eqRefinedS6ABCD}, this is
\begin{equation}\label{eqRefinedS6VShift}
  \mathcal A V_{X_n}(ut,s)(1-ut\mathcal D)
  =
  (1+s\mathcal D)V_{X_n}(qut,s).
\end{equation}

By~\eqref{eqRefinedS6VDefinition},
$V_{X_n}(ut,0)_{B,[n]\setminus A}=0$ unless $B\subseteq A$.
Thus, after complementing the column indices and ordering
subsets by increasing cardinality, the matrix is upper
triangular, with diagonal entries
\[
  V_{X_n}(ut,0)_{B,[n]\setminus B}
  =
  (ut)^{|B|}x_B
  q^{\binom{|B|}{2}
     +\#\{(i,a)\in B\times([n]\setminus B):i<a\}}
  \neq0.
\]
Hence $V_{X_n}(ut,0)$ is invertible over
$\mathbb Q(q,t,x_1,\ldots,x_n)(u)$.
Part~\textup{(1)} gives
\[
  V_{X_n}(ut,0)^{-1}\mathcal A^{-1}
  =
  (1-ut\mathcal D)V_{X_n}(qut,0)^{-1}.
\]
Consequently,
\begin{alignat*}{2}
  \mathcal A\mathsf J_{X_n}(s)\mathcal A^{-1}
  &=
  \mathcal A V_{X_n}(ut,s)V_{X_n}(ut,0)^{-1}\mathcal A^{-1}
  &\qquad\text{(by~\eqref{eqRefinedS6VFactorization})}&\\
  &=
  \mathcal A V_{X_n}(ut,s)(1-ut\mathcal D)
  V_{X_n}(qut,0)^{-1}
  &\qquad\text{(by Lemma \ref{lemTransferMatrixRelations} \textup{(1)})}&\\
  &=
  (1+s\mathcal D)V_{X_n}(qut,s)V_{X_n}(qut,0)^{-1}
  &\qquad\text{(by~\eqref{eqRefinedS6VShift})}&\\
  &=
  (1+s\mathcal D)\mathsf J_{X_n}(s)
  &\qquad\text{(by~\eqref{eqRefinedS6VFactorization})}&\rlap{.}
\end{alignat*}
\end{proof}

For a scalar parameter $z$, define
  \begin{equation}\label{eqRefinedS6TDefinition}T_n(z)_{B,I}
  :=
  \mathbf 1_{\{B\subseteq I\}}
  (-z)^{|I\setminus B|}
  q^{-\#\{(i,b)\in(I\setminus B)\times B:i<b\}}
  \prod_{h=0}^{|B|-1}(1+zq^{-h}),\end{equation}
  \begin{equation}\label{eqRefinedS6RDefinition}R_{X_n}(ut,r)_{C,A}
  :=
  \operatorname{ev}_A\left(
    \frac{g_{|C|}(ut,r)}{(q-1)^{|C|}}
    Q_CG_r(S)^{-1}
  \right).\end{equation}

\begin{lemma}\label{lemRiseTransferFactor}
We have
\[
  R_{X_n}(ut,r)
  =
  T_n\left(\frac{r}{ut}\right)V_{X_n}(ut,0).
\]
The product on the right is polynomial in $ut,r$.
\end{lemma}

\begin{proof}
Applying Proposition~\ref{propGIdentities}\textup{(2)}
with parameters $r$ and $0$ gives
\begin{equation}\label{eqRefinedS6GRRatioShift}
  G_r(S)^{-1}G_0(S)
  =
  (1+rS)G_r(qS)^{-1}G_0(qS).
\end{equation}
By~\eqref{eqRefinedS6PowerSumQ},
$1,S,\ldots,S^n$ are linearly independent.
Comparing coefficients in~\eqref{eqRefinedS6GRRatioShift}
therefore gives, for $1\leq j\leq n$,
\[
  (1-q^j)[S^j]\bigl(G_r(S)^{-1}G_0(S)\bigr)
  =
  rq^{j-1}[S^{j-1}]\bigl(G_r(S)^{-1}G_0(S)\bigr).
\]
Starting from the constant term $1$, we obtain
\[
  [S^j]\bigl(G_r(S)^{-1}G_0(S)\bigr)
  =
  \prod_{h=1}^{j}\frac{rq^{h-1}}{1-q^h}
  =
  \frac{r^jq^{\binom{j}{2}}}{(q;q)_j}.
\]
Thus~\eqref{eqRefinedS6PowerSumQ} yields
\begin{equation}\label{eqRefinedS6GRRatio}
  G_r(S)^{-1}G_0(S)
  =
  \sum_{j=0}^{n}
  \frac{r^jq^{\binom{j}{2}}}{(q;q)_j}S^j
  =
  \sum_{I\subseteq[n]}
  \frac{(-r)^{|I|}q^{\binom{|I|}{2}}}{(q-1)^{|I|}}Q_I.
\end{equation}

Substituting~\eqref{eqRefinedS6GRRatio}
into~\eqref{eqRefinedS6RDefinition} and
using~\eqref{eqRefinedS6OrderedProduct}, we obtain
\begin{align*}
  R_{X_n}(ut,r)_{C,A}
  &=
  \operatorname{ev}_A\left(
    \frac{g_{|C|}(ut,r)}{(q-1)^{|C|}}
    Q_C\bigl(G_r(S)^{-1}G_0(S)\bigr)G_0(S)^{-1}
  \right)\\
  &=
  \sum_{I\supseteq C}
  \frac{
    g_{|C|}(ut,r)(-r)^{|I\setminus C|}
    q^{\binom{|I\setminus C|}{2}
      +\#\{(c,i)\in C\times(I\setminus C):c<i\}}
  }{(q-1)^{|I|}}
  \operatorname{ev}_A\bigl(Q_IG_0(S)^{-1}\bigr).
\end{align*}

For $C\subseteq I$, counting pairs within $C$, within
$I\setminus C$, and between these two sets gives
\[
  \binom{|I|}{2}
  =
  \binom{|C|}{2}+\binom{|I\setminus C|}{2}
  +\#\{(c,i)\in C\times(I\setminus C):c<i\}
  +\#\{(i,c)\in(I\setminus C)\times C:i<c\}.
\]
Consequently, \eqref{eqRefinedS6gDefinition}
and~\eqref{eqRefinedS6TDefinition} give
\begin{align*}
  \frac{
    g_{|C|}(ut,r)(-r)^{|I\setminus C|}
    q^{\binom{|I\setminus C|}{2}
      +\#\{(c,i)\in C\times(I\setminus C):c<i\}}
  }{(ut)^{|I|}q^{\binom{|I|}{2}}}
  &=
  \left(-\frac{r}{ut}\right)^{|I\setminus C|}
  q^{-\#\{(i,c)\in(I\setminus C)\times C:i<c\}}
  \prod_{h=0}^{|C|-1}
  \left(1+\frac{r}{ut}q^{-h}\right)\\
  &=T_n\left(\frac{r}{ut}\right)_{C,I}.
\end{align*}
Using~\eqref{eqNoRiseMatrixEvaluation}, we conclude that
\[
  R_{X_n}(ut,r)_{C,A}
  =
  \sum_{I\supseteq C}
  T_n\left(\frac{r}{ut}\right)_{C,I}
  V_{X_n}(ut,0)_{I,A}.
\]

Finally, \eqref{eqRefinedS6TDefinition} shows that
$T_n(z)_{C,I}$ has degree at most $|I|$ in $z$.
By~\eqref{eqNoRiseMatrixEvaluation} and
\eqref{eqRefinedS6G0InverseSeries},
$V_{X_n}(ut,0)_{I,A}$ is a polynomial in $ut$
divisible by $(ut)^{|I|}$.
Hence every summand in the matrix product is polynomial
in $ut,r$, as required.
\end{proof}

We can now extend Proposition~\ref{propNoRiseTransferMatrix} to allow decorated rises.
\begin{theorem}\label{thmRiseTransferMatrix}
For every word $\alpha$ as in
Proposition~\ref{propNoRiseTransferMatrix}
and every $A\subseteq[n]$, we have
\[
  \operatorname{ev}_A\mathcal H_\alpha(ut;r,s)
  =
  \sum_{B\subseteq[n]}
  e_B(\alpha)M_{X_n}(ut,r,s)_{B,A},
\]
where
\begin{equation}\label{eqRefinedS6MDefinition}
  M_{X_n}(ut,r,s)
  :=
  \mathsf J_{X_n}(s)R_{X_n}(ut,r)
  =
  \mathsf J_{X_n}(s)
  T_n\left(\frac{r}{ut}\right)V_{X_n}(ut,0).
\end{equation}
\end{theorem}

\begin{proof}
Recall from~\eqref{eqRefinedS6HDefinition} that
$\mathcal H_\alpha(ut;r,s)=\sum_Es^{|E|}K_{\alpha,E}(ut,r)$.
Each coefficient in $s$ is a finite sum, as noted after
\eqref{eqRefinedS6HDefinition}. Applying
Theorem~\ref{thmRiseDifference} coefficient by coefficient gives
\begin{alignat*}{2}
\mathcal H_\alpha(ut;r,s)G_r(S) &= \sum_{j=0}^{n}q^{-\binom{j}{2}}g_j(ut,r) \left. [\mathcal H_\alpha(ut;0,s)G_0(S)]_j \right|_{ut=1}  &\quad\makebox[78.7pt][r]{\text{(by Theorem~\ref{thmRiseDifference})}} &\\
&= \sum_{B,C\subseteq[n]} e_B(\alpha)\mathsf J_{X_n}(s)_{B,C} \frac{g_{|C|}(ut,r)}{(q-1)^{|C|}}Q_C  &\quad\makebox[78.7pt][r]{\text{(by~\eqref{eqNoRiseTransferLift})}}\rlap{.} &
\end{alignat*}
In the second equality, $\mathsf J_{X_n}(s)$ is independent
of $ut,r$.

Multiplying on the right by $G_r(S)^{-1}$ and applying
$\operatorname{ev}_A$, we obtain
\begin{alignat*}{2}
\operatorname{ev}_A\mathcal H_\alpha(ut;r,s) &= \sum_{B,C\subseteq[n]} e_B(\alpha)\mathsf J_{X_n}(s)_{B,C} \operatorname{ev}_A\left( \frac{g_{|C|}(ut,r)}{(q-1)^{|C|}} Q_CG_r(S)^{-1} \right)\\
&= \sum_{B,C\subseteq[n]} e_B(\alpha)\mathsf J_{X_n}(s)_{B,C} R_{X_n}(ut,r)_{C,A}  &\quad\makebox[0.0pt][r]{\text{(by~\eqref{eqRefinedS6RDefinition})}} &\\
&= \sum_{B\subseteq[n]} e_B(\alpha)M_{X_n}(ut,r,s)_{B,A}  &\quad\makebox[0.0pt][r]{\text{(by~\eqref{eqRefinedS6MDefinition})}}\rlap{.} &
\end{alignat*}
Finally, Lemma~\ref{lemRiseTransferFactor} gives
\begin{equation*}
  M_{X_n}(ut,r,s)
  =
  \mathsf J_{X_n}(s)R_{X_n}(ut,r)
  =
  \mathsf J_{X_n}(s)
  T_n\left(\frac{r}{ut}\right)V_{X_n}(ut,0).\qedhere
\end{equation*}
\end{proof}

\begin{example}
Recall $\alpha=(1^-,2^+,6^+)$ from
Example~\ref{exDeletingMainDiagonal}.
Take $n=8$ and $A=[8]\setminus\{1\}$.
We calculate the coefficient of $sx_8$ on both sides of
Theorem~\ref{thmRiseTransferMatrix}.

By~\eqref{lemEvaluationOrderedMonomials},
$\operatorname{ev}_A(Q_I)=0$ unless
$I=\varnothing$ or $I=\{1\}$.
Thus only $E=\{8^-\}$ and $E=\{8^+\}$ contribute to this
coefficient.
For either sign, inserting the decorated valley after
$1^-$, after $2^+$, or after $6^+$ gives factors $1,q,q^2$.
If a head of absolute value $1$ is present, it must immediately
follow $1^-$; the corresponding values of $\nu_1$ are
$2,-1,-1$.
Since $\wt(8^-)+\wt(8^+)=(q-1)x_8$ and
$\operatorname{ev}_A(Q_1)=(q-1)q^7x_1$, \eqref{eqRefinedS6KDefinition}
and~\eqref{eqRefinedS6HDefinition} give
\begin{align*}
  [sx_8]\operatorname{ev}_A\mathcal H_\alpha(ut;r,s)
  =
  (q^3-1)\bigl(1+(q-1)q^7(ut+r)x_1\bigr).
\end{align*}

For the right-hand side, \eqref{eqRefinedS6RDefinition} gives
\[
  R_{X_8}(ut,r)_{\varnothing,A}
  =1-(ut+r)q^7x_1,
  \qquad
  R_{X_8}(ut,r)_{\{1\},A}
  =(ut+r)q^7x_1,
\]
and every other entry in this column is zero.
Using Lemma~\ref{lemTransferMatrixRelations}\textup{(2)}, compare the
coefficients of $sx_8$ in
$V_{X_8}(ut,s)=\mathsf J_{X_8}(s)V_{X_8}(ut,0)$.
The factors $v_i$ in Proposition~\ref{propNoRiseTransferMatrix} give
\[
  [sx_8]\mathsf J_{X_8}(s)z_\varnothing
  =z_{\{8\}}-z_\varnothing,
  \qquad
  [sx_8]\mathsf J_{X_8}(s)z_{\{1\}}
  =z_{\{1,8\}}-z_{\{1\}}.
\]
Here $e_{\{1\}}(\alpha)=q$, $e_{\{8\}}(\alpha)=q^3$,
and $e_{\{1,8\}}(\alpha)=q^4$.
Therefore,
\begin{alignat*}{2}
[sx_8]\sum_{B\subseteq[8]} e_B(\alpha)M_{X_8}(ut,r,s)_{B,A} &= (q^3-1)\bigl(1-(ut+r)q^7x_1\bigr) +(q^4-q)(ut+r)q^7x_1\\
&= (q^3-1)\bigl(1+(q-1)q^7(ut+r)x_1\bigr).
\end{alignat*}
The two coefficients agree.
\end{example}

\section{The combinatorial recurrence}\label{secCombinatorialRecurrence}
In this section, we prove
Theorem~\ref{thmIntroCombinatorialExpansion} and complete the
proof of Theorem~\ref{thmIntroMain}.
We first fix the labels on the undecorated roots and sum
over the corresponding parking functions, using the results
of Section~\ref{secInsertionsWithRise}.
We then sum over these root labels to obtain the expansion
of $U_qR$, as stated in the Introduction.
By studying exchanges of adjacent letters in the alphabet,
we derive the recurrence for $\mathcal R_I$ in
Theorem~\ref{thmMainCombinatorialRecursion}.
Finally, we compare the recurrence for $\mathcal L_I$ in
Theorem~\ref{thmAlgMainRecursion} with that for $\mathcal R_I$
in Theorem~\ref{thmMainCombinatorialRecursion}.
Their common initial condition gives $\mathcal L_I=\mathcal R_I$.
Using the expansions in Theorems~\ref{thmAlgMainRecursion}
and~\ref{thmCombinatorialMarkedExpansion} and the invertibility of $U_q$
from Proposition~\ref{thmUqProperties}\textup{(2)},
we obtain $R=1+L$, proving Theorem~\ref{thmIntroMain}.
We then deduce Corollary~\ref{corThetaDeltaConsequences}.

\subsection{Path generating functions}\label{subsecNormalizedPathWeights}
Fix the finite alphabet $X_n=x_1+\cdots+x_n$.
Recall the signed labels $1^-,1^+,\ldots,n^-,n^+$ and
their weights from Definition~\ref{defRefinedS13SignedWeights}.
All generating functions below are formal power series in
$u,r,s$ with coefficients in $\mathbb Q(q,t)[x_1,\ldots,x_n]$.

For $\pi\in\PF_{d+k+\ell,\pm}^{\ast k,\bullet\ell}(X_n)$
with labels $(w_i)$, its contribution is
\[
  \wt(\pi)
  :=
  u^dr^ks^\ell q^{\dinv(\pi)}t^{\area(\pi)}
  \prod_{i=1}^{d+k+\ell}\wt(w_i).
\]
Thus $u$ counts north steps carrying neither decoration,
while $r$ and $s$ count decorated rises and decorated valleys,
respectively.
We include the empty path with weight $1$.

Fix a finite signed word $\alpha$.
We sum over all signed parking functions, of arbitrary size
and with arbitrary decorations, whose roots outside $\DD(\pi)$
have labels $\alpha$ in their order along the path.
Define
\begin{equation}\label{eqRefinedS7PathGeneratingFunction}
  F_\alpha(u,r,s)
  :=
  \sum_{\substack{
    \pi\in\PF_{d+k+\ell,\pm}^{\ast k,\bullet\ell}(X_n)\\
    (w_i:\,a_i=0,\ i\notin\DD(\pi))=\alpha
  }}
  \frac{\wt(\pi)}
       {u^{|\alpha|}q^{I(\alpha)}
        \displaystyle\prod_{j=1}^{|\alpha|}\wt(\alpha_j)}.
\end{equation}
The division is taken termwise.
It removes the row and label weights of these roots,
together with their mutual inversion factor.
All other weights remain.
For fixed $\alpha$, the coefficient of $u^dr^ks^\ell$ involves paths
with exactly $|\alpha|+d+k+\ell$ north steps. Thus the termwise
normalization still gives a well-defined formal power series.
Since the first north step of every nonempty path is an
undecorated root, $F_\varnothing(u,r,s)=1$.

For finite signed words $\alpha,\beta$, consider the choices
of $W,\delta,h$ satisfying the conditions of
Proposition~\ref{propMainDiagonalDeletion} with $\chi=0$,
where
$(c_j:\delta_j=0)=\alpha,
  h_1+\cdots+h_L=|\beta|$,
and $\beta$ is the ordered list of head labels.
Define
\begin{equation}\label{eqCombinatorialTransferWeight}
  \mathcal W_{\alpha,\beta}(ut,r,s)
  :=
  q^{I(\beta)}
  \left(\prod_{b=1}^{|\beta|}\wt(\beta_b)\right)
  \sum_{W,\delta,h}
  q^{\dinv_0(W)-I(\alpha)+C(W;\delta,h;\beta)}
  \prod_{j:\delta_j=1}s\wt(c_j)
  \prod_{j:h_j>0}
  (ut)^{h_j-1}(ut+rq^{\nu_j}).
\end{equation}
For each root assigned at least one head, the first head
may be undecorated or a decorated rise.
These choices contribute $ut$ and $rq^{\nu_j}$, respectively;
the other $h_j-1$ heads each contribute $ut$.
The factor preceding the sum restores the label weights
and mutual inversion factor of $\beta$.
The empty data contribute $1$.

\begin{example}
Recall the paths in Example~\ref{exDeletingMainDiagonal}.
The words of undecorated root labels of $\pi$ and $\pi'$ are
$\alpha=(1^-,2^+,6^+)$ and $\beta=(1^-,3^+)$, respectively.
Since $I(\alpha)=3$ and $I(\beta)=1$, their contributions to
$F_\alpha(u,r,s)$ and $F_\beta(u,r,s)$ are
\begin{alignat*}{2}
  \frac{\wt(\pi)} {u^3q^3(-x_1)(qx_2)(qx_6)} &= ur^3s^3q^{12}t^3x_1x_3x_4x_5x_7^2x_8, 
  \quad
  \frac{\wt(\pi')} {u^2q(-x_1)(qx_3)} = r^2s^2q^8x_4x_5x_7^2.
\end{alignat*}
For the deletion data in Example~\ref{exDeletingMainDiagonal},
$\dinv_0(W)=4$, $C(W;\delta,h;\beta)=0$,
$h=(1,1,0,0)$, and $(\nu_1,\nu_2)=(-1,1)$.
Their contribution to~\eqref{eqCombinatorialTransferWeight} is
$q^3s x_1x_3x_8(ut+rq^{-1})(ut+rq)$.
The displayed choice $\chi=(0,1,0,0)$ contributes
$q^4utrs\,x_1x_3x_8$.
Multiplying it by the contribution of $\pi'$ evaluated at
$(ut,r,ts)$ gives
\[
  (q^4utrs\,x_1x_3x_8)
  (r^2t^2s^2q^8x_4x_5x_7^2)
  =
  ur^3s^3q^{12}t^3x_1x_3x_4x_5x_7^2x_8,
\]
which is the contribution of $\pi$ to $F_\alpha(u,r,s)$.
\end{example}

\begin{proposition}\label{propNormalizedPathRecursion}
For every finite signed word $\alpha$,
\begin{equation}\label{eqNormalizedPathRecursion}
  F_\alpha(u,r,s)
  =
  \sum_\beta
  \mathcal W_{\alpha,\beta}(ut,r,s)F_\beta(ut,r,ts),
\end{equation}
where $\beta$ ranges over all finite signed words.
Together with $F_\varnothing=1$, these identities uniquely
determine the family $F_\alpha$.
\end{proposition}

\begin{proof}
By Proposition~\ref{propMainDiagonalDeletion}, every path
counted by $F_\alpha$ is obtained uniquely from a path $\pi'$
and admissible deletion data.
Let $\beta$ be the ordered list of labels of the undecorated
roots of $\pi'$, which are its heads.

Fix such a path $\pi'$.
The weights of its heads are omitted in $F_\beta$.
For every other north step of $\pi'$, reconstruction increases
its area entry by one.
An undecorated step with area entry $a_i$ therefore contributes
$(ut)t^{a_i}$ instead of $ut^{a_i}$, and a decorated valley
contributes $(ts)t^{a_i}$ instead of $st^{a_i}$.
A decorated rise still contributes $r$.
Replacing $u$ by $ut$ and $s$ by $ts$ in $F_\beta$ makes
exactly these changes, while the head weights remain omitted.

The remaining weights are those of the heads and the inserted
roots, together with the change in $\dinv$.
The weights removed for the word $\alpha$ remain excluded,
while those removed for $\beta$ must be restored.
By Theorem~\ref{thmMainDiagonalStatistics} and
Proposition~\ref{propFiniteHeadFactor}, summing these remaining
factors over the admissible deletion data and the choices
of decorated rises gives exactly
$\mathcal W_{\alpha,\beta}(ut,r,s)$.
Thus, for fixed $\beta$, the reconstructed paths contribute
$\mathcal W_{\alpha,\beta}(ut,r,s)F_\beta(ut,r,ts)$.
Summing over $\beta$ proves~\eqref{eqNormalizedPathRecursion}.

For uniqueness, each summand of $\mathcal W_{\alpha,\beta}$
has total degree $|\beta|+\sum_j\delta_j$ in $u,r,s$,
so only finitely many words and deletion data contribute
in any fixed degree.
When $\beta\neq\varnothing$, this degree is positive;
hence a coefficient of $F_\alpha$ depends only on
coefficients of lower total degree in the other $F_\beta$.
The term with $\beta=\varnothing$ is determined by
$F_\varnothing=1$.
Induction on total degree proves uniqueness.
\end{proof}

Recall the notation $e_A(\alpha)$ from~\eqref{eqRefinedS45Evaluation},
with occurrences in $\alpha$ counted with multiplicity.
The following proposition relates
$\mathcal W_{\alpha,\beta}$ to the matrix in
Theorem~\ref{thmRiseTransferMatrix}.

\begin{proposition}\label{propCombinatorialTransfer}
For every finite signed word $\alpha$ and every $A\subseteq[n]$,
\[
  \sum_\beta
  \mathcal W_{\alpha,\beta}(ut,r,s)e_A(\beta)
  =
  \sum_{B\subseteq[n]}
  e_B(\alpha)M_{X_n}(ut,r,s)_{B,A}.
\]
\end{proposition}

\begin{proof}
By~\eqref{eqRefinedS45QGenerators} and~\eqref{eqRefinedS45Evaluation}, we have
\[
  \operatorname{ev}_A(Q_{i_1}\cdots Q_{i_m})
  =
  \sum_{\substack{
    \beta=(\beta_1,\ldots,\beta_m)\\
    \beta_b\in\{i_b^-,i_b^+\}
  }}
  q^{I(\beta)}e_A(\beta)\prod_{b=1}^{m}\wt(\beta_b).
\]
Thus applying $\operatorname{ev}_A$ sums over the two possible
signs of each head label, with precisely the label and
inversion factors in
$\mathcal W_{\alpha,\beta}e_A(\beta)$.
These sign choices do not change the admissible insertions
or their comparisons with the deleted roots, since
$c\prec i^-$ if and only if $c\prec i^+$.

For a fixed multiset $E$, equation~\eqref{eqRefinedS6KDefinition}
sums the remaining factors in
\eqref{eqCombinatorialTransferWeight}, except for
$s^{|E|}$.
Multiplying by $s^{|E|}$ and summing over $E$ therefore gives
the first equality below:
\begin{alignat*}{2}
  \sum_\beta \mathcal W_{\alpha,\beta}(ut,r,s)e_A(\beta) &= \sum_E s^{|E|} \operatorname{ev}_A\bigl(K_{\alpha,E}(ut,r)\bigr)\\
   &= \operatorname{ev}_A\mathcal H_\alpha(ut;r,s) &\qquad\makebox[81.68pt][r]{\text{(by~\eqref{eqRefinedS6HDefinition})}}&\\
   &= \sum_{B\subseteq[n]} e_B(\alpha)M_{X_n}(ut,r,s)_{B,A} &\qquad\makebox[81.68pt][r]{\text{(by Theorem~\ref{thmRiseTransferMatrix})}}\rlap{.}&
\end{alignat*}
\end{proof}

\subsection{The summation over the undecorated root labels}\label{subsecSummingRootLabels}
Recall $M_{X_n}(ut,r,s)$ from Theorem~\ref{thmRiseTransferMatrix},
and the notation $\mathcal V_n$, $z_A$, and $\epsilon$ from \eqref{eqRefinedS6EpsilonDefinition}.

\begin{lemma}\label{lemTransferNormalization}
We have
\[
  \epsilon M_{X_n}(ut,r,s)=\epsilon,
  \qquad
  M_{X_n}(0,0,0)z_A=z_\varnothing
  \quad(A\subseteq[n]).
\]
\end{lemma}

\begin{proof}
If $\alpha=\varnothing$, there are no undecorated deleted roots.
Since the first root of every nonempty path is undecorated,
only the empty deletion data are admissible.
Thus $\mathcal W_{\varnothing,\varnothing}=1$ and
$\mathcal W_{\varnothing,\beta}=0$ for $\beta\neq\varnothing$.
Using $e_B(\varnothing)=1$ from~\eqref{eqRefinedS45Evaluation}, we obtain
\begin{alignat*}{2}
  \epsilon M_{X_n}(ut,r,s)z_A &= \sum_{B\subseteq[n]}M_{X_n}(ut,r,s)_{B,A} &\qquad\makebox[88.26pt][r]{\text{(by~\eqref{eqRefinedS6EpsilonDefinition})}}&\\
   &= \sum_\beta \mathcal W_{\varnothing,\beta}(ut,r,s)e_A(\beta) &\qquad\makebox[88.26pt][r]{\text{(by Proposition~\ref{propCombinatorialTransfer})}}&\\
   &= \mathcal W_{\varnothing,\varnothing}(ut,r,s)e_A(\varnothing) =1=\epsilon z_A.
\end{alignat*}
This proves the first identity.

For the second identity, setting $s=0$
in~\eqref{eqRefinedS6VFactorization} and using the invertibility
of $V_{X_n}(ut,0)$ gives $\mathsf J_{X_n}(0)=1$.
Moreover, \eqref{eqRefinedS6gDefinition} and~\eqref{eqdefGr} give
\[
  g_m(0,0)=\mathbf1_{\{m=0\}},
  \qquad
  \left.G_r(S)\right|_{ut=r=0}
  =\left.G_r(S)^{-1}\right|_{ut=r=0}=1.
\]
Therefore,

\begin{alignat*}{2}
  M_{X_n}(0,0,0)z_A &= \mathsf J_{X_n}(0)R_{X_n}(0,0)z_A &\qquad\makebox[81.68pt][r]{\text{(by~\eqref{eqRefinedS6MDefinition})}}&\\
   &= \sum_{C\subseteq[n]} \operatorname{ev}_A\left( \frac{g_{|C|}(0,0)}{(q-1)^{|C|}}Q_C \right)z_C &\qquad\makebox[81.68pt][r]{\text{(by~\eqref{eqRefinedS6RDefinition})}}&\\
   &= \sum_{C\subseteq[n]} \mathbf1_{\{C=\varnothing\}}\operatorname{ev}_A(1)z_C &\qquad\makebox[81.68pt][r]{\text{(by~\eqref{eqRefinedS6gDefinition},
          \eqref{eqRefinedS45OrderedQ})}}&\\
   &=z_\varnothing &\qquad\makebox[81.68pt][r]{\text{(by~\eqref{lemEvaluationOrderedMonomials})}}\rlap{.}&
\end{alignat*}
\end{proof}

Write
$\mathbf b_{X_n}(ut,r,s)
=\sum_{B\subseteq[n]}b_{X_n,B}(ut,r,s)z_B$,
where the coefficients are formal power series in $u,r,s$.
We determine these coefficients by the following equations.

\begin{proposition}\label{propNormalizedTransferVector}
The equations
\begin{equation}\label{eqNormalizedTransferVector}
  \mathbf b_{X_n}(ut,r,s)
  =
  M_{X_n}(ut,r,s)\mathbf b_{X_n}(ut^2,r,ts),
  \qquad
  \epsilon\mathbf b_{X_n}(ut,r,s)=1
\end{equation}
have a unique formal solution.
Its constant term in $u,r,s$ is $z_\varnothing$.
\end{proposition}

\begin{proof}
We construct $\mathbf b_{X_n}$ from the recursion
\[
  \mathbf b_{X_n}(ut,r,s)
  =
  z_\varnothing+
  \bigl(M_{X_n}(ut,r,s)-M_{X_n}(0,0,0)\bigr)
  \mathbf b_{X_n}(ut^2,r,ts).
\]
The difference $M_{X_n}(ut,r,s)-M_{X_n}(0,0,0)$ has
zero constant term, and the substitution
$(u,r,s)\mapsto(ut,r,ts)$ preserves total degree.
Thus the constant term of $\mathbf b_{X_n}$ is
$z_\varnothing$, and each coefficient of positive total degree
is determined by coefficients of smaller degree.
Computing the coefficients in increasing total degree
therefore defines exactly one formal power series.

We now verify~\eqref{eqNormalizedTransferVector}.
By Lemma~\ref{lemTransferNormalization},
$\epsilon(M_{X_n}(ut,r,s)-M_{X_n}(0,0,0))=0$,
and $\epsilon z_\varnothing=1$
by~\eqref{eqRefinedS6EpsilonDefinition}.
Applying $\epsilon$ to the defining recursion gives
$\epsilon\mathbf b_{X_n}(ut,r,s)=1$.
Replacing $u$ by $ut$ and $s$ by $ts$ also gives
$\epsilon\mathbf b_{X_n}(ut^2,r,ts)=1$.

Expanding $\mathbf b_{X_n}$ in the basis $\{z_A:A\subseteq[n]\}$,
we obtain
\begin{alignat*}{2}
  M_{X_n}(0,0,0)\mathbf b_{X_n}(ut^2,r,ts)
  &=\sum_{A\subseteq[n]}
    b_{X_n,A}(ut^2,r,ts)\,M_{X_n}(0,0,0)z_A\\
  &=z_\varnothing\sum_{A\subseteq[n]}
    b_{X_n,A}(ut^2,r,ts)
  &\quad\text{(by Lemma~\ref{lemTransferNormalization})}&\\
  &=z_\varnothing\epsilon\mathbf b_{X_n}(ut^2,r,ts)
  &\quad\text{(by~\eqref{eqRefinedS6EpsilonDefinition})}&\\
  &=z_\varnothing.
\end{alignat*}
The defining recursion now gives
\begin{align*}
  \mathbf b_{X_n}(ut,r,s)
  &=z_\varnothing+
    \bigl(M_{X_n}(ut,r,s)-M_{X_n}(0,0,0)\bigr)
    \mathbf b_{X_n}(ut^2,r,ts)\\
  &=M_{X_n}(0,0,0)\mathbf b_{X_n}(ut^2,r,ts)
    +\bigl(M_{X_n}(ut,r,s)-M_{X_n}(0,0,0)\bigr)
     \mathbf b_{X_n}(ut^2,r,ts)\\
  &=M_{X_n}(ut,r,s)\mathbf b_{X_n}(ut^2,r,ts).
\end{align*}
Thus both equations in~\eqref{eqNormalizedTransferVector} hold.

Conversely, every solution of~\eqref{eqNormalizedTransferVector}
has $\epsilon\mathbf b_{X_n}=1$.
The basis calculation above therefore gives
$$M_{X_n}(0,0,0)\mathbf b_{X_n}(ut^2,r,ts)=z_\varnothing,$$
so its first equation implies the defining recursion.
Since that recursion determines every coefficient uniquely,
the solution is unique.
\end{proof}

\begin{theorem}\label{thmCombinatorialTransferVector}
For every finite word $\alpha$ in the signed alphabet
$1^-,1^+,\ldots,n^-,n^+$, we have
\[
  F_\alpha(u,r,s)
  =
  \sum_{B\subseteq[n]}
  e_B(\alpha)b_{X_n,B}(ut,r,s).
\]
\end{theorem}

\begin{proof}
For $\alpha=\varnothing$, the sum on the right is
$\epsilon\mathbf b_{X_n}=1$.
For any $\alpha$, we have
\begin{alignat*}{2}
  \sum_{B\subseteq[n]} e_B(\alpha)b_{X_n,B}(ut,r,s) &={} \sum_{A,B\subseteq[n]} e_B(\alpha)M_{X_n}(ut,r,s)_{B,A} b_{X_n,A}(ut^2,r,ts) &\qquad\makebox[88.26pt][r]{\text{(by~\eqref{eqNormalizedTransferVector})}}&\\
   &={} \sum_\beta\mathcal W_{\alpha,\beta}(ut,r,s) \sum_{A\subseteq[n]}e_A(\beta)b_{X_n,A}(ut^2,r,ts) &\qquad\makebox[88.26pt][r]{\text{(by Proposition~\ref{propCombinatorialTransfer})}}\rlap{.}&
\end{alignat*}
Thus these sums satisfy~\eqref{eqNormalizedPathRecursion}
with the prescribed value at $\alpha=\varnothing$.
Proposition~\ref{propNormalizedPathRecursion} gives uniqueness,
so they equal $F_\alpha$.
\end{proof}

\begin{example}
Recall Example~\ref{exDeletingMainDiagonal}, and take $\alpha=(1^-,2^+,6^+)$.
We compare the coefficient of the monomial $x_8$ on the two
sides of Theorem~\ref{thmCombinatorialTransferVector}.
Such a term has exactly one north step in addition to the
three prescribed undecorated roots.

If this step is a child of one of the three roots, it may be
undecorated or carry a rise decoration.
The three choices of parent contribute inversion factors
$q^2,q,1$, respectively, and the row factor is $ut+r$.
Alternatively, the additional step is a root carrying a
valley decoration.  Its three possible positions contribute
$1,q,q^2$, respectively, and its row factor is $s$.
The two possible labels, $8^-$ and $8^+$, have total weight
$(q-1)x_8$.  Therefore
\[
  [x_8]F_\alpha(u,r,s)
  =
  (q-1)(1+q+q^2)(ut+r+s)
  =
  (q^3-1)(ut+r+s).
\]

On the matrix side, Theorem~\ref{thmRiseTransferMatrix} and
\eqref{eqRefinedS6RDefinition} give
\[
  \left.
    M_{X_n}(ut,r,s)z_\varnothing
  \right|_{x_i=0\ (i\neq8)}
  =
  \frac{1-(ut+r)x_8}{1+sx_8}z_\varnothing
  +
  \frac{(ut+r+s)x_8}{1+sx_8}z_{\{8\}}.
\]
The part of $M_{X_n}$ independent of all $x_i$ is
$M_{X_n}(0,0,0)$.
Using~\eqref{eqNormalizedTransferVector} and
$\epsilon[x_8]\mathbf b_{X_n}=0$, we obtain
\begin{alignat*}{2}
  [x_8]\mathbf b_{X_n}(ut,r,s) &= \bigl([x_8]M_{X_n}(ut,r,s)\bigr)z_\varnothing + M_{X_n}(0,0,0)[x_8]\mathbf b_{X_n}(ut^2,r,ts) = (ut+r+s)(z_{\{8\}}-z_\varnothing).
\end{alignat*}
Since $e_\varnothing(\alpha)=1$ and $e_{\{8\}}(\alpha)=q^3$,
\[
  [x_8]\sum_{B\subseteq[n]}
  e_B(\alpha)b_{X_n,B}(ut,r,s)
  =
  (q^3-1)(ut+r+s)
  =
  [x_8]F_\alpha(u,r,s),
\]
as asserted.
\end{example}

We now sum over the undecorated root labels.
Let $w$ be a formal variable, and define a linear functional
on $\mathcal V_n$ by
\begin{equation}\label{eqRefinedS7BetaDefinition}
  \beta_{X_n}[w](z_A)
  :=
  \sum_{\alpha\in\{1^-,1^+,\ldots,n^-,n^+\}^{*}}
  w^{|\alpha|}q^{I(\alpha)}
  \left(\prod_{j=1}^{|\alpha|}\wt(\alpha_j)\right)e_A(\alpha),
\end{equation}
where $\{1^-,1^+,\ldots,n^-,n^+\}^{*}$ denotes the set of
all finite words in the signed alphabet, including the empty word.

\begin{proposition}\label{propRootFunctional}
For every $A\subseteq[n]$, we have
\begin{equation}\label{eqRootFunctional}
  \beta_{X_n}[w](z_A)
  =
  \sum_{J\subseteq[n]\setminus A}
  (-w)^{|J|}(q;q)_{|J|}
  q^{\#\{(j,a)\in J\times A:j<a\}}x_J.
\end{equation}
Moreover, the generating function $R$ defined
in~\eqref{eqRefinedS13RIntroduction} satisfies
\begin{equation}\label{eqRootPairing}
  (U_qR)(X_n;u,v,r,s)
  =
  \beta_{X_n}[vu]\,\mathbf b_{X_n}(ut,r,s).
\end{equation}
\end{proposition}

\begin{proof}
We expand the powers in the $Z_c$ before passing to
$\mathcal Q_n$.
Since $\operatorname{ev}_A$ descends to this quotient
and $S^{n+1}=0$, we obtain
\begin{alignat*}{2}
  \beta_{X_n}[w](z_A)
  &=
  \sum_{m\geq0}w^m
  \sum_{\alpha\in\{1^-,1^+,\ldots,n^-,n^+\}^{m}}
  \operatorname{ev}_A(Z_{\alpha_1}\cdots Z_{\alpha_m})
  &\qquad\text{(by~\eqref{eqRefinedS7BetaDefinition},
                   \eqref{eqRefinedS45Evaluation})}&\\
  &=
  \sum_{m\geq0}w^m
  \operatorname{ev}_A\left(
    \left(\sum_{i=1}^{n}(Z_{i^-}+Z_{i^+})\right)^m
  \right)
  &\qquad\text{(by expansion and linearity)}&\\
  &=
  \sum_{m=0}^{n}w^m\operatorname{ev}_A(S^m)
  &\qquad\text{(by~\eqref{eqRefinedS45QGenerators},
                   \eqref{eqRefinedS45HeadElements})}&\\
  &=
  \sum_{m=0}^{n}\frac{w^m(q;q)_m}{(1-q)^m}
  \sum_{\substack{J\subseteq[n]\\|J|=m}}
  \operatorname{ev}_A(Q_J)
  &\qquad\text{(by~\eqref{eqRefinedS6PowerSumQ})}&\\
  &=
  \sum_{J\subseteq[n]\setminus A}
  (-w)^{|J|}(q;q)_{|J|}
  q^{\#\{(j,a)\in J\times A:j<a\}}x_J
  &\qquad\text{(by~\eqref{lemEvaluationOrderedMonomials})}&\rlap{.}
\end{alignat*}

For the second identity, Proposition~\ref{propUqSignedPaths}
allows us to sum over signed paths.
Group them by their undecorated root word $\alpha$.
A root cannot be a rise, so
$\touch(\pi)=|\alpha|$ by~\eqref{eqRefinedS13PathStatistics}.
We have
\begin{alignat*}{2}
  &\quad \quad (U_qR)(X_n;u,v,r,s)\\
  &=
  \sum_{\alpha\in\{1^-,1^+,\ldots,n^-,n^+\}^{*}}
  (vu)^{|\alpha|}q^{I(\alpha)}
  \left(\prod_{j=1}^{|\alpha|}\wt(\alpha_j)\right)
  F_\alpha(u,r,s)
  &\qquad\text{(by~\eqref{eqRefinedS13RIntroduction},
                   \eqref{eqRefinedS7PathGeneratingFunction})}&\\
  &=
  \sum_{\alpha\in\{1^-,1^+,\ldots,n^-,n^+\}^{*}}
  (vu)^{|\alpha|}q^{I(\alpha)}
  \left(\prod_{j=1}^{|\alpha|}\wt(\alpha_j)\right)
  \sum_{A\subseteq[n]}e_A(\alpha)b_{X_n,A}(ut,r,s)
  &\qquad\text{(by Theorem~\ref{thmCombinatorialTransferVector})}&\\
  &=
  \sum_{A\subseteq[n]}b_{X_n,A}(ut,r,s)
  \sum_{\alpha\in\{1^-,1^+,\ldots,n^-,n^+\}^{*}}
  (vu)^{|\alpha|}q^{I(\alpha)}
  \left(\prod_{j=1}^{|\alpha|}\wt(\alpha_j)\right)e_A(\alpha)
  &&\\
  &=
  \sum_{A\subseteq[n]}
  b_{X_n,A}(ut,r,s)\,\beta_{X_n}[vu](z_A)
  &\qquad\text{(by~\eqref{eqRefinedS7BetaDefinition})}&\\
  &=
  \beta_{X_n}[vu]\,\mathbf b_{X_n}(ut,r,s).
\end{alignat*}
\end{proof}

\subsection{Exchanging adjacent letters}\label{subsecAdjacentColors}
Regard $X_n=(x_1,\ldots,x_n)$ as an ordered alphabet, and let
$s_iX_n=(x_1,\ldots,x_{i-1},x_{i+1},x_i,x_{i+2},\ldots,x_n)$
for $1\leq i<n$.

\begin{definition}\label{defRecheckS7AdjacentExchange}
For $1\leq i<n$, define the linear operator
$\mathscr R_i:\mathcal V_n\to\mathcal V_n$ as follows.
If $i,i+1\in A$, or if $i,i+1\notin A$, set
$\mathscr R_i z_A=z_A$.
Thus $\mathscr R_i$ fixes a basis monomial whenever it
contains both variables $z_i,z_{i+1}$ or contains neither.

For $A\subseteq[n]\setminus\{i,i+1\}$, define
\begin{equation}\label{eqAdjacentExchangeBlock}
\begin{alignedat}{2}
  \mathscr R_i z_{A\cup\{i\}} &= \frac{(1-q)x_{i+1}}{x_i-qx_{i+1}}z_{A\cup\{i\}} + \frac{x_i-x_{i+1}}{x_i-qx_{i+1}}z_{A\cup\{i+1\}},\\
  \mathscr R_i z_{A\cup\{i+1\}} &= \frac{q(x_i-x_{i+1})}{x_i-qx_{i+1}}z_{A\cup\{i\}} + \frac{(1-q)x_i}{x_i-qx_{i+1}}z_{A\cup\{i+1\}}.
\end{alignedat}
\end{equation}
\eqref{eqAdjacentExchangeBlock} specifies the action when the monomial contains
exactly one of $z_i,z_{i+1}$.
\end{definition}

Recall from~\eqref{eqRefinedS6EpsilonDefinition} that $\epsilon z_A=1$ for every $A\subseteq[n]$,
so applying $\epsilon$ sums the coefficients.
In each line of~\eqref{eqAdjacentExchangeBlock}, this sum is $1$:
\[
  \frac{(1-q)x_{i+1}+x_i-x_{i+1}}{x_i-qx_{i+1}}
  =
  \frac{q(x_i-x_{i+1})+(1-q)x_i}{x_i-qx_{i+1}}
  =
  1.
\]
Hence $\epsilon\mathscr R_i z_A=\epsilon z_A$ for every
basis element, and linearity gives
$\epsilon\mathscr R_i=\epsilon$.

\begin{example}
For $n=2$, the operator $\mathscr R_1$ fixes $1$ and $z_1z_2$.
Its action on $z_1,z_2$ gives
\begin{alignat*}{2}
  \mathscr R_1(x_1z_1+x_2z_2) &= \frac{(1-q)x_1x_2+qx_2(x_1-x_2)}{x_1-qx_2}z_1 + \frac{x_1(x_1-x_2)+(1-q)x_1x_2}{x_1-qx_2}z_2 =x_2z_1+x_1z_2.
\end{alignat*}
\end{example}

\begin{proposition}\label{propAdjacentExchange}
For $1\leq i<n$, omitting the common parameters $ut,r,s$
in $M$ and $\mathbf b$, we have
\[
  \mathscr R_iM_{X_n}=M_{s_iX_n}\mathscr R_i,
  \qquad
  \beta_{X_n}[w]=\beta_{s_iX_n}[w]\mathscr R_i,
  \qquad
  \mathbf b_{s_iX_n}=\mathscr R_i\mathbf b_{X_n}.
\]
\end{proposition}

\begin{proof}
We first verify
\begin{equation}\label{eqRefinedS7ExchangeVT}
  \mathscr R_iV_{X_n}(ut,s)
  =V_{s_iX_n}(ut,s)\mathscr R_i,
  \qquad
  \mathscr R_iT_n(z)=T_n(z)\mathscr R_i.
\end{equation}

For the first identity, fix
$A\setminus\{i,i+1\}$ and $B\setminus\{i,i+1\}$
in~\eqref{eqRefinedS6VDefinition}.
The factors at the other positions depend on
$A\cap\{i,i+1\}$ and $B\cap\{i,i+1\}$ only through
$|A\cap\{i,i+1\}|$ and $|B\cap\{i,i+1\}|$.
Both numbers are unchanged by $\mathscr R_i$.
Thus it is enough to do the calculation for $(x_1,x_2)$
and then replace $ut$ by
$utq^{|A\cap(i+1,n]|+|B\cap[1,i)|}$.

By~\eqref{eqAdjacentExchangeBlock},
\[
  \mathscr R_1(qz_1+z_2)=qz_1+z_2,
  \qquad
  \mathscr R_1(x_1z_1+x_2z_2)=x_2z_1+x_1z_2.
\]
Since $x_1-qx_2\neq0$, the four vectors
$1,z_1z_2,qz_1+z_2,x_1z_1+x_2z_2$ form a basis.
Equation~\eqref{eqRefinedS6VDefinition} gives
\begin{align*}
  (1+sx_1)(1+sx_2)V_{(x_1,x_2)}(ut,s)z_\emptyset
  &=(1-utx_1)(1-utx_2)+(ut+s)(x_1z_1+x_2z_2)\\
  &{}+(ut+s)x_1x_2
    \bigl((qut+s)z_1z_2-ut(qz_1+z_2)\bigr).
\end{align*}
For $z_1z_2$, we have
\[
  (1+sx_1)(1+sx_2)V_{(x_1,x_2)}(ut,s)(z_1z_2)
  =1+s(x_1z_1+x_2z_2)+s^2x_1x_2z_1z_2.
\]
For the other two vectors,
\begin{align*}
  (1+sx_1)(1+sx_2)V_{(x_1,x_2)}(ut,s)(qz_1+z_2)
  &=q+1-qut(x_1+x_2)
    +(qut+(q+1)s)(x_1z_1+x_2z_2)\\
  &{}+sx_1x_2
    \bigl((q+1)(qut+s)z_1z_2-qut(qz_1+z_2)\bigr),
\end{align*}
and
\begin{align*}
  (1+sx_1)(1+sx_2)V_{(x_1,x_2)}(ut,s)(x_1z_1+x_2z_2)
  &=x_1+x_2-(q+1)utx_1x_2
    +s(x_1+x_2-qutx_1x_2)(x_1z_1+x_2z_2)\\
  &{}+x_1x_2
    \bigl(ut(qz_1+z_2)+s(qut+s)(x_1+x_2)z_1z_2\bigr).
\end{align*}
In these expressions, the coefficients of
$1,z_1z_2,qz_1+z_2,x_1z_1+x_2z_2$
are symmetric in $x_1,x_2$.
The factor $(1+sx_1)(1+sx_2)$ is also symmetric.
The two displayed actions of $\mathscr R_1$ therefore give
\[
  \mathscr R_1V_{(x_1,x_2)}(ut,s)
  =V_{(x_2,x_1)}(ut,s)\mathscr R_1
\]
on these four vectors, proving the first identity
in~\eqref{eqRefinedS7ExchangeVT}.

For the second identity, let
$B\subseteq I\subseteq[n]\setminus\{i,i+1\}$.
Equation~\eqref{eqRefinedS6TDefinition} gives
\begin{equation}\label{eqVTwoS7TExchangeEntries}
\begin{aligned}
  T_n(z)_{B,I\cup\{i\}}
  &=T_n(z)_{B,I\cup\{i+1\}},\qquad
  T_n(z)_{B\cup\{i\},I\cup\{i\}}
  =T_n(z)_{B\cup\{i+1\},I\cup\{i+1\}},\\
  T_n(z)_{B\cup\{i\},I\cup\{i,i+1\}}
  &=qT_n(z)_{B\cup\{i+1\},I\cup\{i,i+1\}},\qquad
  T_n(z)_{B\cup\{i\},I\cup\{i+1\}}
  =T_n(z)_{B\cup\{i+1\},I\cup\{i\}}=0.
\end{aligned}
\end{equation}
In the third equality, the term on the right removes $i$
and keeps $i+1$, adding the pair $(i,i+1)$ to the count
in~\eqref{eqRefinedS6TDefinition}.
The minus sign in the exponent makes this term
$q^{-1}$ times the term on the left.
The last equality follows directly from the same definition.

The two coefficients in each line
of~\eqref{eqAdjacentExchangeBlock} sum to $1$.
Hence~\eqref{eqVTwoS7TExchangeEntries} gives, for
$k\in\{i,i+1\}$,
\begin{align*}
  T_n(z)\mathscr R_i z_{I\cup\{k\}}
  &=\sum_{B\subseteq I}T_n(z)_{B,I\cup\{k\}}z_B
    +\sum_{B\subseteq I}
      T_n(z)_{B\cup\{k\},I\cup\{k\}}z_B(\mathscr R_i z_k)\\
  &=\mathscr R_i\left(
      \sum_{B\subseteq I}T_n(z)_{B,I\cup\{k\}}z_B
      +\sum_{B\subseteq I}
        T_n(z)_{B\cup\{k\},I\cup\{k\}}z_Bz_k
    \right)\\
  &=\mathscr R_iT_n(z)z_{I\cup\{k\}}.
\end{align*}

For $z_{I\cup\{i,i+1\}}$,
\eqref{eqVTwoS7TExchangeEntries} gives
\begin{align*}
  \quad T_n(z)z_{I\cup\{i,i+1\}}
  &=\sum_{B\subseteq I}
    T_n(z)_{B,I\cup\{i,i+1\}}z_B+\sum_{B\subseteq I}
    T_n(z)_{B\cup\{i+1\},I\cup\{i,i+1\}}
    z_B(qz_i+z_{i+1})\\
  &\quad+\sum_{B\subseteq I}
    T_n(z)_{B\cup\{i,i+1\},I\cup\{i,i+1\}}
    z_Bz_iz_{i+1}.
\end{align*}
Since $\mathscr R_i1=1,$
  $\mathscr R_i(qz_i+z_{i+1})=qz_i+z_{i+1}$ and 
  $\mathscr R_i(z_iz_{i+1})=z_iz_{i+1}$, each term is fixed by $\mathscr R_i$.
Thus
\[
  \mathscr R_iT_n(z)z_{I\cup\{i,i+1\}}
  =T_n(z)z_{I\cup\{i,i+1\}}
  =T_n(z)\mathscr R_i z_{I\cup\{i,i+1\}}.
\]
Finally, every term of $T_n(z)z_I$ contains neither $z_i$
nor $z_{i+1}$, so
\[
  \mathscr R_iT_n(z)z_I
  =T_n(z)z_I
  =T_n(z)\mathscr R_i z_I.
\]
This proves the second identity in~\eqref{eqRefinedS7ExchangeVT}.

Using the factorization in
Lemma~\ref{lemTransferMatrixRelations}\textup{(2)}, equation~\eqref{eqRefinedS7ExchangeVT}
gives
\[
  \mathscr R_i\mathsf J_{X_n}(s)V_{X_n}(ut,0)
  =
  \mathsf J_{s_iX_n}(s)\mathscr R_iV_{X_n}(ut,0).
\]
Since $V_{X_n}(ut,0)$ is invertible over the rational-function
field, we obtain
\[
  \mathscr R_i\mathsf J_{X_n}(s)
  =\mathsf J_{s_iX_n}(s)\mathscr R_i.
\]
Using~\eqref{eqRefinedS6MDefinition} in the rational function
field, we obtain
\begin{alignat*}{2}
  \mathscr R_iM_{X_n} &= \mathscr R_i\mathsf J_{X_n}(s) T_n\left(\frac{r}{ut}\right)V_{X_n}(ut,0) = \mathsf J_{s_iX_n}(s) T_n\left(\frac{r}{ut}\right) V_{s_iX_n}(ut,0)\mathscr R_i = M_{s_iX_n}\mathscr R_i,
\end{alignat*}
which proves the first assertion.

For the second assertion, use~\eqref{eqRootFunctional}.
If $A$ contains both $i,i+1$, the sum does not involve
$x_i,x_{i+1}$.
If $A$ contains neither, exchanging $i,i+1$ in the summation
sets leaves the powers of $q$ unchanged, so the sum is
symmetric in $x_i,x_{i+1}$.
Since $\mathscr R_i z_A=z_A$ in both cases, the assertion follows.

Now fix $A\subseteq[n]\setminus\{i,i+1\}$.
In~\eqref{eqRootFunctional}, combine the terms indexed by
$J$ and $J\cup\{i+1\}$ for $z_{A\cup\{i\}}$,
and those indexed by $J$ and $J\cup\{i\}$ for
$z_{A\cup\{i+1\}}$.
Using
$(q;q)_{|J|+1}=(1-q^{|J|+1})(q;q)_{|J|}$
and
$|J\cap[1,i)|=|J\cap[1,i+1)|$,
we obtain
\begin{gather*}
  \beta_{X_n}[w](z_{A\cup\{i\}})
  =\sum_{J\subseteq[n]\setminus(A\cup\{i,i+1\})}
    (-w)^{|J|}(q;q)_{|J|}
    q^{\#\{(j,a)\in J\times A:j<a\}+|J\cap[1,i)|}x_J
    \left(
      1-w(1-q^{|J|+1})q^{|A\cap(i+1,n]|}x_{i+1}
    \right),
\end{gather*}
and
\begin{gather*}
  \beta_{X_n}[w](z_{A\cup\{i+1\}})
  =\sum_{J\subseteq[n]\setminus(A\cup\{i,i+1\})}
    (-w)^{|J|}(q;q)_{|J|}
    q^{\#\{(j,a)\in J\times A:j<a\}+|J\cap[1,i)|}x_J
    \left(
      1-w(1-q^{|J|+1})q^{1+|A\cap(i+1,n]|}x_i
    \right).
\end{gather*}
The extra $q$ in the second expression comes from the pair
$(i,i+1)$ in~\eqref{eqRootFunctional}.

For each $J$, the two expressions have the same factors
apart from $x_{i+1}$ and $qx_i$.
After exchanging $x_i,x_{i+1}$, these become
$x_i$ and $qx_{i+1}$.
Equation~\eqref{eqAdjacentExchangeBlock} gives
\begin{align*}
  \frac{
    (1-q)x_{i+1}x_i
    +(x_i-x_{i+1})qx_{i+1}
  }{x_i-qx_{i+1}}
  &=\frac{x_{i+1}(x_i-qx_{i+1})}{x_i-qx_{i+1}}
   =x_{i+1},\\
  \frac{
    q(x_i-x_{i+1})x_i
    +(1-q)x_i\,qx_{i+1}
  }{x_i-qx_{i+1}}
  &=\frac{qx_i(x_i-qx_{i+1})}{x_i-qx_{i+1}}
   =qx_i.
\end{align*}
The constant terms in parentheses are also unchanged,
since the two coefficients in each line
of~\eqref{eqAdjacentExchangeBlock} sum to $1$.
Thus, term by term in the sums above,
\[
  \bigl(\beta_{s_iX_n}[w]\mathscr R_i\bigr)
    (z_{A\cup\{k\}})
  =\beta_{X_n}[w](z_{A\cup\{k\}}),
  \qquad k\in\{i,i+1\}.
\]
Together with the cases already proved, this gives
\[
  \beta_{X_n}[w]=\beta_{s_iX_n}[w]\mathscr R_i.
\]

Finally, \eqref{eqNormalizedTransferVector} and the first assertion give
\begin{alignat*}{2}
  \mathscr R_i\mathbf b_{X_n}(ut,r,s) &= \mathscr R_iM_{X_n}(ut,r,s) \mathbf b_{X_n}(ut^2,r,ts) = M_{s_iX_n}(ut,r,s) \mathscr R_i\mathbf b_{X_n}(ut^2,r,ts).
\end{alignat*}
Moreover,
$\epsilon\mathscr R_i\mathbf b_{X_n}
=\epsilon\mathbf b_{X_n}=1$.
The operator $\mathscr R_i$ is independent of $u,r,s$, so it commutes
with the parameter shift in \eqref{eqNormalizedTransferVector}.
Thus $\mathscr R_i\mathbf b_{X_n}$ satisfies the recurrence
and normalization for the alphabet $s_iX_n$.
Proposition~\ref{propNormalizedTransferVector} gives
$\mathbf b_{s_iX_n}=\mathscr R_i\mathbf b_{X_n}$.
\end{proof}

\subsection{The main combinatorial recurrence}\label{subsecMainCombinatorialRecurrence}
Recall the ordered alphabet $X_n=(x_1,\ldots,x_n)$ from
Subsection~\ref{subsecAdjacentColors}, and the notation
$\mathcal V_n$ and $\epsilon$ from \eqref{eqRefinedS6EpsilonDefinition}.
We write $z_x$ for the coordinate corresponding to a letter
$x\in X_n$.

\begin{definition}\label{defRefinedS7Deletion}
Extend scalars to $\mathbb Q(q,t,x_1,\ldots,x_n)$.
For $x_i\in X_n$, define
$\partial_{\{x_i\}}f
  :=
  \left.
    \mathscr R_{n-1}\cdots\mathscr R_i f
  \right|_{z_n=0}$.
The factors act from right to left, and each factor is evaluated
at the current ordered alphabet: before applying $\mathscr R_k$,
this alphabet is $s_{k-1}\cdots s_iX_n$, with $X_n$ used when $k=i$.
Thus each exchange moves $x_i$ one position to the right.
After setting $z_n=0$, we delete the last position and retain
the relative order of the surviving letters.
When $i=n$, the product of exchanges is empty.
For $I=\{x_{i_1},\ldots,x_{i_k}\}\subseteq X_n$, define
$\partial_I
  :=
  \partial_{\{x_{i_k}\}}\cdots\partial_{\{x_{i_1}\}}
  :
  \mathcal V_n\longrightarrow\mathcal V_{n-k}$,
where $\partial_{\{x_{i_j}\}}$ acts on the alphabet remaining
after $x_{i_1},\ldots,x_{i_{j-1}}$ have been deleted.
Set $\partial_\varnothing=\mathrm{id}$ and define
$\epsilon_I^{X_n}(f)
  :=
  \epsilon(\partial_I f)$.
In other words, after deleting the letters of $I$, we set
all remaining coordinates to $1$.
In particular, $\epsilon_\varnothing^{X_n}=\epsilon$.
\end{definition}

\begin{proposition}\label{propDeletionIntertwining}
For $1\leq i<n$ and $I\subseteq X_n$, respectively, we have
\begin{equation}\label{eqBDeletionIntertwining}
 \mathscr R_i\mathcal D_{X_n}=\mathcal D_{s_iX_n}\mathscr R_i,
 \qquad \partial_I\mathcal D_{X_n}
       =\mathcal D_{X_n\setminus I}\partial_I.
\end{equation}
Here $\partial_I$ is defined using any fixed order of the letters of $I$
as in Definition~\ref{defRefinedS7Deletion}, and $X_n\setminus I$
retains the induced order.
\end{proposition}

\begin{proof}
Fix $A\subseteq[n]\setminus\{i,i+1\}$.
Consider first the terms of $\mathcal D$ that add $i$ or $i+1$.
After removing the common factor
$(q-1)q^{|A\cap(i+1,n]|}$,
\eqref{eqAdjacentExchangeBlock} and~\eqref{eqRefinedS6ABCD}
give the following contributions on both sides of the first equality:
\begin{align*}
  z_A
  \longmapsto
  x_{i+1}z_{A\cup\{i\}}+x_i z_{A\cup\{i+1\}}, \quad
  z_{A\cup\{i\}}
  \longmapsto
  x_{i+1}z_{A\cup\{i,i+1\}},\quad
  z_{A\cup\{i+1\}}
  \longmapsto
  qx_i z_{A\cup\{i,i+1\}}, \quad
  z_{A\cup\{i,i+1\}}
  \longmapsto 0.
\end{align*}
For a term adding $k\notin\{i,i+1\}$, the coefficients
in~\eqref{eqAdjacentExchangeBlock} are unchanged when $k$
is adjoined to $A$, and
\[
  |(A\cup\{i\})\cap(k,n]|
  =
  |(A\cup\{i+1\})\cap(k,n]|.
\]
Thus these terms also agree.
This proves the first equality
in~\eqref{eqBDeletionIntertwining}.

For deletion of the last letter, \eqref{eqRefinedS6ABCD}
gives, for every $A\subseteq[n-1]$,
\begin{align*}
  \mathcal D_{X_n}z_A
  =
  \mathcal D_{X_n\setminus\{x_n\}}z_A
  +(q-1)x_nz_{A\cup\{n\}},\quad
  \mathcal D_{X_n}z_{A\cup\{n\}}
  =
  qz_n\mathcal D_{X_n\setminus\{x_n\}}z_A.
\end{align*}
Setting $z_n=0$ in these two identities and using linearity yields
\[
  \left.\mathcal D_{X_n}f\right|_{z_n=0}
  =
  \mathcal D_{X_n\setminus\{x_n\}}
  \left(\left.f\right|_{z_n=0}\right).
\]

Now use Definition~\ref{defRefinedS7Deletion}, with each
exchange evaluated at the current ordered alphabet.
The first equality already proved and the preceding restriction
identity give
\begin{align*}
  \partial_{\{x_i\}}\mathcal D_{X_n}f
  &=
  \left.
    \mathscr R_{n-1}\cdots\mathscr R_i
    \mathcal D_{X_n}f
  \right|_{z_n=0}\\
  &=
  \left.
    \mathcal D_{s_{n-1}\cdots s_iX_n}
    \mathscr R_{n-1}\cdots\mathscr R_i f
  \right|_{z_n=0}\\
  &=
  \mathcal D_{X_n\setminus\{x_i\}}
  \left(
    \left.\mathscr R_{n-1}\cdots\mathscr R_i f\right|_{z_n=0}
  \right)\\
  &=
  \mathcal D_{X_n\setminus\{x_i\}}\partial_{\{x_i\}}f.
\end{align*}
Here $x_i$ is the last letter of
$s_{n-1}\cdots s_iX_n$, and deleting it preserves the
relative order of the remaining letters.
Composing this identity in the chosen deletion order gives
\[
  \partial_I\mathcal D_{X_n}
  =
  \mathcal D_{X_n\setminus I}\partial_I,
\]
as required.
\end{proof}

The proof of the following lemma is given in Appendix~\ref{app:deletion-cyclic}.

\begin{lemma}\label{lemMarkedDeletionIdentities}
The map $\partial_I$ is independent of the order of deletion, and
$\epsilon_I^{X_n}=\epsilon_I^{s_iX_n}\mathscr R_i$,
where $I$ denotes the same set of letters in both alphabets.
Moreover:
\begin{enumerate}[(1)]
  \item For $I\subseteq X_n$,
  \begin{equation}\label{eqRefinedS7DeletionD}
    \epsilon_I^{X_n}\mathcal D
    =
    (q-1)\sum_{h\in X_n\setminus I}
    h\eta_h(X_n\setminus I)\epsilon_{I\cup\{h\}}^{X_n}.
  \end{equation}

  \item Let $X_n=(x,X_{n-1})$ and $I\subseteq X_{n-1}$,
  with $\mathcal A$ acting on the coordinates of $X_{n-1}$.
  For $f\in\mathcal V_n$,
  \begin{equation}\label{eqRefinedS7DeletionA}
\begin{alignedat}{2}
  \epsilon_I^{X_{n-1}}\mathcal A(f|_{z_x=0}) &= \prod_{h\in X_{n-1}\setminus I}\frac{x-qh}{x-h}\, \epsilon_{I\cup\{x\}}^{X_n}f + (1-q)\sum_{h\in X_{n-1}\setminus I} \frac{h\eta_h(X_{n-1}\setminus I)}{h-x}\, \epsilon_{I\cup\{h\}}^{X_n}f.
\end{alignedat}
\end{equation}
\end{enumerate}
\end{lemma}

Recall the functional $\beta_{X_n}[w]$ from
Proposition~\ref{propRootFunctional}.
We next express it in terms of the functionals
$\epsilon_I^{X_n}$.

\begin{proposition}\label{propRootMarkedExpansion}
We have
\begin{equation}\label{eqRootMarkedExpansion}
  \beta_{X_n}[w]
  =
  \sum_{I\subseteq X_n}
  (-w)^{|I|}(q;q)_{|I|}
  \left(\prod_{x\in I}x\right)
  \eta_I(X_n)\epsilon_I^{X_n}.
\end{equation}
\end{proposition}

\begin{proof}
By~\eqref{eqRefinedS6ABCD}, each application of $\mathcal D$
adds an index not already in the subset.
For a fixed set of added indices, a pair contributes $q$
exactly when its larger index is added first.
Inserting the largest index into an order of the other
$j-1$ indices gives the weights $1,q,\ldots,q^{j-1}$.
Thus summing over all orders gives
$\prod_{h=1}^{j}(1+q+\cdots+q^{h-1})$.
The comparisons with the original indices of $A$ are
independent of the order, so
\[
  \epsilon\mathcal D^j(z_A)
  =
  (-1)^j(q;q)_j
  \sum_{\substack{J\subseteq[n]\setminus A\\|J|=j}}
  q^{\#\{(i,a)\in J\times A:i<a\}}x_J.
\]
Comparing with~\eqref{eqRootFunctional} gives
\begin{equation}\label{eqRefinedS7BetaD}
  \beta_{X_n}[w]
  =\sum_{j=0}^{n}w^j\epsilon\mathcal D^j.
\end{equation}

We prove by induction on $j$ that, for every finite alphabet
$X_n$ and $0\leq j\leq n$,
\begin{equation}\label{eqRefinedS7IteratedDeletion}
  \epsilon\mathcal D^j
  =
  (-1)^j(q;q)_j
  \sum_{\substack{I\subseteq X_n\\|I|=j}}
  \left(\prod_{x\in I}x\right)
  \eta_I(X_n)\epsilon_I^{X_n}.
\end{equation}
The case $j=0$ follows from $\eta_\varnothing(X_n)=1$ and
$\epsilon_\varnothing^{X_n}=\epsilon$.

For $j\geq1$, apply~\eqref{eqRefinedS7DeletionD}
with $I=\varnothing$, and then use the induction hypothesis
on $X_n\setminus\{x\}$:
\begin{alignat*}{2}
  \epsilon\mathcal D^j
  &=(q-1)\sum_{x\in X_n}
    x\eta_x(X_n)\epsilon_{\{x\}}^{X_n}\mathcal D^{j-1}
  &\quad\text{(by~\eqref{eqRefinedS7DeletionD})}&\\
  &=(q-1)\sum_{x\in X_n}
    x\eta_x(X_n)\epsilon
    \partial_{\{x\}}\mathcal D_{X_n}^{j-1}
  &\quad\text{(by Definition~\ref{defRefinedS7Deletion})}&\\
  &=(q-1)\sum_{x\in X_n}
    x\eta_x(X_n)\epsilon
    \mathcal D_{X_n\setminus\{x\}}^{j-1}\partial_{\{x\}}
  &\quad\text{(by~\eqref{eqBDeletionIntertwining})}&\\
  &=(-1)^{j-1}(q-1)(q;q)_{j-1}
    \sum_{x\in X_n}
    \sum_{\substack{J\subseteq X_n\setminus\{x\}\\|J|=j-1}}
    x\left(\prod_{h\in J}h\right)
    \eta_x(X_n)\eta_J(X_n\setminus\{x\})
    \epsilon_J^{X_n\setminus\{x\}}\partial_{\{x\}}
  &\quad\text{(by~\eqref{eqRefinedS7IteratedDeletion})}&\rlap{.}
\end{alignat*}

Here $\partial_J$ acts on $X_n\setminus\{x\}$.
The composition $\partial_J\partial_{\{x\}}$ therefore
deletes $x$ first and then the letters of $J$.
By Lemma~\ref{lemMarkedDeletionIdentities}, we may compute
$\partial_{J\cup\{x\}}$ by deleting $x$ first and then the
letters of $J$.
Definition~\ref{defRefinedS7Deletion} therefore gives
\[
  \epsilon_J^{X_n\setminus\{x\}}\partial_{\{x\}}
  =
  \epsilon\partial_J\partial_{\{x\}}
  =
  \epsilon\partial_{J\cup\{x\}}
  =
  \epsilon_{J\cup\{x\}}^{X_n},
\]
where $\partial_J$ acts on $X_n\setminus\{x\}$.

Substitute this identity into the preceding sum.
For each $I=J\cup\{x\}$, the corresponding pairs are exactly
$(x,I\setminus\{x\})$ with $x\in I$.
Hence
\begin{alignat*}{2}
  \epsilon\mathcal D^j
  &=(-1)^{j-1}(q-1)(q;q)_{j-1}
    \sum_{\substack{I\subseteq X_n\\|I|=j}}
    \left(\prod_{x\in I}x\right)
    \left(
      \sum_{x\in I}
      \eta_x(X_n)
      \eta_{I\setminus\{x\}}(X_n\setminus\{x\})
    \right)\epsilon_I^{X_n}
  &\quad\text{(by $I=J\cup\{x\}$)}&\\
  &=(-1)^{j-1}(q-1)(q;q)_{j-1}
    \sum_{\substack{I\subseteq X_n\\|I|=j}}
    \left(\prod_{x\in I}x\right)\eta_I(X_n)
    \left(\sum_{x\in I}\eta_x(I)\right)\epsilon_I^{X_n}
  &\quad\text{(by~\eqref{eqEtaRemoval})}&\\
  &=(-1)^j(q;q)_j
    \sum_{\substack{I\subseteq X_n\\|I|=j}}
    \left(\prod_{x\in I}x\right)
    \eta_I(X_n)\epsilon_I^{X_n}
  &\quad\text{(by~\eqref{eqEtaSubsetSum})}&\rlap{.}
\end{alignat*}
This proves~\eqref{eqRefinedS7IteratedDeletion}.
Substituting it into~\eqref{eqRefinedS7BetaD}
gives~\eqref{eqRootMarkedExpansion}.
\end{proof}

Recall from Theorem~\ref{thmRiseTransferMatrix} that
$M_{X_n}(ut,r,s)=\mathsf J_{X_n}(s)R_{X_n}(ut,r)$.
By Proposition~\ref{propNormalizedTransferVector},
$\mathbf b_{X_n}$ is the unique formal vector satisfying
\[
  \mathbf b_{X_n}(ut,r,s)
  =
  M_{X_n}(ut,r,s)\mathbf b_{X_n}(ut^2,r,ts),
  \qquad
  \epsilon\mathbf b_{X_n}(ut,r,s)=1.
\]
For $I\subseteq X_n$, define
\begin{equation}\label{eqRefinedS7RDefinition}
  \mathcal R_I(X_n;ut,r,s)
  :=
  \epsilon_I^{X_n}\mathbf b_{X_n}(ut,r,s).
\end{equation}
In particular, $\mathcal R_\varnothing=1$.

\begin{theorem}\label{thmCombinatorialMarkedExpansion}
For a finite alphabet $X_n$ and $0\leq j\leq n$,
\begin{equation}\label{eqCombinatorialMarkedExpansion}
  [v^j](U_qR)(X_n;u,v,r,s)
  =
  (-u)^j(q;q)_j
  \sum_{\substack{I\subseteq X_n\\|I|=j}}
  \left(\prod_{x\in I}x\right)
  \eta_I(X_n)\mathcal R_I(X_n;ut,r,s).
\end{equation}
For $j>n$, the coefficient is zero.
\end{theorem}

\begin{proof}
We have
\begin{alignat*}{2}
  (U_qR)(X_n;u,v,r,s) &= \beta_{X_n}[vu]\mathbf b_{X_n}(ut,r,s) &\qquad\makebox[44.58pt][r]{\text{(by~\eqref{eqRootPairing})}}&\\
   &= \sum_{I\subseteq X_n} (-uv)^{|I|}(q;q)_{|I|} \left(\prod_{x\in I}x\right)\eta_I(X_n) \epsilon_I^{X_n}\mathbf b_{X_n}(ut,r,s) &\qquad\makebox[44.58pt][r]{\text{(by~\eqref{eqRootMarkedExpansion})}}&\\
   &= \sum_{I\subseteq X_n} (-uv)^{|I|}(q;q)_{|I|} \left(\prod_{x\in I}x\right)\eta_I(X_n) \mathcal R_I(X_n;ut,r,s) &\qquad\makebox[44.58pt][r]{\text{(by~\eqref{eqRefinedS7RDefinition})}}\rlap{.}&
\end{alignat*}
Taking the coefficient of $v^j$ proves the assertion.
\end{proof}

The proof of the following lemma is given in Appendix~\ref{app:deletion-cyclic}.

\begin{lemma}\label{lemCyclicTransferVector}
Let $X_n=(X_{n-1},x)$, and let $\mathcal A,\mathcal D$
act on the coordinates of $X_{n-1}$. Then
\begin{equation}\label{eqCyclicTransferVector}
\begin{gathered}
  \mathbf b_{X_n}|_{z_x=1}
  =
  \mathbf b_{(tx,X_{n-1})}|_{z_{tx}=1},\\
  \begin{alignedat}{2}
  (1+rx)(1+sx)\mathbf b_{X_n}|_{z_x=0} &= (1-rsx\mathcal D) 
  \mathbf b_{(tx,X_{n-1})}|_{z_{tx}=1} + tx(rsx-u)\mathcal A \mathbf b_{(tx,X_{n-1})}|_{z_{tx}=0},
\end{alignedat}
\end{gathered}
\end{equation}
where the common parameters $ut,r,s$ are omitted.
\end{lemma}

\begin{theorem}\label{thmMainCombinatorialRecursion}
Fix a finite alphabet $X_n$, a nonempty subset $I\subseteq X_n$,
and $x\in I$.  Then
\begin{equation}\label{eqMainCombinatorialRecursion}
\begin{aligned}
  &(1+rx)(1+sx)\mathcal R_I(X_n;ut,r,s)
  =
  \mathcal R_{I\setminus\{x\}}(X_n^{x\to tx};ut,r,s)\\
  &\quad+
  tx(rsx-u)
  \prod_{h\in X_n\setminus I}\frac{tx-qh}{tx-h}\,
  \mathcal R_{(I\setminus\{x\})\cup\{tx\}}
  (X_n^{x\to tx};ut,r,s)\\
  &\quad+
  (1-q)x\sum_{h\in X_n\setminus I}
  \frac{h\eta_h(X_n\setminus I)}{h-tx}(rsh-ut)\,
  \mathcal R_{(I\setminus\{x\})\cup\{h\}}
  (X_n^{x\to tx};ut,r,s).
\end{aligned}
\end{equation}
Together with $\mathcal R_\varnothing=1$, these identities
uniquely determine the family $\mathcal R_I$ for all finite
alphabets as formal power series in $u,r,s$ with coefficients
in $\mathbb Q(q,t,x_1,\ldots,x_n)$.
\end{theorem}

\begin{proof}
For an adjacent exchange, equation~\eqref{eqRefinedS7RDefinition}, Proposition~\ref{propAdjacentExchange}
and Lemma~\ref{lemMarkedDeletionIdentities} give
\begin{alignat*}{2}
  \mathcal R_I(s_iX_n;ut,r,s) &= \epsilon_I^{s_iX_n}\mathbf b_{s_iX_n}(ut,r,s) = \epsilon_I^{s_iX_n}\mathscr R_i(X_n) \mathbf b_{X_n}(ut,r,s) = \epsilon_I^{X_n}\mathbf b_{X_n}(ut,r,s) = \mathcal R_I(X_n;ut,r,s).
\end{alignat*}
Here $I$ denotes the same set of letters in both alphabets.
Since any reordering is a sequence of adjacent exchanges,
we may place $x$ last and write $X_n=(X_{n-1},x)$.

By Lemma~\ref{lemMarkedDeletionIdentities}, we may delete $x$
before the other letters of $I$.
As $x$ already occupies the last position, this deletion
requires no exchanges: $\partial_{\{x\}}$ simply sets $z_x=0$.
Thus
\begin{align*}
  &(1+rx)(1+sx)\mathcal R_I(X_n;ut,r,s)\\
  =\quad &
  (1+rx)(1+sx)\epsilon_{I\setminus\{x\}}^{X_{n-1}}
  \left(\partial_{\{x\}}\mathbf b_{X_n}(ut,r,s)\right)\\
  =\quad &
  (1+rx)(1+sx)\epsilon_{I\setminus\{x\}}^{X_{n-1}}
  \left(\mathbf b_{X_n}(ut,r,s)|_{z_x=0}\right)\\
  =\quad &
  \epsilon_{I\setminus\{x\}}^{X_{n-1}}
  (1-rsx\mathcal D)
  \left(\mathbf b_{(tx,X_{n-1})}(ut,r,s)|_{z_{tx}=1}\right)+
  tx(rsx-u)\epsilon_{I\setminus\{x\}}^{X_{n-1}}
  \mathcal A
  \left(\mathbf b_{(tx,X_{n-1})}(ut,r,s)|_{z_{tx}=0}\right),
\end{align*}
where \eqref{eqCyclicTransferVector} gives the last equality.

Apply Lemma~\ref{lemMarkedDeletionIdentities}\textup{(1)}
on $X_{n-1}$ with the subset $I\setminus\{x\}$, and
equation~\eqref{eqRefinedS7DeletionA} on $(tx,X_{n-1})$ with the same subset.
In both cases, the letters in the products and sums form
$X_{n-1}\setminus(I\setminus\{x\})=X_n\setminus I$.
Since deletion inside $X_{n-1}$ commutes with evaluation
at $z_{tx}=1$, we obtain
\begin{align*}
  &(1+rx)(1+sx)\mathcal R_I(X_n;ut,r,s)\\
  =\quad&
  \epsilon_{I\setminus\{x\}}^{X_{n-1}}
  \left(\mathbf b_{(tx,X_{n-1})}(ut,r,s)|_{z_{tx}=1}\right)-
  rsx\epsilon_{I\setminus\{x\}}^{X_{n-1}}\mathcal D
  \left(\mathbf b_{(tx,X_{n-1})}(ut,r,s)|_{z_{tx}=1}\right)\\
  &\quad+
  tx(rsx-u)\epsilon_{I\setminus\{x\}}^{X_{n-1}}\mathcal A
  \left(\mathbf b_{(tx,X_{n-1})}(ut,r,s)|_{z_{tx}=0}\right)\\
  =\quad&
  \mathcal R_{I\setminus\{x\}}(X_n^{x\to tx};ut,r,s) -
  (q-1)rsx\sum_{h\in X_n\setminus I}
  h\eta_h(X_n\setminus I)\,
  \mathcal R_{(I\setminus\{x\})\cup\{h\}}
  (X_n^{x\to tx};ut,r,s)\\
  &\quad+
  tx(rsx-u)
  \prod_{h\in X_n\setminus I}\frac{tx-qh}{tx-h}\,
  \mathcal R_{(I\setminus\{x\})\cup\{tx\}}
  (X_n^{x\to tx};ut,r,s)\\
  &\quad+
  (1-q)tx(rsx-u)\sum_{h\in X_n\setminus I}
  \frac{h\eta_h(X_n\setminus I)}{h-tx}\,
  \mathcal R_{(I\setminus\{x\})\cup\{h\}}
  (X_n^{x\to tx};ut,r,s).
\end{align*}
Here the term containing $\mathcal D$ gives the first sum,
while the term containing $\mathcal A$ gives the product
and the second sum.
The two sums involve the same functions, and their scalar
coefficients combine as
\[
  -(q-1)rsx+\frac{(1-q)tx(rsx-u)}{h-tx}
  =
  (1-q)x\left(rs+\frac{t(rsx-u)}{h-tx}\right) = (1-q)x \frac{rsh-ut}{h-tx}.
\]
Consequently,
\begin{align*}
  &(1+rx)(1+sx)\mathcal R_I(X_n;ut,r,s)
  =
  \mathcal R_{I\setminus\{x\}}(X_n^{x\to tx};ut,r,s)\\
  &\quad+
  tx(rsx-u)
  \prod_{h\in X_n\setminus I}\frac{tx-qh}{tx-h}\,
  \mathcal R_{(I\setminus\{x\})\cup\{tx\}}
  (X_n^{x\to tx};ut,r,s)\\
  &\quad+
  (1-q)x\sum_{h\in X_n\setminus I}
  \frac{h\eta_h(X_n\setminus I)}{h-tx}(rsh-ut)\,
  \mathcal R_{(I\setminus\{x\})\cup\{h\}}
  (X_n^{x\to tx};ut,r,s),
\end{align*}
which is just \eqref{eqMainCombinatorialRecursion}.

For uniqueness, suppose that all coefficients of total degree
less than $m$ are known, for every finite alphabet and subset.
In degree $m$, every term on the right of
\eqref{eqMainCombinatorialRecursion}, except
$\mathcal R_{I\setminus\{x\}}(X_n^{x\to tx};ut,r,s)$,
uses only these known coefficients.
The nonconstant terms of $(1+rx)(1+sx)$ on the left
also multiply coefficients of lower degree.
Thus the coefficients of degree $m$ in $\mathcal R_I$
are determined by those for subsets with one fewer letter.
Induction on $|I|$, starting from $\mathcal R_\varnothing=1$,
determines all coefficients of degree $m$.
Induction on $m$, beginning with $m=0$, proves uniqueness.
\end{proof}

\subsection{Proof of Theorem~\ref{thmIntroMain}}\label{subsecMainProof}
We now prove our main theorem.

\begin{theorem}\label{thmTouchingTheta}
For integers $d\geq j\geq1$ and $k,\ell\geq0$, we have
\[
  \Theta_{e_\ell}\Theta_{e_k}\nabla E_{d,j}
  =
  \sum_{\substack{
    \pi\in\PF_{d+k+\ell}^{\ast k,\bullet\ell}\\
    \touch(\pi)=j
  }}
  q^{\dinv(\pi)}t^{\area(\pi)}x^\pi.
\]
\end{theorem}

\begin{proof}
By Theorems~\ref{thmAlgMainRecursion}
and~\ref{thmMainCombinatorialRecursion},
$\mathcal L_I$ and $\mathcal R_I$ satisfy the same recurrence
and the initial condition
$\mathcal L_\varnothing=\mathcal R_\varnothing=1$.
Uniqueness in the latter theorem gives
\[
  \mathcal R_I(X_n;ut,r,s)
  =
  \mathcal L_I(X_n;ut,r,s)
  \qquad(n\geq0,\ I\subseteq X_n).
\]
Substituting into~\eqref{eqCombinatorialMarkedExpansion}
and using Theorem~\ref{thmAlgMainRecursion}, we obtain
\begin{align*}
  (U_qR)(X_n;u,v,r,s)
  &=
  1+\sum_{j=1}^{n}(-uv)^j(q;q)_j
  \sum_{\substack{I\subseteq X_n\\|I|=j}}
  \left(\prod_{x\in I}x\right)
  \eta_I(X_n)\mathcal L_I(X_n;ut,r,s)\\
  &=1+L[(q-1)X_n;u,v,r,s].
\end{align*}

For fixed $d,j,k,\ell$, the coefficients of $u^dv^jr^ks^\ell$
on both sides are homogeneous quasisymmetric functions
of degree $d+k+\ell$.
Their specializations agree for every $X_n$, so the equality
also holds for $X$.

Since $L$ has symmetric coefficients,
Proposition~\ref{thmUqProperties}\textup{(1)} gives
\[
  (U_qR)(X;u,v,r,s)
  =
  1+L[(q-1)X;u,v,r,s]
  =
  U_q\bigl(1+L[X;u,v,r,s]\bigr).
\]
Applying Proposition~\ref{thmUqProperties}\textup{(2)}
coefficientwise yields
$R(X;u,v,r,s)=1+L[X;u,v,r,s]$.
Comparing the coefficients of $u^dv^jr^ks^\ell$
in~\eqref{eqThetaalg}
and~\eqref{eqRefinedS13RIntroduction}
proves the theorem.
\end{proof}

\begin{proof}[Proof of Corollary~\ref{corThetaDeltaConsequences}]
\textup{(1)}
Sum Theorem~\ref{thmTouchingTheta} over $j=1,\ldots,d$
and use \eqref{eqRefinedS13ESum}.

\textup{(2)}
Set $k=0$ in part~\textup{(1)}.
By Theorem~\ref{thmThetaDeltaIdentity},
\[
  \Delta'_{e_{d-1}}e_{d+\ell}
  =
  \Theta_{e_\ell}\nabla e_d
  =
  \sum_{\pi\in\PF_{d+\ell}^{\ast0,\bullet\ell}}
  q^{\dinv(\pi)}t^{\area(\pi)}x^\pi,
\]
as required.

\textup{(3)}
Setting $k=0$ in Theorem~\ref{thmTouchingTheta}
gives the hypothesis of~\cite[Theorem~5.2]{IW22}.
Consequently, for every $m\geq0$,
\[
  \Delta_{h_m}\Theta_{e_\ell}\nabla E_{d,j}
  =
  \sum_{\substack{\pi\\\touch(\pi)=j}}
  q^{\dinv(\pi)}t^{\area(\pi)}x^\pi,
\]
where $\pi$ ranges over the partial parking functions
in Corollary~\ref{corThetaDeltaConsequences}\textup{(3)}.
Here $\touch(\pi)$ counts the roots with positive labels
that carry no decoration, as in~\cite[Definition~3.8]{IW22}.
Summing over $j$ and using
Theorem~\ref{thmThetaDeltaIdentity}, we obtain
\begin{alignat*}{2}
  \Delta_{h_m}\Delta'_{e_{d-1}}e_{d+\ell} = \Delta_{h_m}\Theta_{e_\ell}\nabla e_d = \sum_{j=1}^{d} \Delta_{h_m}\Theta_{e_\ell}\nabla E_{d,j}
   = \sum_{j=1}^{d} \sum_{\substack{\pi\\\touch(\pi)=j}} q^{\dinv(\pi)}t^{\area(\pi)}x^\pi = \sum_\pi q^{\dinv(\pi)}t^{\area(\pi)}x^\pi.
\end{alignat*}
This proves part~\textup{(3)}.
\end{proof}

\appendix
\numberwithin{equation}{section}

\section{Proofs of Lemmas~\ref{lemHomogeneousFactor}--\ref{lemPhiOperatorCalculation}}
\label{appRiseCalculations}
We first establish the auxiliary Lemmas~\ref{lemPhiRecurrence},
\ref{lemPhiZeroRelation}, and \ref{lemPhiWeightedRelation}
needed for these proofs.

\begin{lemma}\label{lemPhiRecurrence}
For $E\neq\varnothing$, we have
\[
  \Phi_c(A,E;r)
  =
  \sum_{A=c'A'}\Phi_{c'}(A',E;r)
  +\sum_{d \in \operatorname{supp}(E)}\Gamma_{c,d}(r)\Phi_d(A,E-d;r),
\]
where
\begin{equation}\label{eqRefinedS6GammaDefinition}
  \Gamma_{c,d}(r)
  :=q^{-\lambda_d(A)}
  \Bigl[
    \mathbf{1}_{\{c\prec d\}}
      \bigl(1+rq^{-|E|}\sigma_{A+E}(S-T_c)\bigr)+\bigl(
      ut(1-\mathbf{1}_{\{c\prec d\}})
      +rq^{1-|E|+\lambda_d(A)+v_d}
    \bigr)\sigma_{A+E}T_c
  \Bigr].
\end{equation}
\end{lemma}

\begin{proof}
We prove
\begin{equation}\label{eqRisePhiGTarget}
\begin{aligned}
  &G_{r/q^{|E|}}(\sigma_c\sigma_{A+E}S)\Phi_c(A,E;r)\\
  &=
  G_{r/q^{|E|}}(\sigma_c\sigma_{A+E}S)
  \left(
    \sum_{A=c'A'}\Phi_{c'}(A',E;r)
    +\sum_{d\in\operatorname{supp}(E)}
      \Gamma_{c,d}(r)\Phi_d(A,E-d;r)
  \right).
\end{aligned}
\end{equation}
The common left factor is invertible by~\eqref{eqdefGr},
so cancelling it gives the required recurrence.
We calculate the left-hand side using~\eqref{eqPhiNormalization}
and the sum defining $K_{cA,E}$ in~\eqref{eqRefinedS6KDefinition}.

We split the sum in~\eqref{eqRefinedS6KDefinition} according
to the second root of $W=(c_1,\ldots,c_L)$, where $c_1=c$.
If this root is undecorated, write $A=c'A'$; removing $c$
and its heads leaves data for $K_{c'A',E}$.
If it is decorated with label $d$, the remaining data are
those for $K_{dA,E-d}$, with $d$ made undecorated.
Both operations are reversible.
Equal copies of $d$ are not distinguished, so each
$d\in\operatorname{supp}(E)$ is counted once.

For $h\geq1$ heads assigned to $c$, only the first label is
restricted by $c$.
Their comparisons with the later roots are given by
$\sigma_{A+E}$.
Thus
\begin{alignat}{2}
  \sum_{h\geq1}(ut)^{h-1}(ut+rq^{\nu_1})
    \sigma_{A+E}(T_cS^{h-1})
  &=(ut+rq^{\nu_1})\sigma_{A+E}(T_c)
    \sum_{h\geq1}\bigl(ut\,\sigma_{A+E}(S)\bigr)^{h-1}
  &\quad\text{(by~\eqref{eqRefinedS45HeadElements},
                   \eqref{eqRefinedS6SigmaWords})}&\notag\\
  &=(ut+rq^{\nu_1})\sigma_{A+E}T_c
    (1-ut\sigma_{A+E}S)^{-1}.
    \label{eqRiseFirstHeadSum}
\end{alignat}
By~\eqref{eqRefinedS45RootStatistics} and
\eqref{eqRefinedS45RiseCorrection}, making the first remaining
root undecorated does not change the comparisons or rise
corrections of the remaining heads.

Suppose first that the second root is $c'$.
The empty choice of heads is allowed, and \eqref{eqRefinedS45RiseCorrection} gives
\begin{equation*}
  \nu_1 = -\sum_{k=3}^{L}\delta_k = -|E|
\end{equation*}
The comparisons between $c$ and the letters of $A$ are
cancelled by $I(cA)$ in~\eqref{eqRefinedS6KDefinition}.
For each occurrence $e\in E$, prepending $c$ contributes
$q^{\mathbf1_{\{c\prec e\}}}$.
Using~\eqref{eqPhiNormalization}, we obtain
\begin{align*}
  q^{\mathbf1_{\{c\prec e\}}}
  \frac{\wt(e)q^{\lambda_e(A)-1}}
       {\wt(e)q^{\lambda_e(cA)-1}}
  =q^{\mathbf1_{\{c\prec e\}}+\lambda_e(A)-\lambda_e(cA)}
    =1
\end{align*}

Now suppose that the second root is a decorated root $d$.
The empty choice of heads is allowed exactly when $c\prec d$,
by Proposition~\ref{propMainDiagonalDeletion}\textup{(2)}.
Since $\delta_2=1$, we have
\begin{alignat*}{2}
  \nu_1
  &=\sum_{\substack{3\leq k\leq L\\\delta_k=0}}
      \mathbf1_{\{c_k\prec d\}}
    +\sum_{\substack{3\leq k\leq L\\\delta_k=1}}
      \mathbf1_{\{d\prec c_k\}}
    -\sum_{k=3}^{L}\delta_k
  &\quad\text{(by~\eqref{eqRefinedS45RiseCorrection})}&\\
  &=\lambda_d(A)+v_d-(|E|-1)
  &\quad\text{(by~\eqref{eqRefinedS45InsertionCounts},
                   \eqref{eqRefinedS6RemainingCounts})}&\\
  &=1-|E|+\lambda_d(A)+v_d.
\end{alignat*}
The letter $d$ contributes
$\wt(d)q^{\mathbf1_{\{c\prec d\}}-1}$.
For each occurrence $e\in E-d$, prepending $c$ adds the
comparison from $c$ to $e$, while decorating $d$ removes the
comparison from $d$ to $e$.
The comparisons between undecorated roots are cancelled by
$I(cA)$ and $I(dA)$.
Consequently, substitution into~\eqref{eqPhiNormalization} gives
\begin{alignat*}{2}
  \frac{\wt(d)q^{\mathbf1_{\{c\prec d\}}-1}}
       {\wt(d)q^{\lambda_d(cA)-1}}
  &=q^{\mathbf1_{\{c\prec d\}}-\lambda_d(cA)}
    =q^{-\lambda_d(A)}
  &\quad\text{(by~\eqref{eqRefinedS45InsertionCounts})}&\rlap{,}\\
  q^{\mathbf1_{\{c\prec e\}}-\mathbf1_{\{d\prec e\}}}
  \frac{\wt(e)q^{\lambda_e(dA)-1}}
       {\wt(e)q^{\lambda_e(cA)-1}}
  &=q^{\mathbf1_{\{c\prec e\}}-\mathbf1_{\{d\prec e\}}
        +\lambda_e(dA)-\lambda_e(cA)}
    =1
  &\quad\text{(by~\eqref{eqRefinedS45InsertionCounts})}&\rlap{.}
\end{alignat*}

Combining these calculations with~\eqref{eqRiseFirstHeadSum}
and using
$\sigma_{c'}\sigma_{A'+E}
 =\sigma_d\sigma_{A+(E-d)}
 =\sigma_{A+E}$, we obtain
\begin{align}
  &G_{r/q^{|E|}}(\sigma_c\sigma_{A+E}S)\Phi_c(A,E;r)
  =\sum_{A=c'A'}
    \Bigl(
      1+(ut+rq^{-|E|})
        \sigma_{A+E}T_c(1-ut\sigma_{A+E}S)^{-1}
    \Bigr)
    G_{r/q^{|E|}}(\sigma_{A+E}S)\Phi_{c'}(A',E;r)
    \notag\\
  &\quad
    +\sum_{d\in\operatorname{supp}(E)}q^{-\lambda_d(A)}
    \Bigl[
      \mathbf1_{\{c\prec d\}}
      +(ut+rq^{1-|E|+\lambda_d(A)+v_d})
        \sigma_{A+E}T_c(1-ut\sigma_{A+E}S)^{-1}
    \Bigr]
    G_{rq^{1-|E|}}(\sigma_{A+E}S)\Phi_d(A,E-d;r).
    \label{eqRisePhiGExpansion}
\end{align}

We next give the identities needed to simplify
\eqref{eqRisePhiGExpansion}.
Every index occurring in $S-T_c$ is smaller than every
index occurring in $T_c$.
Hence
\begin{alignat}{2}
  \sigma_c(S)T_c
  &=(S-T_c)T_c+qT_c^2
  &\quad\text{(by~\eqref{eqRefinedS45HeadElements})}&\notag\\
  &=qT_c(S-T_c)+qT_c^2=qT_cS
  &\quad\text{(by Proposition~\ref{propHeadLinearFunctional}
                   \textup{(1)})}&\rlap{.}
    \label{eqRiseSTCommutation}
\end{alignat}
For $h=0$, $(\sigma_cS)^hT_c=T_c(qS)^h$.
Assuming this identity for $h-1$, we have
\begin{alignat}{2}
  (\sigma_cS)^hT_c
  &=(\sigma_cS)^{h-1}\bigl(\sigma_c(S)T_c\bigr)
    =q(\sigma_cS)^{h-1}T_cS
  &\quad\text{(by~\eqref{eqRiseSTCommutation})}&\notag\\
  &=qT_c(qS)^{h-1}S=T_c(qS)^h
  &\quad\text{(by the induction hypothesis)}&\rlap{.}
    \label{eqRisePowerCommutation}
\end{alignat}

Since $G_r(S)$ is a polynomial in $S$, it commutes with
$1+rS$.
Proposition~\ref{propGIdentities} gives
\begin{alignat}{2}
  G_r(S)(1+rS)
  &=(1+rS)G_r(S)=G_{qr}(S)
  &\quad\text{(by Proposition~\ref{propGIdentities}
                   \textup{(1)})}&\notag\\
  (1-utS)^{-1}G_{qr}(S)
  &=(1-utS)^{-1}(1+rS)G_r(S)=G_r(qS)
  &\quad\text{(by Proposition~\ref{propGIdentities}
                   \textup{(1)--(2)})}&\rlap{.}
    \label{eqRiseGParameterRelations}
\end{alignat}
We can therefore calculate
\begin{alignat}{2}
  G_r(\sigma_cS)T_c
  &=\sum_{h=0}^{n}
    \frac{(-1)^hg_h(ut,r)}{(q;q)_h}(\sigma_cS)^hT_c
  &\quad\text{(by~\eqref{eqdefGr})}&\notag\\
  &=T_c\sum_{h=0}^{n}
    \frac{(-1)^hg_h(ut,r)}{(q;q)_h}(qS)^h
  &\quad\text{(by~\eqref{eqRisePowerCommutation})}&\notag\\
  &=T_cG_r(qS)
  &\quad\text{(by~\eqref{eqdefGr})}&\notag\\
  &=T_c(1-utS)^{-1}G_{qr}(S)
  &\quad\text{(by~\eqref{eqRiseGParameterRelations})}&\rlap{,}
    \label{eqRiseHeadGShift}
\end{alignat}
and
\begin{alignat}{2}
  &\quad \quad G_r(\sigma_cS)\bigl(1+r(S-T_c)-utT_c\bigr)
    \notag\\
  &=G_r(\sigma_cS)(1+rS)
    -(ut+r)G_r(\sigma_cS)T_c
  &&\notag\\
  &=G_r(\sigma_cS)(1+rS)
    -(ut+r)T_cG_r(qS)
  &\quad\text{(by~\eqref{eqRiseHeadGShift})}&\notag\\
  &=\Bigl(
      G_r(S)+(ut+r)T_c(1-utS)^{-1}G_r(S)
    \Bigr)(1+rS)
    -(ut+r)T_cG_r(qS)
  &\quad\text{(by Proposition~\ref{propGIdentities}
                   \textup{(3)})}&\notag\\
  &=G_{qr}(S)+(ut+r)T_cG_r(qS)-(ut+r)T_cG_r(qS)
  &\quad\text{(by~\eqref{eqRiseGParameterRelations})}&\notag\\
  &=G_{qr}(S).
    \label{eqRiseBranchGShift}
\end{alignat}

Apply $\sigma_{A+E}$ to
Proposition~\ref{propGIdentities}\textup{(3)}, with its
parameter equal to $r/q^{|E|}$.
Since $\sigma_{A+E}$ commutes with $\sigma_c$, 
Proposition~\ref{propGIdentities} (3) gives
\begin{equation}\label{eqRiseUndecoratedG}
  \Bigl(
      1+(ut+rq^{-|E|})
      \sigma_{A+E}T_c(1-ut\sigma_{A+E}S)^{-1}
    \Bigr)
    G_{r/q^{|E|}}(\sigma_{A+E}S)
  =G_{r/q^{|E|}}(\sigma_c\sigma_{A+E}S).
\end{equation}

Next replace $r$ by $r/q^{|E|}$ in
\eqref{eqRiseHeadGShift} and~\eqref{eqRiseBranchGShift},
and apply $\sigma_{A+E}$.
The resulting identities are
\begin{equation}\label{eqRiseGAfterSigma}
\begin{aligned}
  G_{r/q^{|E|}}(\sigma_c\sigma_{A+E}S)\sigma_{A+E}T_c
  &=\sigma_{A+E}T_c(1-ut\sigma_{A+E}S)^{-1}
      G_{rq^{1-|E|}}(\sigma_{A+E}S),\\
  G_{r/q^{|E|}}(\sigma_c\sigma_{A+E}S)
    \bigl(
      1+rq^{-|E|}\sigma_{A+E}(S-T_c)
      &-ut\sigma_{A+E}T_c
    \bigr)
  =G_{rq^{1-|E|}}(\sigma_{A+E}S).
\end{aligned}
\end{equation}
Using these identities, we obtain
\begin{alignat}{2}
  &\quad \quad G_{r/q^{|E|}}(\sigma_c\sigma_{A+E}S)\Gamma_{c,d}(r)
    \notag\\
  &=q^{-\lambda_d(A)}
    G_{r/q^{|E|}}(\sigma_c\sigma_{A+E}S)
    \Bigl[
      \mathbf1_{\{c\prec d\}}
      \bigl(
        1+rq^{-|E|}\sigma_{A+E}(S-T_c)
        -ut\sigma_{A+E}T_c
      \bigr)
      \notag\\
  &\qquad\qquad\qquad
      +(ut+rq^{1-|E|+\lambda_d(A)+v_d})\sigma_{A+E}T_c
    \Bigr]
  &\quad\text{(by~\eqref{eqRefinedS6GammaDefinition})}&\notag\\
  &=q^{-\lambda_d(A)}
    \Bigl[
      \mathbf1_{\{c\prec d\}}
      G_{rq^{1-|E|}}(\sigma_{A+E}S)
      \notag\\
  &\qquad
      +(ut+rq^{1-|E|+\lambda_d(A)+v_d})
       \sigma_{A+E}T_c(1-ut\sigma_{A+E}S)^{-1}
       G_{rq^{1-|E|}}(\sigma_{A+E}S)
    \Bigr]
  &\quad\text{(by~\eqref{eqRiseGAfterSigma})}&\notag\\
  &=q^{-\lambda_d(A)}
    \Bigl[
      \mathbf1_{\{c\prec d\}}
      +(ut+rq^{1-|E|+\lambda_d(A)+v_d})
       \sigma_{A+E}T_c(1-ut\sigma_{A+E}S)^{-1}
    \Bigr]
 G_{rq^{1-|E|}}(\sigma_{A+E}S).
    \label{eqRiseDecoratedG}
\end{alignat}

Substituting~\eqref{eqRiseUndecoratedG}
and~\eqref{eqRiseDecoratedG}
into~\eqref{eqRisePhiGExpansion}
proves~\eqref{eqRisePhiGTarget}, and hence the required recurrence.
\end{proof}

\begin{lemma}\label{lemPhiZeroRelation}
For $E\neq\varnothing$, we have
\[
\begin{alignedat}{2}
\sum_{A=c'A'}\Phi_{c'}(A',E;0) +\sum_{d \in \operatorname{supp}(E)}q^{-\lambda_d(A)}\Phi_d(A,E-d;0) &= q^{-|E|}\sum_{A=c'A'}\tau\Phi_{c'}(A',E;0) +\sum_{d \in \operatorname{supp}(E)}q^{v_d+1-|E|}\tau\Phi_d(A,E-d;0).
\end{alignedat}
\]
\end{lemma}

\begin{proof}
Since $1^-\prec z$ for every relevant letter $z$,
prepending $1^-$ to $A$ increases each $\rho_z$ by one.
Separating the terms that use the new alphabet entry $1$
from those that do not gives
\begin{equation}\label{eqRefinedS6InsertionFactorShift}
  [\rho+1]_q=q[\rho]_q+1,\qquad
  \varphi_{\mu,\epsilon}(\rho+1)
  =q^{-\mu}\varphi_{\mu,\epsilon}(\rho)
    +q^{-(1-\epsilon)(\mu-1)}
      \varphi_{\mu-1,\epsilon}(\rho+\epsilon),
  \qquad \mu\geq1.
\end{equation}
For $\epsilon=0$, the new entry can occur at most once.
For $\epsilon=1$, removing one occurrence leaves a choice
of $\mu-1$ entries from the same enlarged alphabet.

In the following products, $\rho_z$ and $n_{>z}$ are
computed from $A,I$, and $\mu_z$ is the multiplicity in $E$.
Starting with $W_I(1^-A,E)$, replace
$[\rho_{\xi_i}+1]_q$ by $q[\rho_{\xi_i}]_q$ and
$\varphi_{\mu_d,\epsilon_d}(\rho_d+1)$ by
$q^{-\mu_d}\varphi_{\mu_d,\epsilon_d}(\rho_d)$,
one factor at a time in increasing label order.
After all replacements, \eqref{eqRefinedS45InsertionProduct} gives
\[
  \prod_{i\in I}q[\rho_{\xi_i}]_q
  \prod_{d\in\operatorname{supp}(E)}
  q^{\mu_d(n_{>d}-1)}
  \varphi_{\mu_d,\epsilon_d}(\rho_d)
  =q^{|I|-|E|}W_I(A,E).
\]
We calculate the product before each replacement minus the
product after it. At that step, the factors with smaller
labels have already been replaced, while those with larger
labels remain unchanged.

At the factor indexed by $\xi_i$, the difference is
$[\rho_{\xi_i}+1]_q-q[\rho_{\xi_i}]_q=1$.
By~\eqref{eqRefinedS45InsertionCounts},
\[
  \#\{k\in I\setminus\{i\}:\xi_k>z\}
  =n_{>z}-\mathbf1_{\{z<\xi_i\}},\qquad
  \lambda_z(1^-A)
  +\#\{k\in I\setminus\{i\}:\xi_k>z\}
  =\rho_z+1-\mathbf1_{\{z<\xi_i\}}.
\]
Thus the difference between the two products is
\begin{align*}
  \prod_{\substack{k\in I\\k<i}}q[\rho_{\xi_k}]_q
   \prod_{\substack{k\in I\\k>i}}[\rho_{\xi_k}+1]_q
   \prod_{\substack{z\in\operatorname{supp}(E)\\z<\xi_i}}
   q^{\mu_z(n_{>z}-1)}
   \varphi_{\mu_z,\epsilon_z}(\rho_z)
   \prod_{\substack{z\in\operatorname{supp}(E)\\z>\xi_i}}
   q^{\mu_zn_{>z}}
   \varphi_{\mu_z,\epsilon_z}(\rho_z+1)
  =q^{\#\{k\in I:k<i\}}
    W_{I\setminus\{i\}}(1^-A,E).
\end{align*}

At the factor indexed by $d\in\operatorname{supp}(E)$,
\eqref{eqRefinedS6InsertionFactorShift} gives
\begin{align*}
  q^{\mu_dn_{>d}}\varphi_{\mu_d,\epsilon_d}(\rho_d+1)
  -q^{\mu_d(n_{>d}-1)}\varphi_{\mu_d,\epsilon_d}(\rho_d)
  =q^{\mu_dn_{>d}-(1-\epsilon_d)(\mu_d-1)}
  \varphi_{\mu_d-1,\epsilon_d}(\rho_d+\epsilon_d).
\end{align*}
Also, \eqref{eqRefinedS45InsertionCounts} gives
\[
  \lambda_z(dA)+n_{>z}
  =
  \begin{cases}
    \rho_z,&z<d,\\
    \rho_d+\epsilon_d,&z=d,\\
    \rho_z+1,&z>d.
  \end{cases}
\]
Hence the difference between the two products is
\begin{alignat*}{2}
  &\prod_{\substack{i\in I\\\xi_i<d}}q[\rho_{\xi_i}]_q
   \prod_{\substack{i\in I\\\xi_i>d}}[\rho_{\xi_i}+1]_q
   \prod_{\substack{z\in\operatorname{supp}(E)\\z<d}}
   q^{\mu_z(n_{>z}-1)}
   \varphi_{\mu_z,\epsilon_z}(\rho_z)\\
  &\quad{}\times
   q^{\mu_dn_{>d}-(1-\epsilon_d)(\mu_d-1)}
   \varphi_{\mu_d-1,\epsilon_d}(\rho_d+\epsilon_d)
   \prod_{\substack{z\in\operatorname{supp}(E)\\z>d}}
   q^{\mu_zn_{>z}}
   \varphi_{\mu_z,\epsilon_z}(\rho_z+1)\\
  &=q^{\#\{i\in I:\xi_i<d\}-\#\{e\in E:e<d\}
        +n_{>d}-(1-\epsilon_d)(\mu_d-1)}
    W_I(dA,E-d)
  &\quad\text{(by~\eqref{eqRefinedS45InsertionProduct})}&\rlap{.}
\end{alignat*}
Since no head label equals $d$,
$\#\{i\in I:\xi_i<d\}+n_{>d}=|I|$.
Also,
$|E|=\#\{e\in E:e<d\}+\mu_d+\#\{e\in E:e>d\}$.
Thus
\begin{alignat*}{2}
  &\quad \#\{i\in I:\xi_i<d\}-\#\{e\in E:e<d\}
    +n_{>d}-(1-\epsilon_d)(\mu_d-1)\\
  &=|I|-\#\{e\in E:e<d\}-\mu_d+1
    +\epsilon_d(\mu_d-1)\\
  &=|I|-|E|+1+\#\{e\in E:e>d\}
    +\epsilon_d(\mu_d-1)\\
  &=|I|-|E|+1+v_d
  &\quad\text{(by~\eqref{eqRefinedS6RemainingCounts})}&\rlap{.}
\end{alignat*}
Summing the differences between consecutive products gives
\begin{equation}\label{eqRefinedS6InsertionProductShift}
  W_I(1^-A,E)
  =q^{|I|-|E|}W_I(A,E)
    +\sum_{i\in I}q^{\#\{k\in I:k<i\}}
      W_{I\setminus\{i\}}(1^-A,E)
      +\sum_{d\in\operatorname{supp}(E)}
      q^{|I|-|E|+1+v_d}W_I(dA,E-d).
\end{equation}

We next express the sum over $i\in I$ as right multiplication
by $utS$:
\begin{alignat}{2}
  \sum_{I\subseteq[n]}\sum_{i\in I}
    (ut)^{|I|}q^{\#\{k\in I:k<i\}}
    W_{I\setminus\{i\}}(1^-A,E)Q_I
  &=\sum_{I\subseteq[n]}\sum_{i\in I}
    (ut)^{|I|}
    W_{I\setminus\{i\}}(1^-A,E)
    Q_{I\setminus\{i\}}Q_i
  &\quad\text{(by~\eqref{eqRefinedS6OrderedProduct})}&\notag\\
  &=\sum_{I\subseteq[n]}\sum_{i\notin I}
    (ut)^{|I|+1}W_I(1^-A,E)Q_IQ_i
  &\quad\text{(by reindexing)}&\notag\\
  &=\sum_{I\subseteq[n]}
    (ut)^{|I|}W_I(1^-A,E)Q_I
    \left(ut\sum_{i=1}^{n}Q_i\right)
  &\quad\text{(by~\eqref{eqRefinedS6OrderedProduct})}&\notag\\
  &=\left(
      \sum_{I\subseteq[n]}(ut)^{|I|}W_I(1^-A,E)Q_I
    \right)utS
  &\quad\text{(by~\eqref{eqRefinedS45HeadElements})}&\rlap{.}
    \label{eqHeadRemovalSum}
\end{alignat}

For any word $cA$ and multiset $E$, we have
\begin{alignat}{2}
  \sum_{I\subseteq[n]}(ut)^{|I|}W_I(cA,E)Q_I
  &=\frac{K_{cA,E}(ut,0)}
    {\displaystyle
     \prod_{d\in\operatorname{supp}(E)}
       \wt(d)^{\mu_d}q^{\mu_d(\lambda_d(cA)-1)}}
  &\quad\text{(by Theorem~\ref{thmRootHeadInsertion})}&\notag\\
  &=G_0(\sigma_c\sigma_{A+E}S)
    \Phi_c(A,E;0)G_0(S)^{-1}
  &\quad\text{(by~\eqref{eqPhiNormalization})}&\rlap{.}
    \label{eqPhiWExpansion}
\end{alignat}
When $A=\varnothing$ and $E\neq\varnothing$,
$W_I(A,E)=0$: the largest letter among $E$ and the heads
has $\rho_z=0$, so its factor in
\eqref{eqRefinedS45InsertionProduct} vanishes.
This agrees with the empty sum over $A=c'A'$.

Since $\tau$ is an algebra automorphism fixing the scalar
parameters, \eqref{eqdefGr} and \eqref{eqRefinedS6TauDefinition} give
\begin{equation}\label{eqTauG0Shift}
  \tau\bigl(G_0(\sigma_{A+E}S)\bigr)
  =G_0\bigl(\tau(\sigma_{A+E}S)\bigr)
  =G_0(q\sigma_{A+E}S), \quad
  \tau\bigl(G_0(S)^{-1}\bigr)
  =\tau\bigl(G_0(S)\bigr)^{-1}
  =G_0(\tau S)^{-1}
  =G_0(qS)^{-1}
\end{equation}

We now calculate
\begin{flalign}
  &\quad G_0(q\sigma_{A+E}S)\Phi_{1^-}(A,E;0)G_0(qS)^{-1}
  &&\notag\\
  &=G_0(q\sigma_{A+E}S)\Phi_{1^-}(A,E;0)
    G_0(S)^{-1}(1-utS)
  &\mathllap{\text{(by Proposition~\ref{propGIdentities}\textup{(2)})}}
  &\notag\\
  &=G_0(\sigma_{1^-}\sigma_{A+E}S)\Phi_{1^-}(A,E;0)
    G_0(S)^{-1}(1-utS)
  &\mathllap{\text{(by~\eqref{eqRefinedS45HeadElements})}}
  &\notag\\
  &=\left(
      \sum_{I\subseteq[n]}(ut)^{|I|}W_I(1^-A,E)Q_I
    \right)(1-utS)
  &\mathllap{\text{(by~\eqref{eqPhiWExpansion})}}
  &\notag\\
  &=\sum_{I\subseteq[n]}(ut)^{|I|}
    \left(
      W_I(1^-A,E)
      -\sum_{i\in I}q^{\#\{k\in I:k<i\}}
        W_{I\setminus\{i\}}(1^-A,E)
    \right)Q_I
  &\mathllap{\text{(by~\eqref{eqHeadRemovalSum})}}
  &\notag\\
  &=\sum_{I\subseteq[n]}(ut)^{|I|}
      q^{|I|-|E|}W_I(A,E)Q_I
    +\sum_{d\in\operatorname{supp}(E)}
     \sum_{I\subseteq[n]}(ut)^{|I|}
      q^{|I|+v_d+1-|E|}W_I(dA,E-d)Q_I
  &\mathllap{\text{(by~\eqref{eqRefinedS6InsertionProductShift})}}
  &\notag\\
  &=q^{-|E|}
    \tau\left(
      \sum_{I\subseteq[n]}(ut)^{|I|}W_I(A,E)Q_I
    \right)
    +\sum_{d\in\operatorname{supp}(E)}q^{v_d+1-|E|}
    \tau\left(
      \sum_{I\subseteq[n]}(ut)^{|I|}
      W_I(dA,E-d)Q_I
    \right)
  &&\notag\\[-.5ex]
  &&\mathllap{\text{(by~\eqref{eqRefinedS6TauDefinition})}}
  &\notag\\
  &\mathrlap{{}=q^{-|E|}\sum_{\mathclap{A=c'A'}}
    \tau\left(
      G_0(\sigma_{c'}\sigma_{A'+E}S)
      \Phi_{c'}(A',E;0)G_0(S)^{-1}
    \right)
    +\sum_{\mathclap{d\in\operatorname{supp}(E)}}q^{v_d+1-|E|}
    \tau\left(
      G_0(\sigma_d\sigma_{A+(E-d)}S)
      \Phi_d(A,E-d;0)G_0(S)^{-1}
    \right)}
  &&\notag\\[-.5ex]
  &&\mathllap{\text{(by~\eqref{eqPhiWExpansion})}}
  &\notag\\
  &\mathrlap{{}=q^{-|E|}\sum_{\mathclap{A=c'A'}}
    \tau\left(
      G_0(\sigma_{A+E}S)
      \Phi_{c'}(A',E;0)G_0(S)^{-1}
    \right)
    +\sum_{\mathclap{d\in\operatorname{supp}(E)}}q^{v_d+1-|E|}
    \tau\left(
      G_0(\sigma_{A+E}S)
      \Phi_d(A,E-d;0)G_0(S)^{-1}
    \right)}
  &&\notag\\[-.5ex]
  &&\mathllap{\text{(by~\eqref{eqRefinedS6SigmaWords})}}
  &\notag\\
  &\mathrlap{{}=q^{-|E|}\sum_{\mathclap{A=c'A'}}
    G_0(q\sigma_{A+E}S)
    \tau\Phi_{c'}(A',E;0)G_0(qS)^{-1}
    +\sum_{\mathclap{d\in\operatorname{supp}(E)}}q^{v_d+1-|E|}
    G_0(q\sigma_{A+E}S)
    \tau\Phi_d(A,E-d;0)G_0(qS)^{-1}}
  &&\notag\\[-.5ex]
  &&\mathllap{\text{(by~\eqref{eqTauG0Shift})}}
  &\notag\\
  &=G_0(q\sigma_{A+E}S)
    \left(
      q^{-|E|}\sum_{A=c'A'}\tau\Phi_{c'}(A',E;0)
      +\sum_{d\in\operatorname{supp}(E)}
       q^{v_d+1-|E|}\tau\Phi_d(A,E-d;0)
    \right)
    G_0(qS)^{-1}.
  &&\label{first-last-expressions}
\end{flalign}

Since $1^-\prec d$, \eqref{eqRefinedS6GammaDefinition}
gives $\Gamma_{1^-,d}(0)=q^{-\lambda_d(A)}$.
The factors $G_0(q\sigma_{A+E}S)$ and $G_0(qS)$ are invertible,
so
\begin{alignat*}{2}
  &\quad \sum_{A=c'A'}\Phi_{c'}(A',E;0)
    +\sum_{d\in\operatorname{supp}(E)}
      q^{-\lambda_d(A)}\Phi_d(A,E-d;0)\\
  &=\Phi_{1^-}(A,E;0)
  &\quad\text{(by Lemma~\ref{lemPhiRecurrence})}&\\
  &=q^{-|E|}\sum_{A=c'A'}\tau\Phi_{c'}(A',E;0)
    +\sum_{d\in\operatorname{supp}(E)}
      q^{v_d+1-|E|}\tau\Phi_d(A,E-d;0)
  &\quad\text{(by~\eqref{first-last-expressions})}&\rlap{,}
\end{alignat*}
as required.
\end{proof}

\begin{proof}[Proof of Lemma~\ref{lemHomogeneousFactor}]
We take coefficients in the basis $\{Q_I\}$ throughout,
so nonzero scalar polynomials can be canceled.
For $j=0$, $[H]_0$ has degree zero in $ut,r$, and the
conclusion is immediate.
For $j\geq1$, we have $\tau[H]_j=q^j[H]_j$, so the given
equation in $Q$-degree $j$ becomes
\begin{equation}\label{eqRefinedS6HomogeneousEquation}
  (qr+utq^j)[H(qr)]_j
  =q^j(qr+ut)[H(r)]_j.
\end{equation}

Suppose that~\eqref{eqRefinedS6HomogeneousEquation} holds.
Substituting $r=-utq^{h-1}$ gives, for $0\leq h<j$,
\[
  ut(q^j-q^h)[H(-utq^h)]_j
  =q^jut(1-q^h)[H(-utq^{h-1})]_j.
\]
For $h=0$, the right-hand side is zero.
Since $ut(q^j-q^h)\neq0$ for $0\leq h<j$, we obtain successively
$[H(-utq^h)]_j=0$, for $0\leq h<j$.
As a polynomial in $r$, $[H(r)]_j$ has degree at most $j$.
These $j$ distinct roots therefore give
\[
  [H(r)]_j
  =\left(\prod_{h=0}^{j-1}(r+utq^h)\right)C_j
  =g_j(ut,r)C_j.
\]
Here $C_j$ is the coefficient of $r^j$ in $[H(r)]_j$.
Since every monomial has total degree $j$ in $ut,r$,
this coefficient is independent of $ut,r$.

Conversely, suppose that $[H(r)]_j=g_j(ut,r)C_j$,
where $C_j$ is independent of $ut,r$.
Then
\begin{alignat*}{2}
  (qr+utq^j)[H(qr)]_j
  &=(qr+utq^j)\prod_{h=0}^{j-1}(qr+utq^h)C_j
  &\quad\text{(by~\eqref{eqRefinedS6gDefinition})}&\\
  &=q^j(qr+ut)(r+utq^{j-1})
    \prod_{h=0}^{j-2}(r+utq^h)C_j\\
  &=q^j(qr+ut)g_j(ut,r)C_j
  &\quad\text{(by~\eqref{eqRefinedS6gDefinition})}&\\
  &=q^j(qr+ut)[H(r)]_j,
\end{alignat*}
which proves~\eqref{eqRefinedS6HomogeneousEquation}.
\end{proof}

\begin{lemma}\label{lemPhiWeightedRelation}
Let $E\neq\varnothing$.
Suppose that, for every $j\geq0$,
$[G_{r/q^{|E|}}(\sigma_{A+E}S)\Phi_{c'}(A',E;r)]_j$
for $A=c'A'$, and
$[G_{rq^{1-|E|}}(\sigma_{A+E}S)\Phi_d(A,E-d;r)]_j$
for $d\in\operatorname{supp}(E)$,
are divisible by $g_j(ut,r)$, with quotients independent of $ut,r$.
Then
\begin{alignat*}{2}
&(1-q^{-|E|}\tau) \sum_{A=c'A'}\Phi_{c'}(A',E;r) +\sum_{d \in \operatorname{supp}(E)} \bigl(q^{-\lambda_d(A)}-q^{v_d+1-|E|}\tau\bigr) \Phi_d(A,E-d;r)\\
&= r(\sigma_{A+E}S) \sum_{d \in \operatorname{supp}(E)} q^{v_d+1-|E|} \bigl(q^{1-|E|}\tau-1\bigr)\Phi_d(A,E-d;r).
\end{alignat*}
\end{lemma}

\begin{proof}
By~\eqref{eqdefGr} and~\eqref{eqRefinedS6TauDefinition},
$\tau(G_r(\sigma_{A+E}S))=G_r(q\sigma_{A+E}S)$.
Hence
\begin{flalign}
  &\quad \left(1-q^{-|E|}\tau
    +q^{-|E|}(\sigma_{A+E}S)(r+ut\tau)\right)
    \left(G_{r/q^{|E|}}(\sigma_{A+E}S)
    \Phi_{c'}(A',E;r)\right)
  &&\notag\\
  &=\Bigl(
    (1+rq^{-|E|}\sigma_{A+E}S)
      G_{r/q^{|E|}}(\sigma_{A+E}S)
    -q^{-|E|}(1-ut\sigma_{A+E}S)
      G_{r/q^{|E|}}(q\sigma_{A+E}S)\tau
    \Bigr)\Phi_{c'}(A',E;r)
  &&\notag\\
  &=(1+rq^{-|E|}\sigma_{A+E}S)
    G_{r/q^{|E|}}(\sigma_{A+E}S)
    (1-q^{-|E|}\tau)\Phi_{c'}(A',E;r)
  &\mathllap{\text{(by Proposition~\ref{propGIdentities}\textup{(2)})}}
  &\notag\\
  &=G_{rq^{1-|E|}}(\sigma_{A+E}S)
    (1-q^{-|E|}\tau)\Phi_{c'}(A',E;r),
  &\mathllap{\text{(by Proposition~\ref{propGIdentities}\textup{(1)})}}
  &\notag\\[1ex]
  &\quad \left(q^{-\lambda_d(A)}-q^{v_d+1-|E|}\tau
    +q^{v_d+1-|E|}(\sigma_{A+E}S)(r+ut\tau)\right)
    \left(G_{rq^{1-|E|}}(\sigma_{A+E}S)
    \Phi_d(A,E-d;r)\right)
  &&\notag\\
  &\mathrlap{{}=\Bigl(
    (q^{-\lambda_d(A)}+rq^{v_d+1-|E|}\sigma_{A+E}S)
      G_{rq^{1-|E|}}(\sigma_{A+E}S)
    -q^{v_d+1-|E|}(1-ut\sigma_{A+E}S)
      G_{rq^{1-|E|}}(q\sigma_{A+E}S)\tau
    \Bigr)\Phi_d(A,E-d;r)}
  &&\notag\\
  &=G_{rq^{1-|E|}}(\sigma_{A+E}S)
    \Bigl(q^{-\lambda_d(A)}
      +rq^{v_d+1-|E|}\sigma_{A+E}S
      -q^{v_d+1-|E|}
       (1+rq^{1-|E|}\sigma_{A+E}S)\tau
    \Bigr)\Phi_d(A,E-d;r)
  &&\notag\\[-.5ex]
  &&\mathllap{\text{(by Proposition~\ref{propGIdentities}\textup{(2)})}}
  &\notag\\
  &=G_{rq^{1-|E|}}(\sigma_{A+E}S)
    \Bigl(q^{-\lambda_d(A)}-q^{v_d+1-|E|}\tau
      -rq^{v_d+1-|E|}(\sigma_{A+E}S)
       (q^{1-|E|}\tau-1)
    \Bigr)\Phi_d(A,E-d;r).
  &&\label{eqRefinedS6WeightedExpansion}
\end{flalign}

Sum the first calculation over $A=c'A'$ and the second over
$d\in\operatorname{supp}(E)$.
For $I\subseteq[n]$ and $j=|I|$, we have
\begin{align*}
  \tau\bigl(g_j(ut,r)Q_I\bigr)
  =q^j& g_j(ut,r)Q_I,\\
  (\sigma_{A+E}S)(r+ut\tau)\bigl(g_j(ut,r)Q_I\bigr)
  &=(r+utq^j)g_j(ut,r)(\sigma_{A+E}S)Q_I
  =g_{j+1}(ut,r)(\sigma_{A+E}S)Q_I.
\end{align*}
The coefficients of $(\sigma_{A+E}S)Q_I$ are independent of
$ut,r$, and its terms have $Q$-degree $j+1$.
Thus, by hypothesis, the coefficient of each $Q_I$ in the
summed left-hand sides is divisible by $g_{|I|}(ut,r)$,
with quotient independent of $ut,r$.

At $r=0$, Lemma~\ref{lemPhiZeroRelation} gives
\begin{align*}
  G_0(\sigma_{A+E}S)
    \left(
      (1-q^{-|E|}\tau)
        \sum_{A=c'A'}\Phi_{c'}(A',E;0)
      +\sum_{d\in\operatorname{supp}(E)}
        (q^{-\lambda_d(A)}-q^{v_d+1-|E|}\tau)
        \Phi_d(A,E-d;0)
    \right) = 0.
\end{align*}
Each quotient is independent of $ut,r$.
Since the sum vanishes at $ut=1,r=0$ and
$g_j(1,0)=q^{\binom{j}{2}}\neq0$, every quotient is zero.
Thus the summed left-hand sides vanish for all $ut,r$.
Adding the corresponding right-hand sides gives
\[
  0=G_{rq^{1-|E|}}(\sigma_{A+E}S)
  \Biggl[
  \begin{aligned}
    &(1-q^{-|E|}\tau)
      \sum_{A=c'A'}\Phi_{c'}(A',E;r)
      +\sum_{d\in\operatorname{supp}(E)}
       \bigl(q^{-\lambda_d(A)}-q^{v_d+1-|E|}\tau\bigr)
       \Phi_d(A,E-d;r)\\
    &\qquad
      -r(\sigma_{A+E}S)
       \sum_{d\in\operatorname{supp}(E)}
       q^{v_d+1-|E|}
       (q^{1-|E|}\tau-1)\Phi_d(A,E-d;r)
  \end{aligned}
  \Biggr].
\]
The factor $G_{rq^{1-|E|}}(\sigma_{A+E}S)$ has constant
term $1$, so it is invertible.
Canceling it and moving the last sum to the right gives
\begin{flalign*}
  &\quad (1-q^{-|E|}\tau)
    \sum_{A=c'A'}\Phi_{c'}(A',E;r)
    +\sum_{d\in\operatorname{supp}(E)}
     \bigl(q^{-\lambda_d(A)}-q^{v_d+1-|E|}\tau\bigr)
     \Phi_d(A,E-d;r)
  &&\\
  &=r(\sigma_{A+E}S)
    \sum_{d\in\operatorname{supp}(E)}
    q^{v_d+1-|E|}
    (q^{1-|E|}\tau-1)\Phi_d(A,E-d;r),
  &&
\end{flalign*}
which is the asserted relation.
\end{proof}

\begin{proof}[Proof of Lemma~\ref{lemPhiOperatorCalculation}]
We first prove
\begin{equation}\label{eqRefinedS6PhiReduction}
  qr(1-ut\sigma_c\sigma_{A+E}S)\Phi_c(A,E;qr)
   +ut(1+rq^{1-|E|}\sigma_c\sigma_{A+E}S)
     \tau\Phi_c(A,E;qr)
  =(qr+ut)\tau\Phi_c(A,E;r).
\end{equation}

By~\eqref{eqdefGr} and~\eqref{eqRefinedS6TauDefinition},
$\tau(G_r(\sigma_{A+E}S))=G_r(q\sigma_{A+E}S)$.
For $A=c'A'$, the hypotheses give
\begin{flalign*}
  &\quad (qr+ut)G_{r/q^{|E|}}(q\sigma_{A+E}S)
    \tau\Phi_{c'}(A',E;r)
  &&\\
  &=(qr+ut\tau)
    \bigl(G_{rq^{1-|E|}}(\sigma_{A+E}S)
      \Phi_{c'}(A',E;qr)\bigr)
  &\mathllap{\text{(by Lemma~\ref{lemHomogeneousFactor})}}&\\
  &=G_{r/q^{|E|}}(q\sigma_{A+E}S)
    \Bigl(
      qr(1-ut\sigma_{A+E}S)\Phi_{c'}(A',E;qr)
      +ut(1+rq^{1-|E|}\sigma_{A+E}S)
        \tau\Phi_{c'}(A',E;qr)
    \Bigr)
  &&\\[-.5ex]
  &&\mathllap{\text{(by Proposition~\ref{propGIdentities}
    \textup{(1)--(2)})}}&
\end{flalign*}
Canceling $G_{r/q^{|E|}}(q\sigma_{A+E}S)$ gives the first
equality below.
The same calculation for
$G_{rq^{1-|E|}}(\sigma_{A+E}S)\Phi_d(A,E-d;r)$,
canceling $G_{rq^{1-|E|}}(q\sigma_{A+E}S)$, gives the second:
\begin{flalign}
  qr(1-ut\sigma_{A+E}S)\Phi_{c'}(A',E;qr)
    +ut(1+rq^{1-|E|}\sigma_{A+E}S)
      \tau\Phi_{c'}(A',E;qr)
  &=(qr+ut)\tau\Phi_{c'}(A',E;r),
  &&\notag\\[.5ex]
  qr(1-ut\sigma_{A+E}S)\Phi_d(A,E-d;qr)
    +ut(1+rq^{2-|E|}\sigma_{A+E}S)
      \tau\Phi_d(A,E-d;qr)
  &=(qr+ut)\tau\Phi_d(A,E-d;r).
  &&\label{eqRiseChildEquations}
\end{flalign}

By Proposition~\ref{propHeadLinearFunctional}\textup{(1)},
$(S-T_c)T_c=qT_c(S-T_c)$.
Together with~\eqref{eqRefinedS45HeadElements}, this gives
\begin{equation}\label{eqRiseOperatorCommutation}
  \sigma_c(S)T_c=qT_cS,\qquad
  \sigma_c(S)(S-T_c)=(S-T_c)S,\qquad
  \sigma_c(S)=S+(q-1)T_c.
\end{equation}
We apply $\sigma_{A+E}$ to these equalities below.

By~\eqref{eqRefinedS6GammaDefinition}
and~\eqref{eqRefinedS6TauDefinition},
\begin{flalign}
  \Gamma_{c,d}(qr)-\tau\Gamma_{c,d}(r)
  &=-ut(q-1)q^{-\lambda_d(A)}
    (1-\mathbf1_{\{c\prec d\}})\sigma_{A+E}T_c,
  &&\notag\\[.5ex]
  \tau\Gamma_{c,d}(qr)-\tau\Gamma_{c,d}(r)
  &=(q-1)r\Bigl(
    \mathbf1_{\{c\prec d\}}q^{1-|E|-\lambda_d(A)}
      \sigma_{A+E}(S-T_c)
    +q^{v_d+2-|E|}\sigma_{A+E}T_c
    \Bigr).
  &&\label{eqRiseGammaDifferences}
\end{flalign}
Also,
\begin{flalign}
  &\quad (\sigma_c\sigma_{A+E}S)\Gamma_{c,d}(qr)
    -\tau\Gamma_{c,d}(r)\sigma_{A+E}S
  &&\notag\\
  &=\mathbf1_{\{c\prec d\}}q^{-\lambda_d(A)}
      \sigma_{A+E}\bigl(\sigma_c(S)-S\bigr)
    +\mathbf1_{\{c\prec d\}}rq^{1-|E|-\lambda_d(A)}
      \sigma_{A+E}\bigl(
        \sigma_c(S)(S-T_c)-(S-T_c)S
      \bigr)
  &&\notag\\
  &\quad
    +utq^{-\lambda_d(A)}(1-\mathbf1_{\{c\prec d\}})
      \sigma_{A+E}\bigl(\sigma_c(S)T_c-qT_cS\bigr)
    +rq^{v_d+2-|E|}
      \sigma_{A+E}\bigl(\sigma_c(S)T_c-T_cS\bigr)
  &&\notag\\[-.5ex]
  &&\mathllap{\text{(by~\eqref{eqRefinedS6GammaDefinition})}}
  &\notag\\
  &=(q-1)\sigma_{A+E}T_c
    \bigl(
      \mathbf1_{\{c\prec d\}}q^{-\lambda_d(A)}
      +rq^{v_d+2-|E|}\sigma_{A+E}S
    \bigr),
  &\mathllap{\text{(by~\eqref{eqRiseOperatorCommutation})}}
  &\label{eqRiseGammaCommutatorOne}\\[1ex]
  &\quad (\sigma_c\sigma_{A+E}S)\tau\Gamma_{c,d}(qr)
    -q\tau\Gamma_{c,d}(r)\sigma_{A+E}S
  &&\notag\\
  &=\mathbf1_{\{c\prec d\}}q^{-\lambda_d(A)}
      \sigma_{A+E}\bigl(\sigma_c(S)-qS\bigr)
    +\mathbf1_{\{c\prec d\}}rq^{2-|E|-\lambda_d(A)}
      \sigma_{A+E}\bigl(
        \sigma_c(S)(S-T_c)-(S-T_c)S
      \bigr)
  &&\notag\\
  &\quad
    +utq^{1-\lambda_d(A)}(1-\mathbf1_{\{c\prec d\}})
      \sigma_{A+E}\bigl(\sigma_c(S)T_c-qT_cS\bigr)
    +rq^{v_d+3-|E|}
      \sigma_{A+E}\bigl(\sigma_c(S)T_c-T_cS\bigr)
  &&\notag\\[-.5ex]
  &&\mathllap{\text{(by~\eqref{eqRefinedS6GammaDefinition})}}
  &\notag\\
  &=(q-1)\Bigl(
      -\mathbf1_{\{c\prec d\}}q^{-\lambda_d(A)}
        \sigma_{A+E}(S-T_c)
      +rq^{v_d+3-|E|}
        (\sigma_{A+E}T_c)\sigma_{A+E}S
    \Bigr).
  &\mathllap{\text{(by~\eqref{eqRiseOperatorCommutation})}}
  &\label{eqRiseGammaCommutatorTwo}
\end{flalign}

Consequently,
\begin{flalign}
  &\quad (1-ut\sigma_c\sigma_{A+E}S)\Gamma_{c,d}(qr)
    -\tau\Gamma_{c,d}(r)(1-ut\sigma_{A+E}S)
  &&\notag\\
  &=-ut(q-1)q^{-\lambda_d(A)}
      (1-\mathbf1_{\{c\prec d\}})\sigma_{A+E}T_c
    -ut\Bigl(
      (\sigma_c\sigma_{A+E}S)\Gamma_{c,d}(qr)
      -\tau\Gamma_{c,d}(r)\sigma_{A+E}S
    \Bigr)
  &&\notag\\[-.5ex]
  &&\mathllap{\text{(by~\eqref{eqRiseGammaDifferences})}}
  &\notag\\
  &=-ut(q-1)\sigma_{A+E}T_c
    \bigl(q^{-\lambda_d(A)}
      +rq^{v_d+2-|E|}\sigma_{A+E}S\bigr),
  &\mathllap{\text{(by~\eqref{eqRiseGammaCommutatorOne})}}
  &\label{eqRefinedS6GammaShiftOne}\\[1ex]
  &\quad (1+rq^{1-|E|}\sigma_c\sigma_{A+E}S)
      \tau\Gamma_{c,d}(qr)
    -\tau\Gamma_{c,d}(r)
      (1+rq^{2-|E|}\sigma_{A+E}S)
  &&\notag\\
  &=(q-1)r\Bigl(
      \mathbf1_{\{c\prec d\}}q^{1-|E|-\lambda_d(A)}
        \sigma_{A+E}(S-T_c)
      +q^{v_d+2-|E|}\sigma_{A+E}T_c
    \Bigr)
  &&\notag\\
  &\quad
    +rq^{1-|E|}\Bigl(
      (\sigma_c\sigma_{A+E}S)\tau\Gamma_{c,d}(qr)
      -q\tau\Gamma_{c,d}(r)\sigma_{A+E}S
    \Bigr)
  &\mathllap{\text{(by~\eqref{eqRiseGammaDifferences})}}
  &\notag\\
  &=rq(q-1)q^{v_d+1-|E|}\sigma_{A+E}T_c
    (1+rq^{2-|E|}\sigma_{A+E}S).
  &\mathllap{\text{(by~\eqref{eqRiseGammaCommutatorTwo})}}
  &\label{eqRefinedS6GammaShiftTwo}
\end{flalign}

For $A=c'A'$, we obtain
\begin{flalign}
  &\quad qr(1-ut\sigma_c\sigma_{A+E}S)\Phi_{c'}(A',E;qr)
    +ut(1+rq^{1-|E|}\sigma_c\sigma_{A+E}S)
      \tau\Phi_{c'}(A',E;qr)
  &&\notag\\
  &=qr(1-ut\sigma_{A+E}S)\Phi_{c'}(A',E;qr)
    +ut(1+rq^{1-|E|}\sigma_{A+E}S)
      \tau\Phi_{c'}(A',E;qr)
  &&\notag\\
  &\quad
    -utrq(q-1)\sigma_{A+E}T_c
      (1-q^{-|E|}\tau)\Phi_{c'}(A',E;qr)
  &\mathllap{\text{(by~\eqref{eqRiseOperatorCommutation})}}
  &\notag\\
  &=(qr+ut)\tau\Phi_{c'}(A',E;r)
    -utrq(q-1)\sigma_{A+E}T_c
      (1-q^{-|E|}\tau)\Phi_{c'}(A',E;qr).
  &\mathllap{\text{(by~\eqref{eqRiseChildEquations})}}
  &\label{eqRiseUndecoratedCalculation}
\end{flalign}

For $d\in\operatorname{supp}(E)$, we have
\begin{flalign}
  &\quad qr(1-ut\sigma_c\sigma_{A+E}S)
      \Gamma_{c,d}(qr)\Phi_d(A,E-d;qr)
    +ut(1+rq^{1-|E|}\sigma_c\sigma_{A+E}S)
      \tau\Gamma_{c,d}(qr)\tau\Phi_d(A,E-d;qr)
  &&\notag\\
  &=qr\tau\Gamma_{c,d}(r)(1-ut\sigma_{A+E}S)
      \Phi_d(A,E-d;qr)
    +ut(1+rq^{1-|E|}\sigma_c\sigma_{A+E}S)
      \tau\Gamma_{c,d}(qr)\tau\Phi_d(A,E-d;qr)
  &&\notag\\
  &\quad
    -utrq(q-1)\sigma_{A+E}T_c
      \bigl(q^{-\lambda_d(A)}
        +rq^{v_d+2-|E|}\sigma_{A+E}S\bigr)
      \Phi_d(A,E-d;qr)
  &\mathllap{\text{(by~\eqref{eqRefinedS6GammaShiftOne})}}
  &\notag\\
  &=\tau\Gamma_{c,d}(r)
    \Bigl(
      qr(1-ut\sigma_{A+E}S)\Phi_d(A,E-d;qr)
      +ut(1+rq^{2-|E|}\sigma_{A+E}S)
        \tau\Phi_d(A,E-d;qr)
    \Bigr)
  &&\notag\\
  &\quad
    -utrq(q-1)\sigma_{A+E}T_c
      \bigl(q^{-\lambda_d(A)}
        +rq^{v_d+2-|E|}\sigma_{A+E}S\bigr)
      \Phi_d(A,E-d;qr)
  &&\notag\\
  &\quad
    +utrq(q-1)q^{v_d+1-|E|}\sigma_{A+E}T_c
      (1+rq^{2-|E|}\sigma_{A+E}S)
      \tau\Phi_d(A,E-d;qr)
  &\mathllap{\text{(by~\eqref{eqRefinedS6GammaShiftTwo})}}
  &\notag\\
  &=(qr+ut)\tau\Gamma_{c,d}(r)\tau\Phi_d(A,E-d;r)
  &&\notag\\
  &\quad
    -utrq(q-1)\sigma_{A+E}T_c
      \bigl(q^{-\lambda_d(A)}
        +rq^{v_d+2-|E|}\sigma_{A+E}S\bigr)
      \Phi_d(A,E-d;qr)
  &&\notag\\
  &\quad
    +utrq(q-1)q^{v_d+1-|E|}\sigma_{A+E}T_c
      (1+rq^{2-|E|}\sigma_{A+E}S)
      \tau\Phi_d(A,E-d;qr)
  &\mathllap{\text{(by~\eqref{eqRiseChildEquations})}}
  &\notag\\
  &=(qr+ut)\tau\bigl(
      \Gamma_{c,d}(r)\Phi_d(A,E-d;r)
    \bigr)
  &&\notag\\
  &\quad
    -utrq(q-1)\sigma_{A+E}T_c
    \Bigl(
      q^{-\lambda_d(A)}-q^{v_d+1-|E|}\tau
      -rq^{v_d+2-|E|}\sigma_{A+E}S
        (q^{1-|E|}\tau-1)
    \Bigr)\Phi_d(A,E-d;qr).
  &&\label{eqRiseDecoratedRemainder}
\end{flalign}

Sum~\eqref{eqRiseUndecoratedCalculation} over $A=c'A'$
and~\eqref{eqRiseDecoratedRemainder} over
$d\in\operatorname{supp}(E)$.
Lemma~\ref{lemPhiRecurrence} gives
\begin{flalign*}
  &\quad qr(1-ut\sigma_c\sigma_{A+E}S)\Phi_c(A,E;qr)
    +ut(1+rq^{1-|E|}\sigma_c\sigma_{A+E}S)
      \tau\Phi_c(A,E;qr)
  &&\\
  &=(qr+ut)\tau\Phi_c(A,E;r)
    -utrq(q-1)\sigma_{A+E}T_c
    \Biggl(
      (1-q^{-|E|}\tau)
        \sum_{A=c'A'}\Phi_{c'}(A',E;qr)
  &&\\
  &\qquad
      +\sum_{d\in\operatorname{supp}(E)}
       \Bigl(
         q^{-\lambda_d(A)}-q^{v_d+1-|E|}\tau
         -rq^{v_d+2-|E|}\sigma_{A+E}S
           (q^{1-|E|}\tau-1)
       \Bigr)\Phi_d(A,E-d;qr)
    \Biggr)
  &&\\
  &=(qr+ut)\tau\Phi_c(A,E;r).
  &\mathllap{\text{(by Lemma~\ref{lemPhiWeightedRelation} at $qr$)}}
  &
\end{flalign*}
This proves~\eqref{eqRefinedS6PhiReduction}.

Finally,
\begin{flalign*}
  &\quad (qr+ut\tau)
    \bigl(
      G_{rq^{1-|E|}}(\sigma_c\sigma_{A+E}S)
      \Phi_c(A,E;qr)
    \bigr)
  &&\\
  &=qrG_{rq^{1-|E|}}(\sigma_c\sigma_{A+E}S)
      \Phi_c(A,E;qr)
    +utG_{rq^{1-|E|}}(q\sigma_c\sigma_{A+E}S)
      \tau\Phi_c(A,E;qr)
  &&\\[-.5ex]
  &&\mathllap{\text{(by~\eqref{eqRefinedS6TauDefinition})}}
  &\\
  &=G_{r/q^{|E|}}(q\sigma_c\sigma_{A+E}S)
    \Bigl(
      qr(1-ut\sigma_c\sigma_{A+E}S)\Phi_c(A,E;qr)
      +ut(1+rq^{1-|E|}\sigma_c\sigma_{A+E}S)
        \tau\Phi_c(A,E;qr)
    \Bigr)
  &&\\[-.5ex]
  &&\mathllap{\text{(by Proposition~\ref{propGIdentities}
    \textup{(1)--(2)})}}
  &\\
  &=(qr+ut)G_{r/q^{|E|}}(q\sigma_c\sigma_{A+E}S)
    \tau\Phi_c(A,E;r)
  &\mathllap{\text{(by~\eqref{eqRefinedS6PhiReduction})}}
  &\\
  &=(qr+ut)\tau\bigl(
      G_{r/q^{|E|}}(\sigma_c\sigma_{A+E}S)
      \Phi_c(A,E;r)
    \bigr).
  &\mathllap{\text{(by~\eqref{eqRefinedS6TauDefinition})}}
  &
\end{flalign*}
\end{proof}

\section{Proofs of Lemmas~\ref{lemMarkedDeletionIdentities} and~\ref{lemCyclicTransferVector}}
\label{app:deletion-cyclic}
This appendix proves Lemmas~\ref{lemMarkedDeletionIdentities} and
\ref{lemCyclicTransferVector}.

\begin{proof}[Proof of Lemma~\ref{lemMarkedDeletionIdentities}]
We display the ordered alphabet whenever it changes.
We first prove
\begin{equation}\label{eqBExchangeRelations}
\begin{aligned}
  \mathscr R_i(s_iX_n)\mathscr R_i(X_n)
  &=\mathrm{id},\\
  \mathscr R_i(s_{i+1}s_iX_n)
  \mathscr R_{i+1}(s_iX_n)\mathscr R_i(X_n)
  &=
  \mathscr R_{i+1}(s_is_{i+1}X_n)
  \mathscr R_i(s_{i+1}X_n)\mathscr R_{i+1}(X_n).
\end{aligned}
\end{equation}
By~\eqref{eqAdjacentExchangeBlock}, these operators commute
with multiplication by monomials in the other coordinates.
It therefore suffices to consider the two or three indicated
coordinates.

For the first equality, both $1$ and $z_iz_{i+1}$ are fixed.
For $k=i,i+1$, direct substitution in
\eqref{eqAdjacentExchangeBlock} gives
\[
  \mathscr R_i(s_iX_n)\mathscr R_i(X_n)z_k
  =
  \frac{(1-q)^2x_ix_{i+1}-q(x_i-x_{i+1})^2}
       {(x_i-qx_{i+1})(x_{i+1}-qx_i)}z_k
  =z_k.
\]

For the second equality, we use
\begin{align}
  (1-q)x_{i+1}(x_i-qx_{i+2})
    +q(x_i-x_{i+1})(x_{i+1}-x_{i+2})
  &=(x_i-qx_{i+1})(x_{i+1}-qx_{i+2}),
  \notag\\
  (x_i-qx_{i+2})(x_{i+1}-x_{i+2})
    +(1-q)x_{i+2}(x_i-x_{i+1})
  &=(x_i-x_{i+2})(x_{i+1}-qx_{i+2}),
  \notag\\
  (x_i-x_{i+1})(x_i-qx_{i+2})
    +(1-q)x_i(x_{i+1}-x_{i+2})
  &=(x_i-x_{i+2})(x_i-qx_{i+1}).
  \label{eqBExchangeCoefficients}
\end{align}
Substituting~\eqref{eqAdjacentExchangeBlock} three times
and using~\eqref{eqBExchangeCoefficients}, both products
in the second equality of~\eqref{eqBExchangeRelations} give
\begin{align*}
  z_i
  &\longmapsto
    \frac{(1-q)x_{i+2}}{x_i-qx_{i+2}}z_i
    +\frac{(1-q)x_{i+1}(x_i-x_{i+2})}
      {(x_i-qx_{i+1})(x_i-qx_{i+2})}z_{i+1}
    +\frac{(x_i-x_{i+1})(x_i-x_{i+2})}
      {(x_i-qx_{i+1})(x_i-qx_{i+2})}z_{i+2},\\
  z_{i+2}
  &\longmapsto
    \frac{q^2(x_i-x_{i+2})(x_{i+1}-x_{i+2})}
      {(x_i-qx_{i+2})(x_{i+1}-qx_{i+2})}z_i
    +\frac{q(1-q)x_{i+1}(x_i-x_{i+2})}
      {(x_i-qx_{i+2})(x_{i+1}-qx_{i+2})}z_{i+1}
    +\frac{(1-q)x_i}{x_i-qx_{i+2}}z_{i+2},\\
  z_iz_{i+1}
  &\longmapsto
    \frac{(1-q)x_{i+2}}{x_i-qx_{i+2}}z_iz_{i+1}
    +\frac{(1-q)x_{i+2}(x_i-x_{i+2})}
      {(x_i-qx_{i+2})(x_{i+1}-qx_{i+2})}z_iz_{i+2}
    +\frac{(x_i-x_{i+2})(x_{i+1}-x_{i+2})}
      {(x_i-qx_{i+2})(x_{i+1}-qx_{i+2})}z_{i+1}z_{i+2},\\
  z_{i+1}z_{i+2}
  &\longmapsto
    \frac{q^2(x_i-x_{i+1})(x_i-x_{i+2})}
      {(x_i-qx_{i+1})(x_i-qx_{i+2})}z_iz_{i+1}
    +\frac{q(1-q)x_i(x_i-x_{i+2})}
      {(x_i-qx_{i+1})(x_i-qx_{i+2})}z_iz_{i+2}
    +\frac{(1-q)x_i}{x_i-qx_{i+2}}z_{i+1}z_{i+2}.
\end{align*}
Moreover, \eqref{eqAdjacentExchangeBlock} gives
$\mathscr R_i(X_n)(qz_i+z_{i+1})=qz_i+z_{i+1}$.
Thus both $\mathscr R_i$ and $\mathscr R_{i+1}$ fix
  $q^2z_i+qz_{i+1}+z_{i+2}$ and 
  $q^2z_iz_{i+1}+qz_iz_{i+2}+z_{i+1}z_{i+2}$.
The two products therefore agree on $z_{i+1}$ and
$z_iz_{i+2}$ as well.
They also fix $1$ and $z_iz_{i+1}z_{i+2}$.
This proves~\eqref{eqBExchangeRelations}.

For $|i-j|>1$, the operators act on disjoint pairs of
coordinates, so
\[
  \mathscr R_j(s_iX_n)\mathscr R_i(X_n)
  =
  \mathscr R_i(s_jX_n)\mathscr R_j(X_n).
\]
Together with~\eqref{eqBExchangeRelations}, these are the
defining relations for adjacent transpositions
\cite[Proposition~1.5.4]{BB05}.
Hence two sequences giving the same permutation give the
same operator, with the alphabet updated after each exchange.

By~\eqref{eqAdjacentExchangeBlock},
\begin{equation}\label{eqVTwoS7DeletionExchanges}
  \left.\mathscr R_i(X_n)f\right|_{z_i=z_{i+1}=0}
  =
  \left.f\right|_{z_i=z_{i+1}=0},
  \qquad
  \epsilon\mathscr R_i(X_n)f=\epsilon f.
\end{equation}
Keep each deleted letter at the right and postpone setting
its coordinate to zero.
Two deletion orders leave the undeleted letters in the same
order; their final alphabets differ only in the order of the
deleted letters.
The corresponding products of exchanges therefore differ
by exchanges among those final positions.
These disappear after evaluation at zero
by~\eqref{eqVTwoS7DeletionExchanges}.
Thus $\partial_I$ is independent of the deletion order.

For $\epsilon_I$, exchanges among the undeleted letters
also disappear by~\eqref{eqVTwoS7DeletionExchanges}.
Consequently,
\begin{equation}\label{eqBDeletionExchangeEvaluation}
  \epsilon_I^{X_n}f
  =
  \epsilon\partial_I(s_iX_n)\mathscr R_i(X_n)f
  =
  \epsilon_I^{s_iX_n}\mathscr R_i(X_n)f.
\end{equation}
We also obtain, for $h\in X_n\setminus I$,
\begin{equation}\label{eqBDeletionComposition}
  \epsilon_{\{h\}}^{X_n\setminus I}\partial_I f
  =
  \epsilon\partial_{\{h\}}\partial_I f
  =
  \epsilon\partial_{I\cup\{h\}}f
  =
  \epsilon_{I\cup\{h\}}^{X_n}f.
\end{equation}

\textup{(1)}
We first take $I=\varnothing$ and induct on $n$.
For $n=0$, both sides are zero.
Write $X_n=(X_{n-1},x)$.
For $h\in X_{n-1}$, we claim that
\begin{equation}\label{eqVTwoS7SingleDeletionInduction}
\begin{aligned}
  (h-qx)\epsilon_{\{h\}}^{X_n}f
  &=(h-x)\epsilon_{\{h\}}^{X_{n-1}}
    \bigl(qf|_{z_x=1}+(1-q)f|_{z_x=0}\bigr)
    +(1-q)x\epsilon(f|_{z_x=0}).
\end{aligned}
\end{equation}
Using~\eqref{eqBDeletionExchangeEvaluation}, move $h$ to
the end of $X_{n-1}$.
These exchanges commute with evaluation at $z_x=0,1$
and preserve $\epsilon$.
On $1,z_h,z_x,z_hz_x$, both sides of
\eqref{eqVTwoS7SingleDeletionInduction} equal
$h-qx$, $(1-q)x$, $q(h-x)$ and $0$, respectively,
by~\eqref{eqAdjacentExchangeBlock}.
This proves~\eqref{eqVTwoS7SingleDeletionInduction}.

We shall also use
$\epsilon_{\{x\}}^{X_n}f=\epsilon(f|_{z_x=0})$.
For $n>1$, we have
\begin{flalign}
  1-(1-q)\sum_{h\in X_{n-1}}
    \frac{h\eta_h(X_{n-1})}{h-x}
  &=1-(1-q)\sum_{h\in X_{n-1}}\eta_h(X_{n-1})
    +(1-q)x\sum_{h\in X_{n-1}}
      \frac{\eta_h(X_{n-1})}{x-h}
  &&\notag\\
  &=q^{n-1}
    +(1-q)x\sum_{h\in X_{n-1}}
      \frac{\eta_h(X_{n-1})}{x-h}
  &\mathllap{\text{(by~\eqref{eqEtaSubsetSum})}}
  &\notag\\
  &=\prod_{h\in X_{n-1}}\frac{x-qh}{x-h}
  &\mathllap{\text{(by~\eqref{eqEtaPartialFractions})}}
  &\notag\\
  &=\eta_x(X_n).
  &\mathllap{\text{(by~\eqref{eqRefinedS13Eta})}}
  &\label{eqBLastLetterEta}
\end{flalign}
For $n=1$, both sides of~\eqref{eqBLastLetterEta} are $1$.
Therefore,
\begin{flalign}
  \epsilon\mathcal D_{X_n}f
  &=\epsilon\mathcal D_{X_{n-1}}
    \bigl(qf|_{z_x=1}+(1-q)f|_{z_x=0}\bigr)
    +(q-1)x\epsilon(f|_{z_x=0})
  &\mathllap{\text{(by~\eqref{eqRefinedS6ABCD})}}
  &\notag\\
  &=(q-1)\sum_{h\in X_{n-1}}h\eta_h(X_{n-1})
    \epsilon_{\{h\}}^{X_{n-1}}
    \bigl(qf|_{z_x=1}+(1-q)f|_{z_x=0}\bigr)
    +(q-1)x\epsilon(f|_{z_x=0})
  &&\notag\\[-.5ex]
  &&\mathllap{\text{(by the induction hypothesis)}}
  &\notag\\
  &=(q-1)\sum_{h\in X_{n-1}}
    \frac{h\eta_h(X_{n-1})}{h-x}
    \Bigl((h-qx)\epsilon_{\{h\}}^{X_n}f
      -(1-q)x\epsilon(f|_{z_x=0})\Bigr)
    +(q-1)x\epsilon(f|_{z_x=0})
  &&\notag\\[-.5ex]
  &&\mathllap{\text{(by~\eqref{eqVTwoS7SingleDeletionInduction})}}
  &\notag\\
  &=(q-1)\sum_{h\in X_{n-1}}
    h\eta_h(X_n)\epsilon_{\{h\}}^{X_n}f
    +(q-1)x
    \left(
      1-(1-q)\sum_{h\in X_{n-1}}
        \frac{h\eta_h(X_{n-1})}{h-x}
    \right)\epsilon(f|_{z_x=0})
  &&\notag\\[-.5ex]
  &&\mathllap{\text{(by~\eqref{eqRefinedS13Eta})}}
  &\notag\\
  &=(q-1)\sum_{h\in X_n}
    h\eta_h(X_n)\epsilon_{\{h\}}^{X_n}f.
  &\mathllap{\text{(by~\eqref{eqBLastLetterEta})}}
  &\label{eqVTwoS7EmptyDeletionD}
\end{flalign}

For arbitrary $I\subseteq X_n$,
\begin{flalign*}
  \epsilon_I^{X_n}\mathcal Df
  &=\epsilon\mathcal D_{X_n\setminus I}\partial_I f
  &\mathllap{\text{(by~\eqref{eqBDeletionIntertwining})}}
  &\\
  &=(q-1)\sum_{h\in X_n\setminus I}
    h\eta_h(X_n\setminus I)
    \epsilon_{\{h\}}^{X_n\setminus I}\partial_I f
  &\mathllap{\text{(by~\eqref{eqVTwoS7EmptyDeletionD})}}
  &\\
  &=(q-1)\sum_{h\in X_n\setminus I}
    h\eta_h(X_n\setminus I)
    \epsilon_{I\cup\{h\}}^{X_n}f.
  &\mathllap{\text{(by~\eqref{eqBDeletionComposition})}}
  &
\end{flalign*}
This proves~\eqref{eqRefinedS7DeletionD}.

\textup{(2)}
First take $I=\varnothing$ and write $X_n=(x,X_{n-1})$.
By~\eqref{eqRefinedS6ABCD},
\begin{equation}\label{eqVTwoS7FirstCoordinateD}
  \epsilon\mathcal D_{X_n}f
  =
  \epsilon\mathcal D_{X_{n-1}}(f|_{z_x=1})
  +(q-1)x\epsilon\mathcal A(f|_{z_x=0}).
\end{equation}
For $h\in X_{n-1}$, deleting $h$ leaves $z_x$ unchanged, so
$\epsilon_{\{h\}}^{X_{n-1}}(f|_{z_x=1})
=\epsilon_{\{h\}}^{X_n}f$.
Hence
\begin{flalign*}
  (q-1)x\epsilon\mathcal A(f|_{z_x=0})
  &=\epsilon\mathcal D_{X_n}f
    -\epsilon\mathcal D_{X_{n-1}}(f|_{z_x=1})
  &\mathllap{\text{(by~\eqref{eqVTwoS7FirstCoordinateD})}}
  &\\
  &=(q-1)x\eta_x(X_n)\epsilon_{\{x\}}^{X_n}f
    +(q-1)\sum_{h\in X_{n-1}}
      h\bigl(\eta_h(X_n)-\eta_h(X_{n-1})\bigr)
      \epsilon_{\{h\}}^{X_n}f.
    &&\notag\\[-.5ex]
  &&\mathllap{\text{(by~\eqref{eqVTwoS7EmptyDeletionD})}}
  &
\end{flalign*}
By~\eqref{eqRefinedS13Eta},
\[
  \eta_h(X_n)-\eta_h(X_{n-1})
  =
  \left(\frac{h-qx}{h-x}-1\right)\eta_h(X_{n-1})
  =
  \frac{(1-q)x}{h-x}\eta_h(X_{n-1}).
\]
Dividing by $(q-1)x$ therefore gives
\begin{equation}\label{eqVTwoS7FirstCoordinateCalculation}
  \epsilon\mathcal A(f|_{z_x=0})
  =\prod_{h\in X_{n-1}}\frac{x-qh}{x-h}\,
      \epsilon_{\{x\}}^{X_n}f
    +(1-q)\sum_{h\in X_{n-1}}
      \frac{h\eta_h(X_{n-1})}{h-x}\,
      \epsilon_{\{h\}}^{X_n}f.
\end{equation}

Now let $I\subseteq X_{n-1}$.
For a basis monomial $z_A$ on $X_{n-1}$, every monomial
occurring in $\partial_I z_A$ has degree $|A|$.
Thus~\eqref{eqRefinedS6ABCD} gives
\[
  \partial_I\mathcal A z_A
  =q^{|A|}\partial_I z_A
  =\mathcal A\partial_I z_A.
\]
These deletions leave $z_x$ unchanged.
Consequently,
\begin{equation}\label{eqBDeletionACompatibility}
  \partial_I\mathcal A=\mathcal A\partial_I,
  \qquad
  \partial_I(f|_{z_x=0})
  =(\partial_I f)|_{z_x=0},
\end{equation}
where $\mathcal A$ acts on the remaining coordinates after
$\partial_I$.
Applying~\eqref{eqVTwoS7FirstCoordinateCalculation} on the
alphabet $(x,X_{n-1}\setminus I)$, we obtain
\begin{flalign*}
  \epsilon_I^{X_{n-1}}\mathcal A(f|_{z_x=0})
  &=\epsilon\mathcal A\bigl((\partial_I f)|_{z_x=0}\bigr)
  &\mathllap{\text{(by~\eqref{eqBDeletionACompatibility})}}
  &\\
  &=\prod_{h\in X_{n-1}\setminus I}
      \frac{x-qh}{x-h}\,
      \epsilon_{\{x\}}^{X_n\setminus I}\partial_I f
    +(1-q)\sum_{h\in X_{n-1}\setminus I}
      \frac{h\eta_h(X_{n-1}\setminus I)}{h-x}\,
      \epsilon_{\{h\}}^{X_n\setminus I}\partial_I f
        &&\notag\\[-.5ex]
  &&\mathllap{\text{(by~\eqref{eqVTwoS7FirstCoordinateCalculation})}}
  &\\
  &=\prod_{h\in X_{n-1}\setminus I}
      \frac{x-qh}{x-h}\,
      \epsilon_{I\cup\{x\}}^{X_n}f
    +(1-q)\sum_{h\in X_{n-1}\setminus I}
      \frac{h\eta_h(X_{n-1}\setminus I)}{h-x}\,
      \epsilon_{I\cup\{h\}}^{X_n}f.
  &&\mathllap{\text{(by~\eqref{eqBDeletionComposition})}}
  &
\end{flalign*}
This proves~\eqref{eqRefinedS7DeletionA}.
\end{proof}

\begin{proof}[Proof of Lemma~\ref{lemCyclicTransferVector}]
We shall prove that
\begin{equation}\label{eqRefinedS7CyclicCandidate}
  z_x\mathbf b_{(tx,X_{n-1})}|_{z_{tx}=1}
  +\frac{1-z_x}{(1+rx)(1+sx)}
   \Bigl(
     (1-rsx\mathcal D)\mathbf b_{(tx,X_{n-1})}|_{z_{tx}=1}
     +tx(rsx-u)\mathcal A\mathbf b_{(tx,X_{n-1})}|_{z_{tx}=0}
   \Bigr)
\end{equation}
satisfies~\eqref{eqNormalizedTransferVector}.
Here the parameters of $\mathbf b$ are $ut,r,s$.
The denominators have constant term $1$, so the coefficients
have formal expansions in $u,r,s$.
Setting all coordinates to $1$ gives
$\epsilon\mathbf b_{(tx,X_{n-1})}=1$.
It therefore remains to verify the recurrence.

We first establish
\begin{equation}\label{eqRiseMatrixCyclicRelations}
\begin{aligned}
  \mathcal B R_{X_{n-1}}(ut,r)
  &=(ut\mathcal A+r)R_{X_{n-1}}(ut,r)\mathcal D,\\
  \mathcal A R_{X_{n-1}}(ut,r)(1-ut\mathcal D)
  &=R_{X_{n-1}}(qut,qr)(1+r\mathcal D).
\end{aligned}
\end{equation}
For these calculations, $S$ belongs to $\mathcal Q_{n-1}$.
By~\eqref{eqRefinedS6OrderedProduct},
$Q_BS=\sum_{i\notin B}q^{|B\cap[1,i)|}Q_{B\cup\{i\}}$.
Thus, for $A,B\subseteq[n-1]$,
\begin{flalign*}
 (\mathcal B R_{X_{n-1}}(ut,r))_{B,A}
  &=\frac{g_{|B|+1}(ut,r)}{(q-1)^{|B|}}
    \operatorname{ev}_A\bigl(Q_BSG_r(S)^{-1}\bigr)
  &\mathllap{\text{(by~\eqref{eqRefinedS6RDefinition})}}&\\
  &=(utq^{|B|}+r)\frac{g_{|B|}(ut,r)}{(q-1)^{|B|}}
    \operatorname{ev}_A\bigl(Q_BG_r(S)^{-1}S\bigr)
  &\mathllap{\text{(by~\eqref{eqRefinedS6gDefinition})}}&\\
  &=\bigl((ut\mathcal A+r)
    R_{X_{n-1}}(ut,r)\mathcal D\bigr)_{B,A}.
  &\mathllap{\text{(by~\eqref{aboveeq})}}&
\end{flalign*}
Here $SG_r(S)^{-1}=G_r(S)^{-1}S$.

Proposition~\ref{propGIdentities}\textup{(2)} gives
\begin{equation}\label{eqBCyclicInverseShift}
  G_r(qS)^{-1}(1+rS)=(1-utS)G_r(S)^{-1}.
\end{equation}
Since $g_j(qut,qr)=q^jg_j(ut,r)$,
\eqref{eqRefinedS6RDefinition} and~\eqref{aboveeq} give
\begin{flalign*}
  \bigl(\mathcal A R_{X_{n-1}}(ut,r)
    (1-ut\mathcal D)\bigr)_{B,A}
  &=\frac{q^{|B|}g_{|B|}(ut,r)}{(q-1)^{|B|}}
    \operatorname{ev}_A
    \bigl(Q_BG_r(S)^{-1}(1-utS)\bigr)
  &&\\
  &=\frac{q^{|B|}g_{|B|}(ut,r)}{(q-1)^{|B|}}
    \operatorname{ev}_A
    \bigl(Q_BG_r(qS)^{-1}(1+rS)\bigr)
  &\mathllap{\text{(by~\eqref{eqBCyclicInverseShift})}}&\\
  &=\bigl(R_{X_{n-1}}(qut,qr)
    (1+r\mathcal D)\bigr)_{B,A}.
  &\mathllap{\text{(by~\eqref{eqRefinedS6RDefinition}
    and~\eqref{aboveeq})}}&
\end{flalign*}
This proves~\eqref{eqRiseMatrixCyclicRelations}.
Consequently,
\begin{flalign}
  R_{X_{n-1}}(qut,qr)
  &=\mathcal A R_{X_{n-1}}(ut,r)(1-ut\mathcal D)
    -rR_{X_{n-1}}(qut,qr)\mathcal D
  &\mathllap{\text{(by~\eqref{eqRiseMatrixCyclicRelations})}}
  &\notag\\
  &=\mathcal C R_{X_{n-1}}(ut,r)
    +r\bigl(R_{X_{n-1}}(ut,r)
      -R_{X_{n-1}}(qut,qr)\bigr)\mathcal D.
  &&\label{eqBCyclicCR}\\[-.5ex]
  &&\mathllap{\text{(by~\eqref{eqRiseMatrixCyclicRelations}
    and~\eqref{eqRefinedS6ABCD})}}
  &\notag
\end{flalign}

We now take $S=Q_1+\cdots+Q_n$.
By~\eqref{eqRefinedS6PowerSumQ}, the elements
$1,S,\ldots,S^n$ are linearly independent.
Comparing coefficients in~\eqref{eqBCyclicInverseShift} gives
\[
  (1-q^j)[S^j]G_r(S)^{-1}
  =(ut+rq^{j-1})[S^{j-1}]G_r(S)^{-1},
  \qquad 1\leq j\leq n.
\]
The constant term is $1$, so
\begin{equation}\label{eqBCyclicInverseSeries}
  G_r(S)^{-1}
  =\sum_{j=0}^{n}
    \frac{\prod_{h=0}^{j-1}(ut+rq^h)}{(q;q)_j}S^j.
\end{equation}
By~\eqref{eqRefinedS6OrderedProduct},
\begin{align*}
  S^j
  =(S-Q_n)^j
    +\sum_{h=0}^{j-1}
      (S-Q_n)^hQ_n(S-Q_n)^{j-1-h}
  =(S-Q_n)^j+[j]_qQ_n(S-Q_n)^{j-1},
\end{align*}
Using
\[
  [j+1]_q
  \frac{\prod_{h=0}^{j}(ut+rq^h)}{(q;q)_{j+1}}
  =
  -\frac{ut+rq^j}{q-1}
  \frac{\prod_{h=0}^{j-1}(ut+rq^h)}{(q;q)_j}
\]
in~\eqref{eqBCyclicInverseSeries}, we obtain
\begin{equation}\label{eqBCyclicInverseEnds}
\begin{aligned}
  G_r(S)^{-1}
  &=G_r(S-Q_n)^{-1}
    -\frac{Q_n}{q-1}
      \Bigl(utG_r(S-Q_n)^{-1}
        +rG_r(q(S-Q_n))^{-1}\Bigr),\\
  G_r(S)^{-1}
  &=G_r(S-Q_1)^{-1}
    -\frac{1}{q-1}
      \Bigl(utG_r(S-Q_1)^{-1}
        +rG_r(q(S-Q_1))^{-1}\Bigr)Q_1.
\end{aligned}
\end{equation}

For $(X_{n-1},x)$ and $A,I\subseteq[n-1]$,
\eqref{lemEvaluationOrderedMonomials} gives
\[
  \operatorname{ev}_A(Q_nQ_I)
    =(q-1)x\operatorname{ev}_A(Q_I),\qquad
  \operatorname{ev}_{A\cup\{n\}}(Q_I)
    =q^{|I|}\operatorname{ev}_A(Q_I),\qquad
  \operatorname{ev}_{A\cup\{n\}}(Q_nQ_I)=0.
\]
For $(tx,X_{n-1})$ and $A,I\subseteq\{2,\ldots,n\}$, it gives
\[
  \operatorname{ev}_A(Q_IQ_1)
    =(q-1)txq^{|A|}\operatorname{ev}_A(Q_I),\qquad
  \operatorname{ev}_{A\cup\{1\}}(Q_I)
    =\operatorname{ev}_A(Q_I),\qquad
  \operatorname{ev}_{A\cup\{1\}}(Q_IQ_1)=0.
\]
In the latter case, we identify the remaining indices
with those of $X_{n-1}$.
We also use
\[
  Q_CQ_n=q^{|C|}Q_nQ_C,\qquad
  g_{|C|}(qut,qr)=q^{|C|}g_{|C|}(ut,r),\qquad
  g_{|C|+1}(ut,r)
    =(ut+rq^{-|C|})g_{|C|}(qut,qr).
\]
Substitution in~\eqref{eqRefinedS6RDefinition}, together with
\eqref{eqBCyclicInverseEnds} and
$f=f|_{z_x=0}+z_x(f|_{z_x=1}-f|_{z_x=0})$, yields
\begin{equation}\label{eqBCyclicREnds}
\begin{aligned}
  (R_{(X_{n-1},x)}(ut,r)f)|_{z_x=0}
  &=
    \Bigl((1-utx\mathcal A)R_{X_{n-1}}(ut,r)
      -(1+rx)R_{X_{n-1}}(qut,qr)\Bigr)f|_{z_x=0}
    +R_{X_{n-1}}(qut,qr)f|_{z_x=1},\\
  (R_{(X_{n-1},x)}(ut,r)f)|_{z_x=1}
  &=
    (1+rx)\bigl(R_{X_{n-1}}(ut,r)
      -R_{X_{n-1}}(qut,qr)\bigr)f|_{z_x=0}
    +R_{X_{n-1}}(qut,qr)f|_{z_x=1},\\
  (R_{(tx,X_{n-1})}(ut,r)f)|_{z_{tx}=0}
  &=
    -tx\bigl(utR_{X_{n-1}}(ut,r)
      +r\mathcal A^{-1}R_{X_{n-1}}(qut,qr)\bigr)
      \mathcal A f|_{z_{tx}=0}
    +R_{X_{n-1}}(ut,r)f|_{z_{tx}=1},\\
  (R_{(tx,X_{n-1})}(ut,r)f)|_{z_{tx}=1}
  &=
    -ut\,tx\bigl(R_{X_{n-1}}(ut,r)
      -R_{X_{n-1}}(qut,qr)\bigr)
      \mathcal A f|_{z_{tx}=0}
    +R_{X_{n-1}}(ut,r)f|_{z_{tx}=1}.
\end{aligned}
\end{equation}

\eqref{eqRefinedS6VDefinition} gives
\[
  v_i(B,A\cup\{n\})
  =\left.v_i(B,A)\right|_{ut\mapsto qut}
  \quad(i<n,\ A,B\subseteq[n-1]),
\]
and
\[
  v_i(B\cup\{1\},A)
  =\left.v_i(B,A)\right|_{ut\mapsto qut}
  \quad(i>1,\ A,B\subseteq\{2,\ldots,n\}).
\]
Separating the last or first coordinate therefore gives
\begin{equation}\label{eqBCyclicVEnds}
\begin{aligned}
  (1+sx)(V_{(X_{n-1},x)}(ut,s)f)|_{z_x=0}
  &=(1-utx\mathcal A)V_{X_{n-1}}(ut,s)f|_{z_x=0}
    +V_{X_{n-1}}(qut,s)(f|_{z_x=1}-f|_{z_x=0}),\\
  (V_{(X_{n-1},x)}(ut,s)f)|_{z_x=1}
  &=V_{X_{n-1}}(ut,s)f|_{z_x=0}
    +V_{X_{n-1}}(qut,s)(f|_{z_x=1}-f|_{z_x=0}),\\
  (1+tsx)(V_{(tx,X_{n-1})}(ut,s)f)|_{z_{tx}=0}
  &=-ut\,txV_{X_{n-1}}(ut,s)\mathcal A f|_{z_{tx}=0}
    +V_{X_{n-1}}(ut,s)f|_{z_{tx}=1},\\
  (1+tsx)(V_{(tx,X_{n-1})}(ut,s)f)|_{z_{tx}=1}
  &=-ut\,tx\bigl(V_{X_{n-1}}(ut,s)
      -V_{X_{n-1}}(qut,s)\bigr)\mathcal A f|_{z_{tx}=0}\\
  &\quad
    +\bigl(V_{X_{n-1}}(ut,s)
      +tsxV_{X_{n-1}}(qut,s)\bigr)f|_{z_{tx}=1}.
\end{aligned}
\end{equation}

We next prove
\begin{equation}\label{eqBCyclicJEnds}
\begin{aligned}
  (\mathsf J_{(X_{n-1},x)}(s)f)|_{z_x=1}
  &=\mathsf J_{X_{n-1}}(s)f|_{z_x=1},\\
  (1+sx)(\mathsf J_{(X_{n-1},x)}(s)f)|_{z_x=0}
  &=(1+s\mathcal D)\mathsf J_{X_{n-1}}(s)f|_{z_x=0}
    -s\mathcal D\mathsf J_{X_{n-1}}(s)f|_{z_x=1},\\
  (1+tsx)(\mathsf J_{(tx,X_{n-1})}(s)f)|_{z_{tx}=0}
  &=\mathsf J_{X_{n-1}}(s)f|_{z_{tx}=0},\\
  (1+tsx)(\mathsf J_{(tx,X_{n-1})}(s)f)|_{z_{tx}=1}
  &=(1+tsx)\mathsf J_{X_{n-1}}(s)f|_{z_{tx}=1}
    -tsx\mathsf J_{X_{n-1}}(s)(1-\mathcal C)f|_{z_{tx}=0}.
\end{aligned}
\end{equation}
Compose on the right with $V_{(X_{n-1},x)}(ut,0)$ or
$V_{(tx,X_{n-1})}(ut,0)$.
By~\eqref{eqRefinedS6VFactorization}, it suffices to recover
\eqref{eqBCyclicVEnds}.
For $z_x=0$, we have
\begin{flalign*}
  &\quad (1+s\mathcal D)\mathsf J_{X_{n-1}}(s)
      \bigl((V_{(X_{n-1},x)}(ut,0)f)|_{z_x=0}\bigr)
    -s\mathcal D\mathsf J_{X_{n-1}}(s)
      \bigl((V_{(X_{n-1},x)}(ut,0)f)|_{z_x=1}\bigr)
  &&\\
  &=\Bigl(
      (1+s\mathcal D)\mathsf J_{X_{n-1}}(s)(1-utx\mathcal A)
      -s\mathcal D\mathsf J_{X_{n-1}}(s)
    \Bigr)V_{X_{n-1}}(ut,0)f|_{z_x=0}\\
  &\quad
    +\mathsf J_{X_{n-1}}(s)V_{X_{n-1}}(qut,0)
      (f|_{z_x=1}-f|_{z_x=0})
  &\mathllap{\text{(by~\eqref{eqBCyclicVEnds})}}&\\
  &=(1-utx\mathcal A)\mathsf J_{X_{n-1}}(s)
      V_{X_{n-1}}(ut,0)f|_{z_x=0}
    +\mathsf J_{X_{n-1}}(s)V_{X_{n-1}}(qut,0)
      (f|_{z_x=1}-f|_{z_x=0})
  &&\\[-.5ex]
  &&\mathllap{\text{(by Lemma~\ref{lemTransferMatrixRelations}\textup{(2)})}}&\\
  &=(1-utx\mathcal A)V_{X_{n-1}}(ut,s)f|_{z_x=0}
    +V_{X_{n-1}}(qut,s)(f|_{z_x=1}-f|_{z_x=0})
  &&\\[-.5ex]
  &&\mathllap{\text{(by~\eqref{eqRefinedS6VFactorization})}}&\\
  &=(1+sx)(V_{(X_{n-1},x)}(ut,s)f)|_{z_x=0}.
  &\mathllap{\text{(by~\eqref{eqBCyclicVEnds})}}&
\end{flalign*}
For $z_{tx}=1$, we obtain
\begin{flalign*}
  &\quad (1+tsx)\mathsf J_{X_{n-1}}(s)
      \bigl((V_{(tx,X_{n-1})}(ut,0)f)|_{z_{tx}=1}\bigr)
    -tsx\mathsf J_{X_{n-1}}(s)(1-\mathcal C)
      \bigl((V_{(tx,X_{n-1})}(ut,0)f)|_{z_{tx}=0}\bigr)
  &&\\
  &=-ut\,tx\mathsf J_{X_{n-1}}(s)
    \Bigl(
      (1+tsx)\bigl(V_{X_{n-1}}(ut,0)-V_{X_{n-1}}(qut,0)\bigr)
      -tsx(1-\mathcal C)V_{X_{n-1}}(ut,0)
    \Bigr)\mathcal A f|_{z_{tx}=0}\\
  &\quad
    +\mathsf J_{X_{n-1}}(s)
    \Bigl((1+tsx)V_{X_{n-1}}(ut,0)
      -tsx(1-\mathcal C)V_{X_{n-1}}(ut,0)\Bigr)f|_{z_{tx}=1}
  &&\\[-.5ex]
  &&\mathllap{\text{(by~\eqref{eqBCyclicVEnds})}}&\\
  &=-ut\,tx\mathsf J_{X_{n-1}}(s)
      \bigl(V_{X_{n-1}}(ut,0)-V_{X_{n-1}}(qut,0)\bigr)
      \mathcal A f|_{z_{tx}=0}\\
  &\quad
    +\mathsf J_{X_{n-1}}(s)
      \bigl(V_{X_{n-1}}(ut,0)+tsxV_{X_{n-1}}(qut,0)\bigr)
      f|_{z_{tx}=1}
  &\mathllap{\text{(by Lemma~\ref{lemTransferMatrixRelations}\textup{(1)})}}&\\
  &=-ut\,tx\bigl(V_{X_{n-1}}(ut,s)-V_{X_{n-1}}(qut,s)\bigr)
      \mathcal A f|_{z_{tx}=0}
    +\bigl(V_{X_{n-1}}(ut,s)+tsxV_{X_{n-1}}(qut,s)\bigr)
      f|_{z_{tx}=1}
  &&\\[-.5ex]
  &&\mathllap{\text{(by~\eqref{eqRefinedS6VFactorization})}}&\\
  &=(1+tsx)(V_{(tx,X_{n-1})}(ut,s)f)|_{z_{tx}=1}.
  &\mathllap{\text{(by~\eqref{eqBCyclicVEnds})}}&
\end{flalign*}
The equalities at $z_x=1$ and $z_{tx}=0$ follow directly
from~\eqref{eqBCyclicVEnds} and~\eqref{eqRefinedS6VFactorization}.
Since the two operators $V(ut,0)$ are invertible over the
rational function field, this proves~\eqref{eqBCyclicJEnds}.

By~\eqref{eqRefinedS6ABCD} and
\eqref{eqRiseMatrixCyclicRelations} at $(qut,qr)$,
\begin{equation}\label{eqBCyclicBShift}
  \mathcal B\mathcal A^{-1}R_{X_{n-1}}(qut,qr)
  =q^{-1}\mathcal A^{-1}\mathcal B R_{X_{n-1}}(qut,qr)
  =(ut+r\mathcal A^{-1})R_{X_{n-1}}(qut,qr)\mathcal D.
\end{equation}
Consequently,
\begin{flalign}
  &\quad (1-\mathcal C)
    \bigl(utR_{X_{n-1}}(ut,r)
      +r\mathcal A^{-1}R_{X_{n-1}}(qut,qr)\bigr)
  &&\notag\\
  &=utR_{X_{n-1}}(ut,r)
    -ut\mathcal A R_{X_{n-1}}(ut,r)(1-ut\mathcal D)
    +utrR_{X_{n-1}}(ut,r)\mathcal D\notag\\
  &\quad
    +r\mathcal A^{-1}R_{X_{n-1}}(qut,qr)
    -rR_{X_{n-1}}(qut,qr)
    +r(ut+r\mathcal A^{-1})R_{X_{n-1}}(qut,qr)\mathcal D
  &&\notag\\[-.5ex]
  &&\mathllap{\text{(by~\eqref{eqRefinedS6ABCD},
    \eqref{eqRiseMatrixCyclicRelations}, and~\eqref{eqBCyclicBShift})}}
  &\notag\\
  &=utR_{X_{n-1}}(ut,r)
    -(ut+r)R_{X_{n-1}}(qut,qr)
    +utrR_{X_{n-1}}(ut,r)\mathcal D
    +r\mathcal A^{-1}R_{X_{n-1}}(qut,qr)(1+r\mathcal D)
  &&\notag\\[-.5ex]
  &&\mathllap{\text{(by~\eqref{eqRiseMatrixCyclicRelations})}}
  &\notag\\
  &=utR_{X_{n-1}}(ut,r)
    -(ut+r)R_{X_{n-1}}(qut,qr)
    +utrR_{X_{n-1}}(ut,r)\mathcal D
    +rR_{X_{n-1}}(ut,r)(1-ut\mathcal D)
  &&\notag\\[-.5ex]
  &&\mathllap{\text{(by~\eqref{eqRiseMatrixCyclicRelations})}}
  &\notag\\
  &=(ut+r)\bigl(R_{X_{n-1}}(ut,r)
    -R_{X_{n-1}}(qut,qr)\bigr).
  &&\label{eqBCyclicAuxiliary}
\end{flalign}

Using $M_{X_n}=\mathsf J_{X_n}R_{X_n}$
from~\eqref{eqRefinedS6MDefinition}, substitute
\eqref{eqBCyclicREnds} into~\eqref{eqBCyclicJEnds}.
Then~\eqref{eqBCyclicAuxiliary} and
Lemma~\ref{lemTransferMatrixRelations}\textup{(2)} give
\begin{equation}\label{eqCyclicMEvaluations}
\begin{aligned}
  (M_{(X_{n-1},x)}(ut,r,s)f)|_{z_x=1}
  &=(1+rx)\mathsf J_{X_{n-1}}(s)
    \bigl(R_{X_{n-1}}(ut,r)-R_{X_{n-1}}(qut,qr)\bigr)
      f|_{z_x=0}\\
  &\quad
    +\mathsf J_{X_{n-1}}(s)R_{X_{n-1}}(qut,qr)f|_{z_x=1},\\
  (1+sx)(M_{(X_{n-1},x)}(ut,r,s)f)|_{z_x=0}
  &=\Bigl(
      (1-rsx\mathcal D-utx\mathcal A)
        \mathsf J_{X_{n-1}}(s)R_{X_{n-1}}(ut,r)
      -(1+rx)\mathsf J_{X_{n-1}}(s)R_{X_{n-1}}(qut,qr)
    \Bigr)f|_{z_x=0}\\
  &\quad
    +\mathsf J_{X_{n-1}}(s)R_{X_{n-1}}(qut,qr)f|_{z_x=1},\\
  (1+tsx)(M_{(tx,X_{n-1})}(ut,r,s)f)|_{z_{tx}=0}
  &=-tx\mathsf J_{X_{n-1}}(s)
    \bigl(utR_{X_{n-1}}(ut,r)
      +r\mathcal A^{-1}R_{X_{n-1}}(qut,qr)\bigr)
      \mathcal A f|_{z_{tx}=0}\\
  &\quad
    +\mathsf J_{X_{n-1}}(s)R_{X_{n-1}}(ut,r)f|_{z_{tx}=1},\\
  (1+tsx)(M_{(tx,X_{n-1})}(ut,r,s)f)|_{z_{tx}=1}
  &=tx(trsx-ut)\mathsf J_{X_{n-1}}(s)
    \bigl(R_{X_{n-1}}(ut,r)-R_{X_{n-1}}(qut,qr)\bigr)
      \mathcal A f|_{z_{tx}=0}\\
  &\quad
    +\mathsf J_{X_{n-1}}(s)(1+tsx\mathcal C)
      R_{X_{n-1}}(ut,r)f|_{z_{tx}=1}.
\end{aligned}
\end{equation}

In~\eqref{eqRefinedS7CyclicCandidate}, replace
$\mathbf b_{(tx,X_{n-1})}$ by
$M_{(tx,X_{n-1})}(ut,r,s)f$.
Compare this with $M_{X_n}(ut,r,s)$ applied to
\eqref{eqRefinedS7CyclicCandidate} after replacing
$u,s,\mathbf b_{(tx,X_{n-1})}$ by $ut,ts,f$, respectively.
Both expressions are linear in $z_x$ and in
$f|_{z_{tx}=0},f|_{z_{tx}=1}$.
We compare the coefficients of these two evaluations
at $z_x=0,1$, using~\eqref{eqCyclicMEvaluations}.

At $z_x=1$, multiply by $1+tsx$.
The coefficients of $f|_{z_{tx}=0}$ agree because
$tx(trsx-ut)=t^2x(rsx-u)$.
For $f|_{z_{tx}=1}$, we have
\begin{flalign*}
  &\quad \bigl(R_{X_{n-1}}(ut,r)-R_{X_{n-1}}(qut,qr)\bigr)
      (1-trsx\mathcal D)
    +(1+tsx)R_{X_{n-1}}(qut,qr)
  &&\\
  &=R_{X_{n-1}}(ut,r)
    +tsx\Bigl(
      R_{X_{n-1}}(qut,qr)
      -r\bigl(R_{X_{n-1}}(ut,r)
        -R_{X_{n-1}}(qut,qr)\bigr)\mathcal D
    \Bigr)
  &&\\
  &=(1+tsx\mathcal C)R_{X_{n-1}}(ut,r).
  &\mathllap{\text{(by~\eqref{eqBCyclicCR})}}&
\end{flalign*}
Multiplying by $\mathsf J_{X_{n-1}}(s)$ on the left gives
the required equality.

At $z_x=0$, multiply by $(1+rx)(1+sx)(1+tsx)$.
Lemma~\ref{lemTransferMatrixRelations}(2) gives 
\begin{equation}\label{eqRefinedS7CyclicJ}
  (1-rsx\mathcal D)\mathsf J_{X_{n-1}}(s)
    +rx\mathcal A\mathsf J_{X_{n-1}}(s)\mathcal A^{-1}
  =(1-rsx\mathcal D)\mathsf J_{X_{n-1}}(s)
    +rx(1+s\mathcal D)\mathsf J_{X_{n-1}}(s)
  =(1+rx)\mathsf J_{X_{n-1}}(s).
\end{equation}
For the coefficients of $f|_{z_{tx}=0}$, calculate
\begin{flalign*}
  &\quad (1-rsx\mathcal D)\mathsf J_{X_{n-1}}(s)
    \bigl(R_{X_{n-1}}(ut,r)-R_{X_{n-1}}(qut,qr)\bigr)
    -x\mathcal A\mathsf J_{X_{n-1}}(s)
      \bigl(utR_{X_{n-1}}(ut,r)
        +r\mathcal A^{-1}R_{X_{n-1}}(qut,qr)\bigr)
  &&\\
  &=(1-rsx\mathcal D-utx\mathcal A)
      \mathsf J_{X_{n-1}}(s)R_{X_{n-1}}(ut,r)
    -\Bigl(
      (1-rsx\mathcal D)\mathsf J_{X_{n-1}}(s)
      +rx\mathcal A\mathsf J_{X_{n-1}}(s)\mathcal A^{-1}
    \Bigr)R_{X_{n-1}}(qut,qr)
  &&\\
  &=(1-rsx\mathcal D-utx\mathcal A)
      \mathsf J_{X_{n-1}}(s)R_{X_{n-1}}(ut,r)
    -(1+rx)\mathsf J_{X_{n-1}}(s)R_{X_{n-1}}(qut,qr).
  &&\mathllap{\text{(by~\eqref{eqRefinedS7CyclicJ})}}&
\end{flalign*}
Multiplying this equality by $t^2x(rsx-u)$ and by
$\mathcal A$ on the right proves that these coefficients agree.

For the coefficients of $f|_{z_{tx}=1}$,
right multiplication of~\eqref{eqRefinedS7CyclicJ}
by $\mathcal A$ gives
\[
  (1+rx)\mathsf J_{X_{n-1}}(s)\mathcal A
  =
  (1-rsx\mathcal D)\mathsf J_{X_{n-1}}(s)\mathcal A
  +rx\mathcal A\mathsf J_{X_{n-1}}(s).
\]
Substitution yields
\begin{equation}\label{eqBCyclicTailRelations}
\begin{aligned}
  &r(1-rsx\mathcal D-utx\mathcal A)\mathsf J_{X_{n-1}}(s)
    +ut(1+rx)\mathsf J_{X_{n-1}}(s)\mathcal A
    =(1-rsx\mathcal D)\mathsf J_{X_{n-1}}(s)(r+ut\mathcal A),\\
  &(1-rsx\mathcal D-utx\mathcal A)\mathsf J_{X_{n-1}}(s)
    +tsx(1+rx)\mathsf J_{X_{n-1}}(s)\mathcal A
    =tx(rsx-u)\mathcal A\mathsf J_{X_{n-1}}(s)
      +(1-rsx\mathcal D)\mathsf J_{X_{n-1}}(s)
        (1+tsx\mathcal A).
\end{aligned}
\end{equation}
Therefore,
\begin{flalign*}
  &\quad \Bigl(
    (1-rsx\mathcal D-utx\mathcal A)
      \mathsf J_{X_{n-1}}(s)R_{X_{n-1}}(ut,r)
    -(1+rx)\mathsf J_{X_{n-1}}(s)R_{X_{n-1}}(qut,qr)
    \Bigr)(1-trsx\mathcal D)\\
  &\quad \quad
    +(1+rx)(1+tsx)\mathsf J_{X_{n-1}}(s)
      R_{X_{n-1}}(qut,qr)
  &&\\
  &=(1-rsx\mathcal D-utx\mathcal A)
      \mathsf J_{X_{n-1}}(s)R_{X_{n-1}}(ut,r)
      (1-trsx\mathcal D)
    +tsx(1+rx)\mathsf J_{X_{n-1}}(s)
      \mathcal A R_{X_{n-1}}(ut,r)(1-ut\mathcal D)
  &\mathllap{\text{(by~\eqref{eqRiseMatrixCyclicRelations})}}&\\
  &=tx(rsx-u)\mathcal A\mathsf J_{X_{n-1}}(s)
      R_{X_{n-1}}(ut,r)
    +(1-rsx\mathcal D)\mathsf J_{X_{n-1}}(s)
      \Bigl(
        (1+tsx\mathcal A)R_{X_{n-1}}(ut,r)
        -tsx(r+ut\mathcal A)R_{X_{n-1}}(ut,r)\mathcal D
      \Bigr)
  &&\notag\\[-.5ex]
  &&\mathllap{\text{(by~\eqref{eqBCyclicTailRelations})}}&\\
  &=tx(rsx-u)\mathcal A\mathsf J_{X_{n-1}}(s)
      R_{X_{n-1}}(ut,r)
    +(1-rsx\mathcal D)\mathsf J_{X_{n-1}}(s)
      (1+tsx\mathcal A-tsx\mathcal B)R_{X_{n-1}}(ut,r)
  &&\mathllap{\text{(by~\eqref{eqRiseMatrixCyclicRelations})}}&\\
  &=tx(rsx-u)\mathcal A\mathsf J_{X_{n-1}}(s)
      R_{X_{n-1}}(ut,r)
    +(1-rsx\mathcal D)\mathsf J_{X_{n-1}}(s)
      (1+tsx\mathcal C)R_{X_{n-1}}(ut,r).
  &&\mathllap{\text{(by~\eqref{eqRefinedS6ABCD})}}&
\end{flalign*}
By~\eqref{eqCyclicMEvaluations}, the last expression is also
the coefficient of $f|_{z_{tx}=1}$ obtained from
\eqref{eqRefinedS7CyclicCandidate}.

Taking $f=\mathbf b_{(tx,X_{n-1})}(ut^2,r,ts)$
and using~\eqref{eqNormalizedTransferVector}, we conclude that
\eqref{eqRefinedS7CyclicCandidate} satisfies the required recurrence.
Its $\epsilon$-evaluation is $1$, so
Proposition~\ref{propNormalizedTransferVector} shows that it equals
$\mathbf b_{X_n}$.
Setting $z_x=0,1$ gives~\eqref{eqCyclicTransferVector}.
\end{proof}

\section*{Use of AI tools}
The theoretical derivations and the proofs in this paper were developed with the help of ChatGPT 6 Astra.  
The authors studied and checked all of this material, and wrote the paper themselves.  
We also checked the main identities by computer for small $n$.  
The authors take full responsibility for the content of the paper.

\providecommand{\bysame}{\leavevmode\hbox to3em{\hrulefill}\thinspace}

\end{document}